\documentclass[a4paper,12pt]{article} 
\usepackage{amssymb}
\usepackage{amscd}
\usepackage{amsmath} 
\usepackage{graphicx}
\usepackage{latexsym} 
\usepackage{enumerate}

\usepackage{theorem}
\theoremstyle{plain}
\theorembodyfont{\itshape}
\newtheorem{theorem}{Theorem}[section]
\newtheorem{proposition}[theorem]{Proposition}
\newtheorem{lemma}[theorem]{Lemma}
\newtheorem{corollary}[theorem]{Corollary}
 
\theorembodyfont{\rmfamily}
\newtheorem{assumption}[theorem]{Assumption}
\newtheorem{definition}[theorem]{Definition}
\newtheorem{remark}[theorem]{Remark}

\newtheorem{example}[theorem]{Example}

\def\C{\mathbb{C}}
\def\N{\mathbb{N}}

\def\R{\mathbb{R}}
\def\Z{\mathbb{Z}}
\def\P{\mathbb{P}}
\def\K{\mathcal{K}}
\def\M{\mathcal{M}}
\def\E{\mathcal{E}}
\def\a{\alpha}
\def\b{\beta}

\def\p{\mathfrak{p}}

\def\si{\sigma}
\def\ga{\gamma}

\def\ol{\overline}

\def\wt{\widetilde}
\def\PVI{\mathrm{P}_{\mathrm{VI}}}
\def\RH{\mathrm{RH}} 
\def\rh{\mathrm{rh}} 
\def\Tr{\mathrm{Tr}}
 
\def\Sol{\mathcal{S}}

\def\la{\langle}
\def\ra{\rangle}
\def\Per{\mathrm{Per}} 
\def\a{\alpha}
\def\b{\beta}
\def\ga{\gamma}
\def\k{\kappa}
\def\l{\lambda}
\def\si{\sigma}
\def\th{\theta}
\def\vD{\varDelta}

\def\Th{\Theta}
\def\ve{\varepsilon}
\def\car{\curvearrowright}
\def\mi{\phantom{-}}
\def\ds{\displaystyle} 
\def\carl{\circlearrowleft}

\title{\bf Singular Cubic Surfaces and \\[2mm]
the Dynamics of Painlev\'e V\!I} 
\author{Katsunori Iwasaki and Takato Uehara\thanks{E-mail addresses: 
{\tt iwasaki@math.kyushu-u.ac.jp} \ and \ 
{\tt t-uehara@math.kyushu-u.ac.jp}} \\ \\
Graduate School of Mathematics, Kyushu University \\
744 Moto-oka, Nishi-ku, Fukuoka 819-0395 Japan} 
\begin{document}
\maketitle
\begin{abstract} 
We develop a dynamical study of the sixth Painlev\'e equation 
for {\sl all} parameters generalizing an earlier work for 
generic parameters. 
Here the main focus of this paper is on non-generic parameters, 
for which the corresponding character variety becomes a cubic 
surface with simple singularities and the Riemann-Hilbert 
correspondence is a minimal resolution of the singular surface, 
not a biholomorphism as in the generic case. 
Introducing a suitable stratification on the parameter space and 
based on geometry of singular cubic surfaces, we establish a 
chaotic nature of the nonlinear monodromy map of Painlev\'e VI and 
give a precise estimate for the number of its isolated 
periodic solutions. 
\end{abstract} 
\section{Introduction} \label{sec:intro}
The aim of this paper is to develop a dynamical study of the sixth 
Painlev\'e equation for {\sl all} parameters generalizing an earlier 
work for generic parameters, where the main two issues are 
establishing a chaotic nature of the nonlinear monodromy map of 
Painlev\'e VI and counting the number of its periodic solutions 
along a given loop. 
For generic parameters these problems were discussed in a 
previous paper \cite{IU1}, so that the main focus of this paper 
is on non-generic parameters. 
A difficulty in the non-generic case stems from the fact that 
the corresponding character variety becomes an affine cubic surface 
with simple singularities and the Riemann-Hilbert correspondence is 
only an analytic minimal resolution of the singular surface, while 
it is a biholomorphism onto a smooth cubic surface in the generic 
case. 
Therefore, in order to apply a general dynamical systems theory 
on a smooth compact complex surface to the present situation, 
one has not only to compactify the affine surface but also to 
take a minimal resolution of the singular surface. 
Another difficulty occurs in counting the number of periodic 
solutions. 
For non-generic parameters there may be periodic solutions 
parametrized by a curve, in which case the counting problem 
obviously fails to make sense in its na\"{i}ve formulation. 
\par
Thus the main problem of this paper is to overcome these 
difficulties. 
As in \cite{IU1}, our principal tool is again a 
Riemann-Hilbert correspondence, but this time being a 
{\sl lifted} one onto a desingularized character variety, 
namely, a desingularized affine cubic surface. 
Through it the monodromy map of Painlev\'e VI is strictly 
conjugate to a biregular map on the latter surface, which 
in turn extends to a birational map on its compactification 
obtained by adding tritangent lines at infinity. 
A close investigation into the last map, especially, 
into its dynamical behaviors around the tritangent lines 
as well as around the exceptional set enables us to establish 
its ergodic properties based on the general theory of 
bimeromorphic surface maps developed in 
\cite{BD, DF, DS, Dujardin} (see also \cite{Cantat, CL}). 
As for counting the number of periodic solutions, the 
difficulty mentioned above is surmounted by our general 
theory of periodic points for area-preserving birational 
maps of surfaces in \cite{IU2}, in which a method of 
counting {\sl isolated} periodic points is developed 
in the presence of periodic curves. 
For all these discussions, a suitable stratification 
is introduced on the parameter space of Painlev\'e 
VI and several arguments are made stratum by stratum, 
becuase the singularities of cubic surfaces depend 
efficiently on the stratification. 
\par
The sixth Painlev\'e equation $\PVI(\k)$ is a Hamiltonian system 
of differential equations with an independent variable 
$z \in Z:= \P^1 \setminus \{ 0,1, \infty \}$ and unknown functions 
$(q,p)=(q(z),p(z))$: 
\[
\dfrac{d q}{d z} = \dfrac{\partial H(\k)}{\partial p}, 
\qquad 
\dfrac{d p}{d z} = -\dfrac{\partial H(\k)}{\partial q}, 
\]
where $\k:=(\k_0,\k_1,\k_2,\k_3,\k_4)$ are complex parameters 
lying in the parameter space
\[
\K := \{\, \k = (\k_0,\k_1,\k_2,\k_3,\k_4) \in \C^5 \,:\, 
2 \k_0 + \k_1 + \k_2 + \k_3 +\k_4 = 1 \,\}
\]
and the Hamiltonian $H(\k) = H(q,p,z;\k)$ is given by 
\[
\begin{array}{rcl}
z(z-1) H(\k) &=&  
(q_0q_1q_z) p^2 - \{\k_1q_1q_z 
+ (\k_2-1)q_0q_1 + \k_3q_0q_z \} p 
+ \k_0(\k_0+\k_4) q_z,
\end{array}
\]
with $q_{\nu} := q - \nu$ for $\nu \in \{0,1,z\}$. 
Let $\M_z(\k)$ be the set of all meromorphic solution germs to 
$\PVI(\k)$ at a point $z \in Z$. In \cite{IIS1,IIS2,IIS3}, 
the set $\M_z(\k)$ is realized as a moduli space of stable parabolic 
connections. Here a {\sl stable parabolic connection} is a rank 
$2$ vector bundle with a Fuchsian connection and a parabolic structure, 
having a Riemann scheme as in Table \ref{tab:riemann} and 
satisfying a sort of stability condition in geometric invariant theory. 
The parameter $\kappa_i$ is the difference of the second 
exponent from the first one at the regular singular point $t_i$ 
and thus $\lambda_i$ is uniquely determined from $\kappa_i$. 
\begin{table}[t]
\begin{center} 
\begin{tabular}{|c||c|c|c|c|}
\hline 
\vspace{-0.42cm} & ~ & ~ & ~ &  \\
singularities & $t_1 = 0 $ & $t_2 = z $ & $t_3 = 1 $ & 
$t_4 = \infty $ \\[1mm]
\hline 
\vspace{-0.42cm} & ~ & ~ & ~ &  \\
first exponent & $-\l_1$ & $-\l_2$ & $-\l_3$ & $-\l_4$ \\[1mm]
\hline 
\vspace{-0.42cm} & ~ & ~ & ~ &  \\
second exponent & $\l_1$ & $\l_2$  & $\l_3$ & $\l_4-1$ \\[1mm]
\hline 
\vspace{-0.42cm} & ~ & ~ & ~ &  \\
difference  & $\k_1$ & $\k_2$ & $\k_3$ & $\k_4$ \\[1mm]  
\hline
\end{tabular}
\end{center}
\caption{Riemann scheme: $\k_i$ is the difference of the second 
exponent from the first.} 
\label{tab:riemann}
\end{table}
This formulation provides the moduli space $\M_z(\k)$ with the 
structure of a smooth quasi-projective rational surface. 
It is known that there exists a natural compactification 
$\overline{\M}_z(\k)$ of the moduli space $\M_z(\k)$. 
The space $\overline{\M}_z(\k)$ has a unique effective anti-canonical 
divisor $\mathcal{Y}_z(\k)$ and $\M_z(\k)$ is obtained from 
$\overline{\M}_z(\k)$ by removing $\mathcal{Y}_z(\k)$. 
Thus there exists a global holomorphic $2$-form 
$\omega_z(\k)$ on $\M_z(\k)$, meromorphic on 
$\overline{\M}_z(\k)$ with pole divisor $\mathcal{Y}_z(\k)$. 
It is unique up to constant multiples and yields 
a natural area form on $\M_z(\k)$. 
\par
The Painlev\'e equation enjoys the Painlev\'e property, namely, 
any solution germ $Q \in \M_z(\k)$ can be continued analytically 
along any loop $\ga \in \pi_1(Z,z)$, so that the automorphism 
\[ 
\ga_* : \M_z(\k) \overset{\sim}{\to} \M_z(\k), \qquad 
Q \mapsto \ga_*Q 
\]
is well defined, where $\ga_*Q$ is the result of analytic 
continuation of $Q$ along $\ga$. 
The map $\ga_*$ is called the {\sl nonlinear monodromy map} 
along $\ga$. 
It preserves the area form $\omega_z(\k)$. 
The dynamical system of this map is what we want to study 
in this paper. 
However this map is too transcendental to deal with directly, 
so that it will be converted to a more tractable map on an 
affine cubic surface $\Sol(\th)$ via the Riemann-Hilbert 
correspondence 
\begin{equation} \label{eqn:RH2}
\RH_{z,\k} : \M_z(\k) \to \Sol(\th), 
\end{equation}
which is an analytic minimal resolution of the (possibly) 
singular surface $\Sol(\th)$ (see Theorem \ref{thm:SolRHP}). 
Through it the monodromy map $\ga_*$ is conjugated to a 
polynomial automorphism on $\Sol(\th)$. 
This last map was studied in \cite{Cantat,CL,Iwasaki1,IU1}, 
but a further investigation is made in this paper. 
\par
In terms of the dynamical behavior of the monodromy map, 
each loop in $\pi_1(Z,z)$ falls into two types, that is, 
an elementary loop and a non-elementary loop. 
\begin{definition} \label{def:elementary} 
Let $\ga_1$, $\ga_2$ and $\ga_3$ be loops in $\pi_1(Z,z)$ 
surrounding $0$, $1$ and $\infty$ respectively 
once anti-clockwise as in Figure \ref{fig:elementary}. 
\begin{figure}[t] 
\begin{center}
\unitlength 0.1in
\begin{picture}( 20.8600, 20.8600)(  3.6000,-27.1600)
%
\special{pn 13}%
\special{ar 1404 1674 1044 1044  0.0000000 6.2831853}%
%
\special{pn 20}%
\special{ar 1400 1178 228 228  0.0000000 6.2831853}%
%
\special{pn 20}%
\special{sh 1.000}%
\special{ar 1404 1678 34 34  0.0000000 6.2831853}%
%
\special{pn 20}%
\special{pa 1400 1404}%
\special{pa 1400 1648}%
\special{fp}%
%
\special{pn 20}%
\special{ar 1012 2008 228 228  0.0000000 6.2831853}%
%
\special{pn 20}%
\special{pa 1186 1860}%
\special{pa 1372 1702}%
\special{fp}%
%
\special{pn 20}%
\special{ar 1798 2004 228 228  0.0000000 6.2831853}%
%
\special{pn 20}%
\special{pa 1624 1858}%
\special{pa 1438 1702}%
\special{fp}%
%
\special{pn 20}%
\special{sh 1.000}%
\special{ar 1396 1174 34 34  0.0000000 6.2831853}%
%
\special{pn 20}%
\special{sh 1.000}%
\special{ar 1012 2010 34 34  0.0000000 6.2831853}%
%
\special{pn 20}%
\special{sh 1.000}%
\special{ar 1794 2002 34 34  0.0000000 6.2831853}%
\put(13.7000,-13.5000){\makebox(0,0)[lb]{$0$}}%
\put(9.8000,-22.0000){\makebox(0,0)[lb]{$1$}}%
\put(17.4000,-21.5000){\makebox(0,0)[lb]{$\infty$}}%
\put(13.6000,-18.5000){\makebox(0,0)[lb]{$z$}}%
\put(13.4300,-25.3300){\makebox(0,0)[lb]{$Z$}}%
%
\special{pn 13}%
\special{ar 1398 1176 268 268  4.3714100 5.1359243}%
%
\special{pn 13}%
\special{pa 1318 920}%
\special{pa 1288 938}%
\special{fp}%
\special{sh 1}%
\special{pa 1288 938}%
\special{pa 1356 922}%
\special{pa 1334 912}%
\special{pa 1336 888}%
\special{pa 1288 938}%
\special{fp}%
%
\special{pn 13}%
\special{ar 1008 2014 266 266  2.0157068 2.7806232}%
%
\special{pn 13}%
\special{pa 884 2252}%
\special{pa 916 2260}%
\special{fp}%
\special{sh 1}%
\special{pa 916 2260}%
\special{pa 858 2224}%
\special{pa 866 2246}%
\special{pa 846 2262}%
\special{pa 916 2260}%
\special{fp}%
%
\special{pn 13}%
\special{ar 1804 2000 266 266  0.2772553 1.0381186}%
%
\special{pn 13}%
\special{pa 2060 2082}%
\special{pa 2062 2050}%
\special{fp}%
\special{sh 1}%
\special{pa 2062 2050}%
\special{pa 2038 2116}%
\special{pa 2060 2104}%
\special{pa 2078 2118}%
\special{pa 2062 2050}%
\special{fp}%
\put(9.1000,-24.3000){\makebox(0,0)[lb]{$\ga_2$}}%
\put(17.2000,-24.0000){\makebox(0,0)[lb]{$\ga_3$}}%
\put(16.7000,-12.1000){\makebox(0,0)[lb]{$\ga_1$}}%
\end{picture}%
\end{center}
\caption{Three basic loops in $Z := \P^1 \setminus \{0,1,\infty\}$} 
\label{fig:elementary} 
\end{figure}
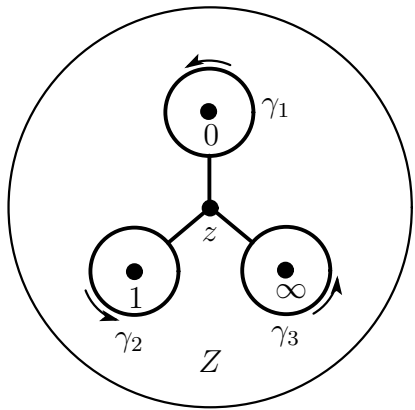
A loop $\ga \in \pi_1(Z,z)$ is said to be {\sl elementary} if 
$\ga$ is conjugate to the loop $\ga_i^m$ for some index 
$i \in \{1,2,3\}$ and some integer $m \in \Z$. 
Otherwise, $\ga$ is said to be {\sl non-elementary}. 
\end{definition} 
\par
If $\ga$ is elementary, the map $\ga_* : \M_z(\k) \carl$ 
preserves a fibration and exhibits a very simple dynamical 
behavior (see Remark \ref{rem:elem}). 
So from now on we assume that $\ga$ is non-elementary. 
\begin{remark} \label{rem:FDD} 
In \cite{IU1} we introduced an algebraic number 
$\lambda(\ga) \ge 1$ called {\sl the first dynamical degree} 
of $\ga \in \pi_1(Z,z)$ (see also Definition \ref{def:birat}) 
and established an algorithm to calculate $\lambda(\ga)$ in terms 
of a reduced word for the loop $\ga$ in the alphabet $\ga_1^{\pm1}$, 
$\ga_2^{\pm1}$, $\ga_3^{\pm1}$ (\cite[Theorem 3]{IU1}). 
It was also shown that $\lambda(\ga)$ is a quadratic unit strictly 
greater than $1$ if and only if $\ga$ is non-elementary. 
\end{remark}
\begin{example} \label{ex:loop}
Let $\ve$ and $\wp$ denote loops conjugate to $\ga_i \ga_j^{-1}$ and 
$[\ga_i, \ga_j^{-1}] = \ga_i \ga_j^{-1} \ga_i^{-1} \ga_j$ 
for some indices $\{i,j,k\} = \{1,2,3\}$ as in Figure \ref{fig:loop}. 
\begin{figure}[t] 
\begin{center}
\unitlength 0.1in
\begin{picture}( 56.5200, 11.0400)(  1.5000,-18.4500)
%
\special{pn 20}%
\special{sh 1.000}%
\special{ar 3822 1294 40 40  0.0000000 6.2831853}%
%
\special{pn 20}%
\special{sh 1.000}%
\special{ar 5256 1294 42 42  0.0000000 6.2831853}%
\put(37.7400,-14.7300){\makebox(0,0)[lb]{$z_i$}}%
\put(52.2000,-14.9000){\makebox(0,0)[lb]{$z_j$}}%
%
\special{pn 20}%
\special{pa 4176 1354}%
\special{pa 4734 1354}%
\special{fp}%
%
\special{pn 20}%
\special{pa 3840 934}%
\special{pa 3792 934}%
\special{fp}%
\special{sh 1}%
\special{pa 3792 934}%
\special{pa 3858 954}%
\special{pa 3844 934}%
\special{pa 3858 914}%
\special{pa 3792 934}%
\special{fp}%
%
\special{pn 20}%
\special{pa 3462 1270}%
\special{pa 3462 1322}%
\special{fp}%
\special{sh 1}%
\special{pa 3462 1322}%
\special{pa 3482 1256}%
\special{pa 3462 1270}%
\special{pa 3442 1256}%
\special{pa 3462 1322}%
\special{fp}%
%
\special{pn 20}%
\special{pa 3798 1658}%
\special{pa 3852 1658}%
\special{fp}%
\special{sh 1}%
\special{pa 3852 1658}%
\special{pa 3784 1638}%
\special{pa 3798 1658}%
\special{pa 3784 1678}%
\special{pa 3852 1658}%
\special{fp}%
%
\special{pn 20}%
\special{pa 4472 1354}%
\special{pa 4508 1354}%
\special{fp}%
\special{sh 1}%
\special{pa 4508 1354}%
\special{pa 4442 1334}%
\special{pa 4456 1354}%
\special{pa 4442 1374}%
\special{pa 4508 1354}%
\special{fp}%
%
\special{pn 20}%
\special{pa 5232 1666}%
\special{pa 5298 1658}%
\special{fp}%
\special{sh 1}%
\special{pa 5298 1658}%
\special{pa 5230 1646}%
\special{pa 5246 1664}%
\special{pa 5234 1686}%
\special{pa 5298 1658}%
\special{fp}%
%
\special{pn 20}%
\special{pa 5634 1318}%
\special{pa 5634 1270}%
\special{fp}%
\special{sh 1}%
\special{pa 5634 1270}%
\special{pa 5614 1336}%
\special{pa 5634 1322}%
\special{pa 5654 1336}%
\special{pa 5634 1270}%
\special{fp}%
%
\special{pn 20}%
\special{pa 5286 926}%
\special{pa 5250 926}%
\special{fp}%
\special{sh 1}%
\special{pa 5250 926}%
\special{pa 5318 946}%
\special{pa 5304 926}%
\special{pa 5318 906}%
\special{pa 5250 926}%
\special{fp}%
%
\special{pn 20}%
\special{pa 3846 1846}%
\special{pa 3804 1838}%
\special{fp}%
\special{sh 1}%
\special{pa 3804 1838}%
\special{pa 3866 1868}%
\special{pa 3856 1848}%
\special{pa 3872 1830}%
\special{pa 3804 1838}%
\special{fp}%
%
\special{pn 20}%
\special{pa 3276 1330}%
\special{pa 3276 1270}%
\special{fp}%
\special{sh 1}%
\special{pa 3276 1270}%
\special{pa 3256 1336}%
\special{pa 3276 1322}%
\special{pa 3296 1336}%
\special{pa 3276 1270}%
\special{fp}%
%
\special{pn 20}%
\special{pa 3804 746}%
\special{pa 3852 746}%
\special{fp}%
\special{sh 1}%
\special{pa 3852 746}%
\special{pa 3784 726}%
\special{pa 3798 746}%
\special{pa 3784 766}%
\special{pa 3852 746}%
\special{fp}%
%
\special{pn 20}%
\special{pa 5238 766}%
\special{pa 5304 758}%
\special{fp}%
\special{sh 1}%
\special{pa 5304 758}%
\special{pa 5236 746}%
\special{pa 5250 764}%
\special{pa 5240 786}%
\special{pa 5304 758}%
\special{fp}%
%
\special{pn 20}%
\special{pa 5796 1262}%
\special{pa 5802 1330}%
\special{fp}%
\special{sh 1}%
\special{pa 5802 1330}%
\special{pa 5816 1262}%
\special{pa 5796 1276}%
\special{pa 5776 1266}%
\special{pa 5802 1330}%
\special{fp}%
%
\special{pn 20}%
\special{pa 5286 1834}%
\special{pa 5232 1834}%
\special{fp}%
\special{sh 1}%
\special{pa 5232 1834}%
\special{pa 5298 1854}%
\special{pa 5284 1834}%
\special{pa 5298 1814}%
\special{pa 5232 1834}%
\special{fp}%
\put(44.8200,-16.4600){\makebox(0,0)[lb]{$\wp$}}%
%
\special{pn 20}%
\special{ar 5262 1298 538 538  3.0378003 6.2831853}%
\special{ar 5262 1298 538 538  0.0000000 2.9271389}%
%
\special{pn 20}%
\special{pa 4356 1406}%
\special{pa 4734 1406}%
\special{fp}%
%
\special{pn 20}%
\special{pa 4176 1238}%
\special{pa 4908 1238}%
\special{fp}%
%
\special{pn 20}%
\special{ar 3822 1298 358 358  0.1545646 6.1148227}%
%
\special{pn 20}%
\special{ar 3822 1298 548 548  0.1999199 6.2831853}%
%
\special{pn 20}%
\special{pa 4374 1298}%
\special{pa 4902 1298}%
\special{fp}%
%
\special{pn 20}%
\special{ar 5268 1298 372 372  3.3067413 6.2831853}%
\special{ar 5268 1298 372 372  0.0000000 3.1415927}%
%
\special{pn 20}%
\special{pa 4614 1238}%
\special{pa 4578 1238}%
\special{fp}%
\special{sh 1}%
\special{pa 4578 1238}%
\special{pa 4646 1258}%
\special{pa 4632 1238}%
\special{pa 4646 1218}%
\special{pa 4578 1238}%
\special{fp}%
%
\special{pn 20}%
\special{pa 4620 1406}%
\special{pa 4578 1406}%
\special{fp}%
\special{sh 1}%
\special{pa 4578 1406}%
\special{pa 4646 1426}%
\special{pa 4632 1406}%
\special{pa 4646 1386}%
\special{pa 4578 1406}%
\special{fp}%
%
\special{pn 20}%
\special{pa 4470 1298}%
\special{pa 4512 1294}%
\special{fp}%
\special{sh 1}%
\special{pa 4512 1294}%
\special{pa 4442 1282}%
\special{pa 4458 1300}%
\special{pa 4448 1322}%
\special{pa 4512 1294}%
\special{fp}%
%
\special{pn 20}%
\special{sh 1.000}%
\special{ar 676 1278 38 38  0.0000000 6.2831853}%
%
\special{pn 20}%
\special{sh 1.000}%
\special{ar 2036 1278 38 38  0.0000000 6.2831853}%
\put(6.3000,-14.6000){\makebox(0,0)[lb]{$z_i$}}%
\put(20.0000,-14.8000){\makebox(0,0)[lb]{$z_j$}}%
%
\special{pn 20}%
\special{pa 2020 776}%
\special{pa 2082 770}%
\special{fp}%
\special{sh 1}%
\special{pa 2082 770}%
\special{pa 2014 758}%
\special{pa 2030 776}%
\special{pa 2018 796}%
\special{pa 2082 770}%
\special{fp}%
%
\special{pn 20}%
\special{pa 2550 1248}%
\special{pa 2556 1312}%
\special{fp}%
\special{sh 1}%
\special{pa 2556 1312}%
\special{pa 2570 1244}%
\special{pa 2550 1260}%
\special{pa 2530 1248}%
\special{pa 2556 1312}%
\special{fp}%
%
\special{pn 20}%
\special{pa 2066 1792}%
\special{pa 2014 1792}%
\special{fp}%
\special{sh 1}%
\special{pa 2014 1792}%
\special{pa 2080 1812}%
\special{pa 2066 1792}%
\special{pa 2080 1772}%
\special{pa 2014 1792}%
\special{fp}%
%
\special{pn 20}%
\special{pa 690 770}%
\special{pa 652 770}%
\special{fp}%
\special{sh 1}%
\special{pa 652 770}%
\special{pa 718 790}%
\special{pa 704 770}%
\special{pa 718 750}%
\special{pa 652 770}%
\special{fp}%
%
\special{pn 20}%
\special{pa 690 1800}%
\special{pa 752 1794}%
\special{fp}%
\special{sh 1}%
\special{pa 752 1794}%
\special{pa 682 1782}%
\special{pa 698 1800}%
\special{pa 688 1820}%
\special{pa 752 1794}%
\special{fp}%
%
\special{pn 20}%
\special{pa 150 1290}%
\special{pa 158 1358}%
\special{fp}%
\special{sh 1}%
\special{pa 158 1358}%
\special{pa 170 1290}%
\special{pa 152 1304}%
\special{pa 130 1294}%
\special{pa 158 1358}%
\special{fp}%
%
\special{pn 20}%
\special{pa 1050 924}%
\special{pa 1644 1596}%
\special{fp}%
%
\special{pn 20}%
\special{pa 1644 970}%
\special{pa 1050 1630}%
\special{fp}%
%
\special{pn 20}%
\special{ar 674 1278 516 516  0.7237480 5.6219733}%
%
\special{pn 20}%
\special{ar 2042 1278 502 502  3.7720777 6.2831853}%
\special{ar 2042 1278 502 502  0.0000000 2.4534819}%
\put(13.3000,-15.9000){\makebox(0,0)[lb]{$\ve$}}%
\end{picture}%
\end{center}
\caption{An eight-loop and a Pochhammer loop, where 
$z_1=0$, $z_2=1$, $z_3=\infty$} 
\label{fig:loop} 
\end{figure}
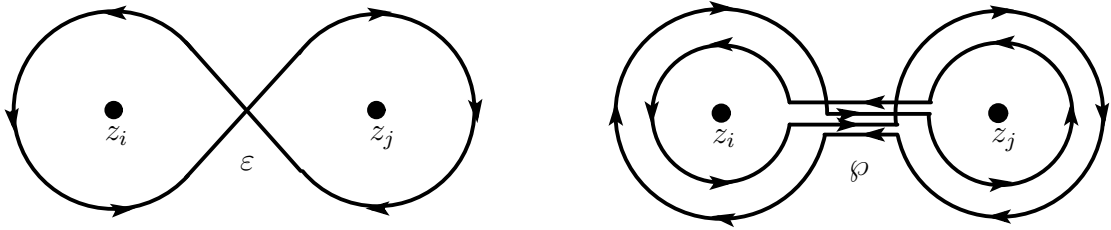 
They are called an {\sl eight-loop} and a {\sl Pochhammer loop} 
respectively. 
Their first dynamical degrees are given by 
\begin{equation} \label{eqn:DDEP}
\lambda(\ve)=3+2\sqrt{2}, \qquad \lambda(\wp)=9+4\sqrt{5}.
\end{equation}
\end{example}
\par
The first main theorem of this paper is concerned with the ergodic 
properties of the monodromy map for all parameters, which 
generalizes a main theorem in \cite{IU1} for generic parameters. 
Before stating it, we review some terminology from \cite{IU1}. 
\begin{definition} \label{def:chaos} 
The dynamical system of a holomorphic map $f : X \to X$ on a 
complex surface $X$ is said to be {\sl chaotic} if there exists an 
$f$-invariant Borel probability measure $\mu$ on $X$ such 
that the following conditions are satisfied: 
\begin{enumerate}
\item the measure $\mu$ has a positive entropy $h_{\mu}(f) > 0$; 
\item the measure $\mu$ is mixing, hyperbolic of saddle type and 
has a product structure with respect to local stable and unstable 
manifolds. 
Moreover, hyperbolic periodic points of $f$ are dense in 
the support of $\mu$. 
\end{enumerate}
\end{definition}
Let $\Omega(f)$ denote the nonwandering set of $f$. 
It is an $f$-invariant set serving as the hub for recurrent 
behaviors of $f$, that is, it contains the support of any 
$f$-invariant probability measure. 
\begin{theorem} \label{thm:main1} 
Let $\ga \in \pi_1(Z,z)$ be a non-elementary loop and 
$\ga_* : \mathcal{M}_z(\k) \carl$ be the monodromy map 
along the loop $\ga$. 
\begin{enumerate}
\item The nonwandering set $\Omega(\ga)$ of $\ga_*$ is compact in 
$\mathcal{M}_z(\k)$ and the trajectory of each initial point 
$Q \in \mathcal{M}_z(\k) \setminus \Omega(\ga)$ 
tends to infinity $\mathcal{Y}_z(\k)$ under the iterations of $\ga_*$. 
\item The map $\ga_* : \Omega(\ga) \carl$ is chaotic, that is, there 
exists a $\ga_*$-invariant Borel probability measure $\mu_{\ga}$ 
satisfying the conditions of Definition $\ref{def:chaos}$. 
The measure-theoretic entropy $h(\ga) := h_{\mu_{\ga}}(\ga_*)$ 
with respect to $\mu_{\ga}$ and the topological entropy 
$h_{\mathrm{top}}(\ga)$ of $\ga_*$ are given by 
\begin{equation} \label{eqn:ent}
h(\ga) = h_{\mathrm{top}}(\ga) = \log \lambda(\ga),
\end{equation}
where $\lambda(\ga)$ is the first dynamical degree of $\ga$ 
$($see Remark $\ref{rem:FDD}$$)$. 
In particular, $\mu_{\ga}$ is a unique maximal entropy measure. 
\end{enumerate}
\end{theorem}
\par
Next the second main theorem is concerned with the number of 
periodic solutions to $\PVI(\k)$. 
A periodic solution to $\PVI(\k)$ of period $n \in \N$ along 
a loop $\gamma$ is defined to be a periodic point of 
period $n$ of the map $\ga_* : \M_z(\k) \carl$. 
For some non-generic parameters, however, $\ga_*$ may admit 
curves of periodic points, in which case the set of 
periodic points is obviously uncountable and one should 
consider the set $\mathrm{Per}_n^i(\ga;\k)$ of {\sl isolated} 
periodic points of period $n$ instead. 
It will turn out that any periodic curve must be an irreducible 
component of the exceptional set $\E_z(\k)$ of the Riemann-Hilbert 
correspondence (\ref{eqn:RH2}) (see Theorem \ref{thm:main2}). 
Each irreducible component of $\E_z(\k)$ is known as a 
{\sl Riccati curve}, since $\E_z(\k)$ parameterizes all 
Riccati solution germs to $\PVI(\k)$ at $z$. 
They are classical special solutions that can be 
expressed in terms of Gauss hypergeoemtric functions 
(see \cite{IIS1,Okamoto,STe,Watanabe}). 
The second main theorem is now stated as follows. 
\begin{theorem} \label{thm:main2} 
Assume that $\ga \in \pi_1(Z,z)$ is non-elementary. 
Then any irreducible periodic curve of $\ga_*$ must be a Riccati curve 
and the set $\mathrm{Per}_n^i(\ga;\k)$ is finite for any $n \in \N$ and 
its cardinality $\# \mathrm{Per}_n^i(\gamma;\k)$ counted with 
multiplicity is estimated as 
\[
|\# \mathrm{Per}_n^i(\gamma;\k) - \lambda(\gamma)^n| < O(1) 
\qquad \mbox{as} \quad n \to \infty, 
\]
where $\lambda(\ga)$ is the first dynamical degree of $\ga$ and 
$O(1)$ stands for a bounded function of $n \in \N$. 
Moreover, if $\mathrm{HPer}_n(\ga;\k)$ denotes the set of saddle 
periodic points of period $n$ for $\ga_*$, then
\begin{equation} \label{eqn:asymp2}
\# \mathrm{Per}_n^i(\gamma;\k) \sim  \# \mathrm{HPer}_n(\ga;\k) 
\sim \lambda(\gamma)^n \qquad \mbox{as} \quad n \to \infty, 
\end{equation}
where $a_n \sim b_n$ means that the ratio $a_n/b_n \to 1$ as 
$n \to \infty$, so that asymptotically almost all points in 
$\mathrm{Per}_n^i(\gamma;\k)$ belong to $\mathrm{HPer}_n(\ga;\k)$. 
\end{theorem}
\begin{remark} \label{rem:multi}
In formula (\ref{eqn:asymp2}) the number 
$\# \mathrm{Per}_n^i(\gamma;\k)$ may be counted {\sl without} 
multiplicity, since $\mathrm{HPer}_n(\ga;\k) \subset \Per_n^i(\ga;\k)$ 
and each element of $\mathrm{HPer}_n(\ga;\k)$ is of simple multiplicity. 
\end{remark}
\begin{example} \label{ex:card}
Here are some explicit formulas for the number 
$\# \mathrm{Per}_n^i(\ga;\k)$. 
\begin{enumerate}
\item If the parameter $\k \in \K$ is generic, then 
$\mathcal{E}_z(\k)$ is empty and the number of isolated periodic 
solutions along any non-elementary loop $\ga \in \pi_1(Z,z)$ is 
calculated in \cite[Theorem 3]{IU1} as
\[
\# \mathrm{Per}_n^i(\ga;\k)= \lambda(\ga)^n +\lambda(\ga)^{-n}+4. 
\]
\item If $\k=(0,0,0,0,1) \in \K$, which is non-generic, then 
$\mathcal{E}_z(\k)$ consists of four Riccati curves and 
the numbers of isolated periodic solutions along an eight-loop 
$\ve$ and a Pochhammer loop $\wp$ are calculated as 
\[
\# \mathrm{Per}_n^i(\ve;\k)= \lambda(\ve)^n +\lambda(\ve)^{-n}, 
\qquad 
\# \mathrm{Per}_n^i(\wp;\k)= \lambda(\wp)^n +\lambda(\wp)^{-n}-4, 
\]
respectively, where $\lambda(\ve)$ and $\lambda(\wp)$ are given by 
formula (\ref{eqn:DDEP}). 
Note that $\ve_*$ fixes exactly one Riccati curve, while $\wp_*$ 
fixes all the four Riccati curves. 
\end{enumerate}
In Section \ref{sec:Ex} the number $\# \mathrm{Per}_n^i(\wp;\k)$ 
will be calculated for various values of $\k \in \K$. 
In general the number $\# \mathrm{Per}_n^i(\ga;\k)$ is 
computable at least in principle, once the data of a 
non-elementary loop $\ga \in \pi_1(Z,z)$, a parameter 
$\k \in \K$ and a period $n \in \N$ is given explicitly. 
\end{example}
\par
The discussion above is about periodic solutions along 
a single loop in $\pi_1(Z,z)$. 
One can also think of priodic solutions with respect to the 
entire $\pi_1(Z,z)$-action, that is, global solutions with 
finitely many branches. 
They are exactly algebraic solutions to Painlev\'e VI 
(\cite{Iwasaki1}), which have been classified by 
\cite{LT} after many contributions by various authors 
(see also \cite{Boalch,Iwasaki2}). 
\par
This article is organized as follows. 
In Section \ref{sec:Pre} we review some general theories 
of complex surface dynamics that will be needed to prove 
our main results. 
In Section \ref{sec:cubic} a four-parameter family of 
affine cubic surfaces is introduced as the target spaces 
of the Riemann-Hilbert correspondence and the possible types 
of singularities on them are classified. 
After investigating polynomial automorphisms on the cubic 
surfaces in Section \ref{sec:dynamics}, we establish 
our main results in Section \ref{sec:MT}. 
Finally, Section \ref{sec:Ex} is devoted to a thorough 
study of isolated periodic solutions along a 
Pochhammer loop for various parameters. 
\section{Preliminaries} \label{sec:Pre} 
In this section we collect some basic results from complex 
surface dynamics that will be needed later. 
In order to establish Theorem \ref{thm:main1} it is necessary 
to construct an invariant measure satisfying the conditions in 
Definition \ref{def:chaos}. 
The first part of this section is a review on the construction 
of such measures for a class of bimeromorphic maps on smooth 
surfaces. 
The second part is a survey on a general theory of isolated 
periodic points for area-preserving surface maps admitting 
periodic curves, which will be used to prove 
Theorem \ref{thm:main2}. 
All along the way the following two concepts for bimeromorphic 
surface maps due to \cite{FS,DF} are important. 
\begin{definition} \label{def:birat}
Let $f : X \to X$ be a bimeromorphic map on a compact 
K\"ahler surface $X$. 
\begin{enumerate}
\item {\sl The first dynamical degree} $\lambda(f)$ is 
defined by 
\[
\lambda(f) := \lim_{n \to \infty} 
|\!|(f^n)^*|_{H^{1,1}(X)}|\!|^{1/n} \ge 1,
\]
where $|| \cdot ||$ is an operator norm on 
$\mathrm{End} H^{1,1}(X)$. 
It is known that the limit exists, $\lambda(f)$ is 
independent of the norm $|| \cdot ||$ chosen and invariant 
under bimeromorphic conjugation. 
\item The map $f$ is said to be {\sl analytically stable} 
(AS for short) if the condition 
$(f^n)^* = (f^*)^n : H^{1,1}(X) \to H^{1,1}(X)$ holds for 
any $n \in \N$. 
It is known that $f$ is AS if and only if 
\begin{equation} \label{eqn:AS}
f^{-m}I(f) \cap f^n I(f^{-1}) = \emptyset 
\qquad \mbox{for every} \quad 
m, n \ge 0, 
\end{equation}
where $I(f)$ is the indeterminacy set of $f$. 
If $f$ is AS then the first dynamical degree $\lambda(f)$ 
coincides with the spectral radius of the linear map 
$f^*|_{H^{1,1}(X)}$. 
\end{enumerate}
\end{definition}
\par 
We begin with a review on invariant measures. 
Under the condition 
\begin{equation} \label{eqn:fdd2}
\l(f) > 1, 
\end{equation}
Bedford and Diller \cite{BD} constructed ``good'' positive 
closed $(1,1)$-currents $\mu^{\pm}$ on $X$ such that 
\[
(f^{\pm1})^* \mu^{\pm} = \l(f) \, \mu^{\pm}, 
\]
where $\mu^+$ and $\mu^-$ are called the stable and unstable 
currents for $f$. 
They represent unique (up to scale) nef classes 
$\th^{\pm} \in H^{1,1}_{\R}(X)$ such that 
$(f^{\pm1})^* \th^{\pm} = \l(f) \, \th^{\pm}$ and 
can be expressed as 
\[
\mu^{\pm}= 
c^{\pm}_{\omega} 
\lim_{k \to \infty} \l(f)^{-k} (f^{\pm k})^* \omega, 
\]
where $\omega$ is any given smooth closed $(1,1)$-form on $X$ 
and $c^{\pm}_{\omega} > 0$ are constants. 
Moreover, under a quantitative condition on the indeterminacy 
sets of the forward and backward maps: 
\begin{equation} \label{eqn:AS3} 
\sum_{N = 0}^{\infty} \l(f)^{-N} 
\log \, \mathrm{dist}(f^N I(f^{-1}), I(f)) > - \infty, 
\end{equation} 
Bedford and Diller \cite{BD} and Dujardin \cite{Dujardin} 
legitimated the wedge product $\mu := \mu^+ \wedge \mu^-$ 
as an $f$-invariant Borel probability measure, where the first 
authors defined it by appealing to pluripotential theory while 
the second author viewed it as geometric intersection. 
The condition (\ref{eqn:AS3}) is slightly stronger than 
(\ref{eqn:AS}) and a map satisfying condition (\ref{eqn:AS3}) 
might be called {\sl quantitatively AS}. 
The measure has good dynamical properties 
as is mentioned in the following. 
\begin{theorem}[\cite{BD, Dujardin}] \label{thm:BD} 
If $f : X \carl$ satisfies conditions $(\ref{eqn:fdd2})$ and 
$(\ref{eqn:AS3})$, then the wedge product $\mu$ of stable 
and unstable currents $\mu^{\pm}$ is well defined and, after a 
suitable normalization, $\mu$ gives an $f$-invariant Borel 
probability measure satisfying all the conditions in 
Definition $\ref{def:chaos}$. 
Moreover, 
\begin{enumerate}
\item the measure-theoretic entropy $h_{\mu}(f)$ with respect 
to the measure $\mu$ and the topological entropy 
$h_{\mathrm{top}}(f)$ of $f$ are expressed as
\[
h_{\mu}(f) = h_{\mathrm{top}}(f) = \log \lambda(f), 
\]
\item the measure $\mu$ puts no mass on any algebraic curve on $X$, 
\item there exists a set $\mathcal{P}_n(f) \subset \mathrm{supp}\, 
\mu$ of saddle periodic points of period $n$ such that 
\[
\# \mathcal{P}_n(f) \sim \lambda(f)^n, \qquad 
\dfrac{1}{\lambda(f)^n} \sum_{p \in \mathcal{P}_n(f)} \delta_p 
\to \mu, \qquad \mbox{as} \quad n \to \infty. 
\]
\end{enumerate}
\end{theorem} 
\begin{remark} \label{rem:nonwander} 
The definition of entropy needs some care for a 
bimeromorphic surface map $f : X \carl$, since it may have 
indeterminacy sets $I(f^{\pm1})$ on which $f^{\pm1}$ are 
not defined. 
To handle this situation, notice that $f$ restricts to a 
well-defined automorphism $f|_{X_f}$ of the space 
\begin{equation} \label{eqn:invset}
X_f := X \setminus \bigr(\bigcup_{n \ge 0} f^n I(f) \cup f^{-n} 
I(f^{-1}) \bigl). 
\end{equation}
If a Borel probability measure $\nu$ on $X$ satisfies 
$\nu(X_f) = 1$, then the measure-theoretic entropy of $f$ with 
respect to $\nu$ can be defined by $h_{\nu}(f):=h_{\nu}(f|_{X_f})$ 
in terms of the map $f : X_f \carl$ (see \cite{Guedj}). 
The measure $\mu$ constructed in Theorem \ref{thm:BD} satisfies 
this condition and so the entropy $h_{\mu}(f)$ is well defined. 
In the same spirit the topological entropy of $f$ is defined by 
$h_{\mathrm{top}}(f):=h_{\mathrm{top}}(f|_{X_f})$, where the 
right-hand side employs Bowen's definition on a non-compact 
space (see \cite{Bowen}). 
Similarly the nonwandering set $\Omega(f)$ of $f$ is that of 
$f|_{X_f}$, i.e., $\Omega(f) := \Omega(f|_{X_f})$. 
\end{remark}
\par
We proceed to a survey of the results of \cite{IU2} on the number 
of isolated periodic points for area-preserving surface maps 
with periodic curves. 
Let $f : X \carl$ be a birational map of a smooth projective 
surface $X$. 
Since $f$ may have indeterminacy sets $I(f^{\pm1})$ on which 
$f^{\pm1}$ are not defined, we should again be careful with 
the definitions of a fixed point and a fixed curve. 
\begin{definition} \label{def:fixed} 
A point $x \in X$ is called a {\sl fixed point} if $x$ is an 
element of the set 
\[
X_0(f) := X_0^{\circ}(f) \cup X_0^{\circ}(f^{-1}), 
\]
where $X_0^{\circ}(f)$ is the set of all points 
$x \in X \setminus I(f)$ fixed by $f$. Moreover let 
$X_1(f)$ be the set of all irreducible curves $C$ in $X$ 
such that $C \setminus I(f)$ is fixed pointwise by $f$. 
An element $C \in X_1(f)$ is called a {\sl fixed curve}. 
It is easy to see that the definition is symmetric, namely, 
$X_1(f) = X_1(f^{-1})$. 
As in \cite{IU2,S} the set $X_1(f)$ of fixed curves is 
divided into two disjoint subsets: 
\begin{equation} \label{eqn:divide}
X_1(f) = X_I(f) \amalg X_{I\!I}(f), 
\end{equation}
namely, into the curves of type $I$ and those of type $I\!I$ 
(see Definition \ref{def:type2}). 
\end{definition}
\par
The set of periodic points of period $n$ for $f$ is defined by 
$\mathrm{Per}_n(f) := X_0(f^n)$. 
The subset of {\sl isolated} ones is denoted by 
$\mathrm{Per}_n^{i}(f)$ and its cardinality counted with 
multiplicity is defined by 
\[
\# \mathrm{Per}_n^{i}(f) :=\sum_{x \in \mathrm{Per}_n^{i}(f)} \nu_x(f^n),
\]
where $\nu_x(f^n)$ is the local index of $f^n$ at $x \in X_0(f^n)$ to be 
defined in Definition \ref{def:index}. 
\begin{theorem}[\cite{IU2}] \label{thm:esti} 
Let $f : X \to X$ be an AS birational map with $\lambda(f) > 1$ on a 
smooth projective surface and assume that $f$ preserves a nontrivial 
meromorphic $2$-form $\omega$ such that no irreducible component of 
the pole divisor $(\omega)_{\infty}$ of $\omega$ is a periodic curve of 
type $I$. Then $f$ has at most finitely many irreducible periodic curves 
and must have infinitely many isolated periodic points. 
Moreover the cardinality $\# \mathrm{Per}_n^{i}(f)$ is estimated as
\[
|\# \mathrm{Per}_n^{i}(f) - \lambda(f)^n| \le 
\left\{ \begin{array}{ll}
O(1) \quad & (\mbox{$X \sim$ no Abelian surface}), \\[2mm] 
4 \, \lambda(f)^{n/2} + O(1) \quad & 
(\mbox{$X \sim$ an Abelian surface}), 
\end{array} \right. 
\]
where $X \sim Y$ indicates that $X$ is birationally equivalent to $Y$ 
and $O(1)$ is a bounded function of $n \in \N$. 
\end{theorem}
\par
Let $\mathrm{HPer}_n(f)$ be the set of all saddle periodic points of 
period $n$ and $\mathcal{P}_n(f)$ the set of saddle periodic 
points mentioned in Theorem \ref{thm:BD}. 
Then one has $\mathcal{P}_n(f) \subset \mathrm{HPer}_n(f) \subset 
\mathrm{Per}^i_n(f)$, since any saddle periodic point is isolated. 
Thus Theorems \ref{thm:BD} and \ref{thm:esti} have the following. 
\begin{corollary} \label{cor:esti} 
If $f$ satisfies the assumptions in Theorems $\ref{thm:BD}$ and 
$\ref{thm:esti}$, then
\begin{equation} \label{eqn:asymp}
\# \mathrm{Per}_n^i(f) \sim  \# \mathrm{HPer}_n(f) 
\sim \lambda(f)^n \qquad \mbox{as} \quad n \to \infty, 
\end{equation}
so that asymptotically almost all points in $\mathrm{Per}_n^i(f)$ 
belong to $\mathrm{HPer}_n(f)$. 
\end{corollary}
\begin{remark} \label{rem:count}
In formula (\ref{eqn:asymp}) the number $\# \mathrm{Per}_n^i(f)$ 
may be counted {\sl without} multiplicity, because $\mathrm{HPer}_n(f)$ 
is a subset of $\mathrm{Per}_n^i(f)$, every point in $\mathrm{HPer}_n(f)$ 
is of simple multiplicity, and the asymptotics 
$\# \mathrm{HPer}_n(f) \sim \lambda(f)^n$ holds with multiplicity 
taken into account.
\end{remark}
\par
The above results are derived from a basic formula representing the 
Lefschetz numbers of iterates of $f$ in terms of suitable local data 
around isolated periodic points as well as around periodic curves 
of $f$ (see Theorem \ref{thm:formula}). 
Here the Lefschetz number of $f$ is defined by 
\[
L(f) := \sum_{i} (-1)^i \mathrm{Tr}[f^* : H^i(X) \to H^i(X)]. 
\]
In order to state that formula, let $P(f)$ be the set of all 
positive integers that arise as primitive periods of some 
irreducible periodic curves of $f$. 
For each $n \in \N$, denote by $P_n(f)$ the set of all elements 
$k \in P(f)$ that divides $n$. 
Moreover, for each $k \in P(f)$ let $\mathrm{PC}_k(f)$ be 
the set of all irreducible periodic curves  of primitive 
period $k$, and $C_k(f)$ the union of all curves in 
$\mathrm{PC}_k(f)$. 
Note that there exist the following decompositions: 
\[
X_0(f^n) = \mathrm{Per}_n^{i}(f) \amalg 
\bigcup_{k \in P_n(f)} C_k(f), \qquad 
X_1(f^n) = \coprod_{k \in P_n(f)} \mathrm{PC}_k(f). 
\]
For each $k \in P(f)$ let $\xi_k(f)$ be the number 
defined by 
\[
\xi_k(f) := \sum_{x \in C_k(f)} \nu_x(f^k) + 
\sum_{C \in \mathrm{PC}_k(f)} \tau_{C} \cdot \nu_{C}(f^k). 
\]
where $\tau_{C}$ is the self-intersection number of $C$. 
Then our formula is stated as follows. 
\begin{theorem}[\cite{IU2}] \label{thm:formula} 
If $f : X \carl$ satisfies the assumptions in 
Theorem $\ref{thm:esti}$, then
\[
L(f^n) = \# \mathrm{Per}_n^{i}(f) + \sum_{k \in P_n(f)} \xi_k(f) 
\qquad (n \in \N).  
\]
\end{theorem}
This formula is used not only to get the general estimate 
in Theorem \ref{thm:esti} but also to calculate the exact 
value of $\# \mathrm{Per}_n^{i}(f)$ for various individual 
maps $f$ (see Section \ref{sec:Ex}). 
\par
Now let us recall the definitions of $\nu_x(f)$, $\nu_C(f)$, 
$X_I(f)$ and $X_{I\!I}(f)$. 
Leaving the general cases in \cite{IU2,S} we put the following 
assumption for the sake of simplicity. 
It will be fulfilled by the birational maps on (desingularized) 
cubic surfaces to be discussed later (see Remark \ref{rem:dynkin}). 
\begin{assumption} \label{ass:fix}
All fixed curves of $f : X \carl$ are smooth, no two of which 
are tangent and no three of which meet in a single point. 
\end{assumption}
\par
For a given point $x \in X_0^{\circ}(f)$ let $A_x$ be the 
completion of the local ring of $X$ at $x$, which can be 
identified with the formal power series ring $\C[\![x_1,x_2]\!]$ 
of two variables, because $X$ is smooth. 
Since $f$ is holomorphic around $x$, $f$ induces an endomorphism 
$f_x^* : A_x \to A_x$ in a natural manner. 
From Assumption \ref{ass:fix}, $f_x^*$ can be expressed in 
suitable coordinates $(x_1,x_2)$ as 
\[
\left\{
\begin{array}{rcl}
f_x^*(x_1) &=& x_1 + x_1^{n_1} \cdot x_2^{n_2} \cdot h_1, \\[1mm]
f_x^*(x_2) &=& x_2 + x_1^{n_1} \cdot x_2^{n_2} \cdot h_2, 
\end{array}
\right. 
\]
with some relatively prime elements $h_1$, $h_2 \in A_x$ 
and some nongegative integers $n_1$, $n_2 \in \Z_{\ge0}$. 
For $\{j,k\}=\{1,2\}$ let $\tau_{(x_j)} : 
\C[\![x_1,x_2]\!] \to \C[\![x_k]\!]$ denote the natural 
projection. 
Write 
\[
\tau_{(x_j)}(h_i) = x_k^{n_{ij}} \cdot h_{ij} \qquad 
(i,j \in \{1,2\}), 
\]
where $h_{ij}$ is a unit in $\C[\![x_k]\!]$ and $n_{ij}$ is 
either an integer or infinity. 
By convention $x_k^{\infty} := 0$. 
\begin{definition} \label{def:type} 
The prime ideal $\p_i := (x_i)$ generated by $x_i$ is said to 
be {\sl of type I} relative to $f_x^*$ if $n_{ii} \neq \infty$; 
otherwise, $\p_i$ is said to be {\sl of type I\!I}. 
\end{definition}
\par
For each $i \in \{1,2\}$ we put 
\begin{eqnarray}
\nu_{\p_i}(f_x^*) &:= & n_i, \label{eqn:nup} 
\\[2mm]
\nu_{A_x}(f_x^*) &:= & \dim_{\C} \, A_x/(h_1, h_2)  + 
\displaystyle \sum_{i=1}^{2} 
\nu_{\p_i}(f_x^*) \cdot \mu_{\p_i}(f_x^*), \label{eqn:nuA} 
\end{eqnarray} 
where $A_x/(h_1, h_2)$ is the quotient vector space of $A_x$ by 
the ideal $(h_1, h_2)$, which is finite-dimensional from 
Assumption \ref{ass:fix}, and the numbers $\mu_{\p_i}(f_x^*)$ 
are defined by 
\[
\mu_{\p_i}(f_x^*) := 
\left\{
\begin{array}{rcll}
n_{ii} \qquad & 
(\,\mbox{if} \,\, \p_i \,\, \mbox{is of type I}), \\[2mm]
n_{ji} \qquad & 
(\,\mbox{if} \,\, \p_i \,\, \mbox{is of type I\!I}), 
\end{array}
\right. 
\]
with $\{i,j\} = \{1,2\}$. 
Similarly for a given point $x \in X_0^{\circ}(f^{-1})$ one can 
define the number $\nu_{A_x}((f^{-1})_x^*)$ via formula 
(\ref{eqn:nuA}). 
Next, given a fixed curve $C \in X_1(f)$, take a point 
$x$ of $C \setminus I(f)$. 
Then one can speak of the endomorphism $f_x^* : A_x \to A_x$. 
Since $C$ is smooth by Assumption \ref{ass:fix}, the prime 
ideal $\p_C \subset A_x$ defining the germ at $x$ of the 
curve $C$ may be written $\p_C = (x_1)$ in suitable 
coordinates $(x_1,x_2)$, so that one can define the number 
$\nu_{\p_C}(f_x^*)$ via formula (\ref{eqn:nup}). 
\begin{definition} \label{def:index} 
The local index $\nu_x(f)$ at a fixed point $x \in X_0(f)$ is 
defined by 
\[
\nu_x(f) := \left\{
\begin{array}{ll}
\nu_{A_x}(f_x^*) \qquad & (\,\mbox{if $x \in X_0^{\circ}(f)$}), 
\\[2mm] 
\nu_{A_x}((f^{-1})_x^*) \qquad & 
(\,\mbox{if $x \in X_0^{\circ}(f^{-1})$}), 
\end{array}\right. 
\]
where the right-hand side is consistent, that is, 
$\nu_{A_x}(f_x^*) = \nu_{A_x}((f^{-1})_x^*)$ for any 
$x \in X_0^{\circ}(f) \cap 
X_0^{\circ}(f^{-1})$ (see \cite{IU2}). 
The index $\nu_C(f)$ at a fixed curve $C \in X_1(f)$ is 
defined by 
\[
\nu_C(f) := \nu_{\p_C}(f_x^*), 
\]
with $x \in C \setminus I(f)$, where the right-hand side does 
not depend on the choice of $x$ (see \cite{S}). 
\end{definition} 
\par
Finally we recall the following definition concerning the types 
of fixed curves. 
\begin{definition} \label{def:type2} 
A fixed curve $C \in X_1(f)$ is said to be 
{\sl of type I} or {\sl of type I\!I} relative to 
$f : X \to X$, according as the prime ideal $\p_C$ is 
of type I or of type I\!I relative to $f_x^* : A_x \to A_x$ 
in the sense of Definition \ref{def:type}. 
This definition does not depend on the choice of 
$x \in C \setminus I(f)$ (see \cite{S}). 
Let $X_I(f)$ and $X_{I\!I}(f)$ denote the set of fixed curves 
of types I and that of type I\!I respectively. 
Then there exists the direct sum decomposition (\ref{eqn:divide}). 
\end{definition} 
\section{Singular Cubic Surfaces} \label{sec:cubic} 
The purpose of this section is to introduce a four-parameter 
family of affine cubic surfaces which are the target spaces 
of the Riemann-Hilbert correspondence (\ref{eqn:RH2}) and to 
classify the types of singularities that can occur on those 
surfaces. 
The affine cubic surfaces we consider are
\[
\Sol(\th) = \{\, x = (x_1,x_2,x_3) \in \C^3 \,:\, 
f(x,\th) = 0 \, \}, 
\]
where $f(x,\th)$ is a cubic polynomial of $x$ with parameters 
$\th = (\th_1,\th_2,\th_3,\th_4) \in \Th := \C^4$: 
\[
f(x,\th) := x_1x_2x_3 + x_1^2 + x_2^2 + x_3^2 
- \th_1 x_1 - \th_2 x_2 - \th_3 x_3 + \th_4. 
\] 
\par
Depending on parameters $\th \in \Th$ the surface $\Sol(\th)$ 
may admit singular points. 
In order to describe its singularity structure, 
it is convenient to introduce a map 
\begin{equation} \label{eqn:rhp}
\mathrm{rh} : \K \to \Th, 
\end{equation}
called {\sl the Riemann-Hilbert correspondence in the parameter 
level}. 
It is the composite of maps 
\begin{equation} \label{eqn:KBATh}
\begin{CD}
\K @> \beta >> B @> \alpha >> A @> \varphi >> \Th,
\end{CD}
\end{equation}
where the intermediate parameter spaces $A$ and $B$ are 
given by 
\[
\begin{array}{rcl}
A &:=& \{ a=(a_1,a_2,a_3,a_4) \in \C_a^4 \}, \\[2mm]
B &:=& \{ b = (b_0,b_1,b_2,b_3,b_4) \in (\C_b^{\times})^5 
\, : \, b_0^2 \, b_1 \, b_2 \, b_3 \, b_4 = 1 \}, 
\end{array}
\]
and the three maps $\varphi$, $\alpha$ and $\beta$ are 
defined respectively by 
\begin{eqnarray}
\th_i &=& \left\{
\begin{array}{ll}
a_i a_4 + a_j a_k \qquad & (\{i,j,k\}=\{1,2,3\}), \\[2mm]
a_1 a_2 a_3 a_4 + a_1^2 + a_2^2 + a_3^2 + a_4^2 - 4 
\qquad & (i = 4),  
\end{array}\right. \nonumber \\[2mm]
a_i &=& b_i + b_i^{-1} \qquad (i=1,2,3,4), \nonumber \\[2mm]
b_i &=& \left\{
\begin{array}{rl}
\exp(\sqrt{-1} \pi \k_i) \qquad & (i=0,1,2,3), \\[2mm]
-\exp(\sqrt{-1} \pi \k_4)\qquad & (i = 4).  
\end{array} \right. \label{eqn:b} 
\end{eqnarray}
\par 
It turns out that the discriminant $\vD(\th)$ of the cubic 
surface $\Sol(\th)$ factors as 
\begin{equation} \label{eqn:vD}
\displaystyle \vD(\th) = \prod_{l=1}^4(b_l-b_l^{-1})^2 
\prod_{\ve \in \{\pm1\}^4}(b^{\ve}-1) 
=\prod_{l=1}^4 \sin^2 \pi \k_l 
\prod_{\ve \in \{\pm1\}^4} \cos \dfrac{\pi (\ve \cdot \k)}{2}, 
\end{equation}
where $b^{\ve} := b_1^{\ve_1}b_2^{\ve_2}b_3^{\ve_3}b_4^{\ve_4}$ 
and $\ve \cdot \k := \ve_1 \k_1+\ve_2 \k_2+\ve_3 \k_3+\ve_4 \k_4$
for each $\ve = (\ve_1,\ve_2,\ve_3,\ve_4) \in 
\{\pm1\}^4$. 
Thus $\Sol(\th)$ is singular if and only if $\k$ satisfies 
at least one of the affine linear relations: 
\begin{equation} \label{eqn:walls}
\k_i = m, \qquad  \ve \cdot \k = 2m+1 
\qquad (m \in \Z, \, i \in \{1,2,3,4\}, \, \ve \in \{\pm1\}^4). 
\end{equation}
\par
Behind formulas (\ref{eqn:vD}) and (\ref{eqn:walls}) there 
exists an affine Weyl group structure on the Riemann-Hilbert 
correspondence in the parameter level (\ref{eqn:rhp}). 
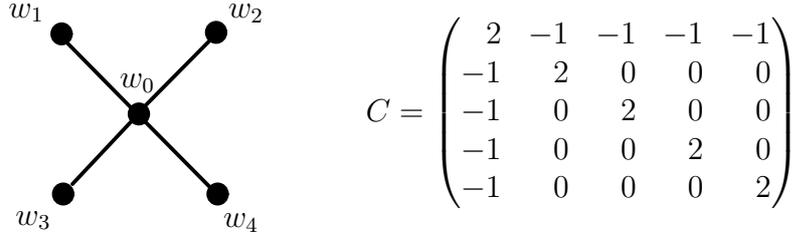
\begin{figure}[t]
\begin{center}
\unitlength 0.1in
\begin{picture}( 42.7000, 11.1600)(  2.2000,-17.6000)
%
\special{pn 20}%
\special{sh 1.000}%
\special{ar 494 876 48 50  0.0000000 6.2831853}%
%
\special{pn 20}%
\special{sh 1.000}%
\special{ar 1294 868 48 50  0.0000000 6.2831853}%
%
\special{pn 20}%
\special{sh 1.000}%
\special{ar 502 1712 48 50  0.0000000 6.2831853}%
%
\special{pn 20}%
\special{sh 1.000}%
\special{ar 1302 1712 48 50  0.0000000 6.2831853}%
%
\special{pn 20}%
\special{sh 1.000}%
\special{ar 894 1294 48 50  0.0000000 6.2831853}%
%
\special{pn 20}%
\special{pa 526 926}%
\special{pa 862 1270}%
\special{fp}%
\special{pa 870 1270}%
\special{pa 870 1270}%
\special{fp}%
%
\special{pn 20}%
\special{pa 926 1334}%
\special{pa 1262 1676}%
\special{fp}%
\special{pa 1262 1676}%
\special{pa 1262 1676}%
\special{fp}%
%
\special{pn 20}%
\special{pa 1262 910}%
\special{pa 926 1252}%
\special{fp}%
%
\special{pn 20}%
\special{pa 854 1334}%
\special{pa 550 1660}%
\special{fp}%
\put(2.2000,-8.1400){\makebox(0,0)[lb]{$w_1$}}%
\put(13.5800,-8.1700){\makebox(0,0)[lb]{$w_2$}}%
\put(2.5700,-18.9700){\makebox(0,0)[lb]{$w_3$}}%
\put(13.3200,-19.1000){\makebox(0,0)[lb]{$w_4$}}%
\put(7.9500,-11.8000){\makebox(0,0)[lb]{$w_0$}}%
%
\put(22.0000,-11.3000){\makebox(0,0)[lb]{}}%
\put(20.7000,-17.9000){\makebox(0,0)[lb]{$C= \begin{pmatrix} \mi 2 & -1 & -1 & -1 & -1 \\ -1 & \mi 2 & \mi 0 & \mi 0 & \mi 0 \\ -1 & \mi 0 & \mi 2 & \mi 0 & \mi 0 \\ -1 & \mi 0 & \mi 0 & \mi 2 & \mi 0 \\ -1 & \mi 0 & \mi 0 & \mi 0 & \mi 2 \end{pmatrix}$}}%
\put(44.9000,-12.1000){\makebox(0,0)[lb]{$\phantom{a}$}}%
\end{picture}%
\end{center}
\caption{Dynkin diagram and Cartan matrix of type $D_4^{(1)}$} 
\label{fig:dynkin} 
\end{figure}
The affine space $\K$ carries the inner product induced from 
the standard complex Euclidean inner product on $\C^4$ via the 
isomorphism $\K \stackrel{\sim}{\to} \C^4$, 
$\k \mapsto (\k_1,\k_2,\k_3,\k_4)$. 
For each $i \in \{0,1,2,3,4\}$ let $w_i : \K \carl$ be the 
orthogonal affine reflection in the affine hyperplane 
$H_i:=\{\k \in \K : \k_i = 0\}$. 
Explicitly $w_i$ is given by 
\[
w_i(\k_j) = \k_j - \k_i c_{ij}, 
\]
where $C = (c_{ij})$ is the Cartan matrix of type $D_4^{(1)}$ 
in Figure \ref{fig:dynkin} (right). 
The group 
\[ 
W(D_4^{(1)}) := \la w_0,w_1,w_2,w_3,w_4 \ra 
\]
generated by $w_0$, $w_1$, $w_2$, $w_3$, $w_4$ is an affine 
Weyl group of type $D_4^{(1)}$. 
The hyperplanes in (\ref{eqn:walls}) are the reflection 
hyperplanes of $W(D_4^{(1)})$ and the map (\ref{eqn:rhp}) 
is a branched $W(D_4^{(1)})$-covering ramifying along these 
hyperplanes. 
The automorphism group of the Dynkin diagram in Figure 
\ref{fig:dynkin} (left) is the symmetric group $S_4$ 
of degree $4$, which acts on $\K$ by permuting 
$\k_1, \k_2, \k_3, \k_4$ and fixing $\k_0$. 
The group $W(D_4^{(1)})$ extends to an affine Weyl group 
of type $F_4^{(1)}$: 
\[
W(F_4^{(1)}) := S_4 \ltimes W(D_4^{(1)}).
\]
\par
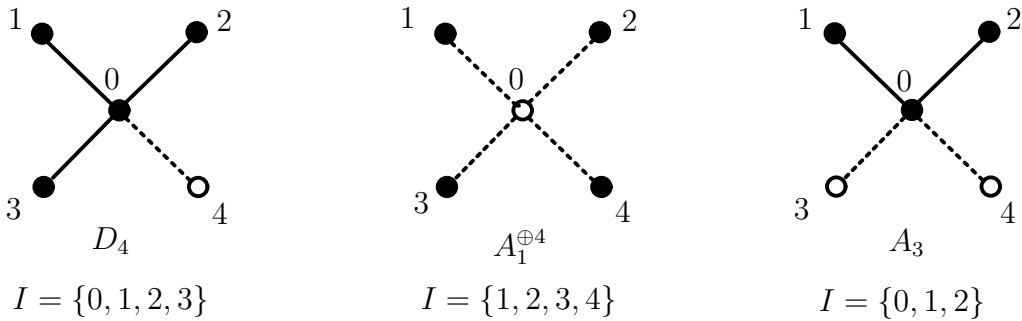
\begin{figure}[b]
\begin{center}
\unitlength 0.1in
\begin{picture}( 51.8000, 15.2000)(  2.6000,-21.7000)
%
\special{pn 20}%
\special{sh 1.000}%
\special{ar 448 852 48 48  0.0000000 6.2831853}%
%
\special{pn 20}%
\special{sh 1.000}%
\special{ar 1248 844 48 48  0.0000000 6.2831853}%
%
\special{pn 20}%
\special{sh 1.000}%
\special{ar 456 1652 48 48  0.0000000 6.2831853}%
%
\special{pn 20}%
\special{ar 1256 1652 48 48  0.0000000 6.2831853}%
%
\special{pn 20}%
\special{sh 1.000}%
\special{ar 848 1252 48 48  0.0000000 6.2831853}%
%
\special{pn 20}%
\special{pa 480 900}%
\special{pa 816 1228}%
\special{fp}%
\special{pa 824 1228}%
\special{pa 824 1228}%
\special{fp}%
%
\special{pn 20}%
\special{pa 880 1292}%
\special{pa 1216 1620}%
\special{dt 0.054}%
\special{pa 1216 1620}%
\special{pa 1216 1620}%
\special{dt 0.054}%
%
\special{pn 20}%
\special{pa 1216 884}%
\special{pa 880 1212}%
\special{fp}%
%
\special{pn 20}%
\special{pa 808 1292}%
\special{pa 504 1604}%
\special{fp}%
\put(2.7000,-8.3000){\makebox(0,0)[lb]{$\sc{1}$}}%
\put(13.5000,-8.4000){\makebox(0,0)[lb]{$\sc{2}$}}%
\put(2.6000,-18.2000){\makebox(0,0)[lb]{$\sc{3}$}}%
\put(13.2800,-18.4400){\makebox(0,0)[lb]{$\sc{4}$}}%
\put(7.6800,-11.4000){\makebox(0,0)[lb]{$\sc{0}$}}%
%
\special{pn 20}%
\special{sh 1.000}%
\special{ar 2536 852 48 48  0.0000000 6.2831853}%
%
\special{pn 20}%
\special{sh 1.000}%
\special{ar 3336 844 48 48  0.0000000 6.2831853}%
%
\special{pn 20}%
\special{sh 1.000}%
\special{ar 2544 1652 48 48  0.0000000 6.2831853}%
%
\special{pn 20}%
\special{sh 1.000}%
\special{ar 3344 1652 48 48  0.0000000 6.2831853}%
%
\special{pn 20}%
\special{ar 2936 1252 48 48  0.0000000 6.2831853}%
%
\special{pn 20}%
\special{pa 2568 900}%
\special{pa 2904 1228}%
\special{dt 0.054}%
\special{pa 2912 1228}%
\special{pa 2912 1228}%
\special{dt 0.054}%
%
\special{pn 20}%
\special{pa 2968 1292}%
\special{pa 3304 1620}%
\special{dt 0.054}%
\special{pa 3304 1620}%
\special{pa 3304 1620}%
\special{dt 0.054}%
%
\special{pn 20}%
\special{pa 3304 884}%
\special{pa 2968 1212}%
\special{dt 0.054}%
%
\special{pn 20}%
\special{pa 2896 1292}%
\special{pa 2592 1604}%
\special{dt 0.054}%
\put(23.0000,-8.2000){\makebox(0,0)[lb]{$\sc{1}$}}%
\put(34.5000,-8.4000){\makebox(0,0)[lb]{$\sc{2}$}}%
\put(23.7000,-17.9000){\makebox(0,0)[lb]{$\sc{3}$}}%
\put(34.1600,-18.3600){\makebox(0,0)[lb]{$\sc{4}$}}%
\put(28.6400,-11.4000){\makebox(0,0)[lb]{$\sc{0}$}}%
%
\special{pn 20}%
\special{sh 1.000}%
\special{ar 4552 850 48 48  0.0000000 6.2831853}%
%
\special{pn 20}%
\special{sh 1.000}%
\special{ar 5352 842 48 48  0.0000000 6.2831853}%
%
\special{pn 20}%
\special{ar 4560 1650 48 48  0.0000000 6.2831853}%
%
\special{pn 20}%
\special{ar 5360 1650 48 48  0.0000000 6.2831853}%
%
\special{pn 20}%
\special{sh 1.000}%
\special{ar 4952 1250 48 48  0.0000000 6.2831853}%
%
\special{pn 20}%
\special{pa 4584 898}%
\special{pa 4920 1226}%
\special{fp}%
\special{pa 4928 1226}%
\special{pa 4928 1226}%
\special{fp}%
%
\special{pn 20}%
\special{pa 4984 1290}%
\special{pa 5320 1618}%
\special{dt 0.054}%
\special{pa 5320 1618}%
\special{pa 5320 1618}%
\special{dt 0.054}%
%
\special{pn 20}%
\special{pa 5320 882}%
\special{pa 4984 1210}%
\special{fp}%
%
\special{pn 20}%
\special{pa 4912 1290}%
\special{pa 4608 1602}%
\special{dt 0.054}%
\put(43.6000,-8.3000){\makebox(0,0)[lb]{$\sc{1}$}}%
\put(54.4000,-8.3000){\makebox(0,0)[lb]{$\sc{2}$}}%
\put(43.4000,-18.2000){\makebox(0,0)[lb]{$\sc{3}$}}%
\put(54.0800,-18.4200){\makebox(0,0)[lb]{$\sc{4}$}}%
\put(48.7200,-11.5400){\makebox(0,0)[lb]{$\sc{0}$}}%
\put(7.0400,-20.0400){\makebox(0,0)[lb]{$D_4$}}%
\put(27.8400,-20.5200){\makebox(0,0)[lb]{$A_1^{\oplus 4}$}}%
\put(48.3200,-20.1000){\makebox(0,0)[lb]{$A_3$}}%
\put(2.9600,-23.3200){\makebox(0,0)[lb]{$I = \{0,1,2,3\}$}}%
\put(24.0800,-23.3200){\makebox(0,0)[lb]{$I = \{1,2,3,4\}$}}%
\put(44.7000,-23.4000){\makebox(0,0)[lb]{$I = \{0,1,2\}$}}%
\end{picture}%
\end{center}
\caption{Some $W(D_4^{(1)})$-strata and their abstract 
Dynkin types} \label{fig:strata}
\end{figure}

Corresponding to the two affine Weyl groups mentioned above, 
we can define two stratifications on $\K$. 
Let $\mathcal{I}$ be the set of all proper subsets 
$I \subset \{0,1,2,3,4\}$ including the empty set $\emptyset$. 
For each $I \in \mathcal{I}$ we denote by $\ol{\K}_I$ the 
$W(D_4^{(1)})$-translates of the affine subspace 
\[
H_I:=\bigcap_{i \in I} H_i, 
\]
and by $\K_I$ the set obtained from $\ol{\K}_I$ by removing all 
the sets $\ol{\K}_J$ such that $\# J = \# I +1$. 
Moreover let $D_I$ be the Dynkin subdiagram of $D_4^{(1)}$ that 
has nodes exactly in $I$. 
Some examples of these are given in Figure \ref{fig:strata}. 
It turns out that for any pair $I$, $I' \in \mathcal{I}$ 
either $\K_I=\K_{I'}$ or $\K_I \cap \K_{I'}=\emptyset$ holds 
(see Remark \ref{rem:stra}), so that one can define a 
stratification of $\K$ by the subsets $\K_I$ with $I \in \mathcal{I}$. 
It is called {\sl the $W(D_4^{(1)})$-stratification}. 
A parameter $\k \in \K$ is said to be {\sl generic} if 
$\k \in \K_{\emptyset}$; otherwise, $\k$ is said to be 
{\sl non-generic}.
\par 
For $I \in \mathcal{I}$ one can speak of the {\sl abstract} 
Dynkin type of the subdiagram $D_I$. 
For example, $D_I$ is of abstract type $A_3$ when $I = \{0,1,2\}$. 
All the possible abstract Dynkin types are those in Figure 
\ref{fig:AdRel} below. 
On the other hand there is a natural action of $S_4$ on the set 
$\mathcal{I}$ induced from its action on the nodes $\{1,2,3,4\}$ 
and the abstract Dynkin type of $D_I$ is represented by the 
$S_4$-orbit of $I$. 
Thus all the abstract Dynkin types are parametrized by the 
quotient set $\mathcal{I}/S_4$. 
\begin{remark}[\cite{Iwasaki1}] \label{rem:stra}
Let $I$ and $I'$ be distinct elements of $\mathcal{I}$. 
Then $\K_{I}=\K_{I'}$ if and only if $D_{I}$ and $D_{I'}$ have 
the same abstract type $A_1$ or $A_2$. 
Here the ``if" part is shown by a direct calculation, while 
the ``only if" part follows from the fact that the map 
(\ref{eqn:rhp}) is $W(D_4^{(1)})$-invariant and 
$\mathrm{rh}(\k) \neq \mathrm{rh}(\k')$ for any 
$\k \in \K_{I}$ and $\k' \in \K_{I'}$ if the condition is 
not fulfilled. 
Therefore, 
\begin{enumerate}
\item there is a unique $W(D_4^{(1)})$-stratum of abstract type 
$\emptyset$, $A_1$, $A_2$ or $A_1^{\oplus 4}$, 
\item there are six $W(D_4^{(1)})$-strata of abstract type 
$A_1^{\oplus 2}$ or $A_3$, 
\item there are four $W(D_4^{(1)})$-strata of abstract type 
$A_1^{\oplus 3}$ or $D_4$. 
\end{enumerate}
\end{remark}
\par 
We proceed to define a coarser stratification, that is, the 
{\sl $W(F_4^{(1)})$-stratification}. 
Observe that the $W(F_4^{(1)})$-translates of $H_I$ depends only 
on the abstract Dynkin type $* = [I] \in \mathcal{I}/S_4$, so 
that it is denoted by $\ol{\K}(*)$. 
Note that $\ol{\K}(*)$ is the $W(D_4^{(1)})$-translates of the union
\begin{equation} \label{eqn:F4H} 
H(*) := \bigcup_{[I] =*} H_I. 
\end{equation}
We say that $**$ is adjacent to $*$ and write $* \to **$ if 
$\ol{\K}(**)$ is a subset of $\ol{\K}(*)$. 
All the adjacency relations are depicted in Figure \ref{fig:AdRel}. 
\begin{figure}[t] 
\begin{center} 
\[
\begin{CD}
\emptyset @>>> A_1 @>>> A_1^{\oplus 2} @>>> A_1^{\oplus 3} @>>> 
A_1^{\oplus 4} \\
 & &   @VVV  @VVV  @VVV \\ 
 & &   A_2 @>>> A_3 @>>> D_4 
\end{CD}
\]
\end{center}
\caption{Adjacency relations among the $W(F_4^{(1)})$-strata}
\label{fig:AdRel} 
\end{figure}
The set $\K(*)$ is obtained from $\ol{\K}(*)$ by removing 
all the sets $\ol{\K}(**)$ such that $* \to **$. 
There is a direct sum decomposition 
\[
\K= \coprod_{* \in \mathcal{I}/S_4} \K(*). 
\]
Note that the $W(F_4^{(1)})$-stratum $\K(*)$ is the union of all 
$W(D_4^{(1)})$-strata of abstract type $*$. 
\par
We compactify the affine cubic surface $\Sol(\th)$ by the 
standard embedding: 
\[
\Sol(\th) \hookrightarrow \ol{\Sol}(\th) \subset \P^3, \qquad 
x = (x_1,x_2,x_3) \mapsto [1:x_1:x_2:x_3], 
\]
where the compactified surface is given by 
$\ol{\Sol}(\th) = \{\,X \in \P^3\,:\, F(X,\th)=0 \,\}$ with 
\[
F(X,\th) = X_1X_2X_3+X_0(X_1^2+X_2^2+X_3^2)-
X_0^2(\th_1X_1+\th_2X_2+\th_3X_3)+\th_4X_0^3. 
\]
The intersection of $\ol{\Sol}(\th)$ with the plane at infinity 
is the union $L$ of three lines 
\[
L_i = \{\, X \in \P^3 \,:\, X_0 = X_i = 0 \,\} \qquad 
(i = 1,2,3). 
\]
The set $L$, called the tritangent lines at infinity (see Figure 
\ref{fig:cubic1}), is independent of $\th \in \Th$ and the surface 
$\ol{\Sol}(\th)$ is smooth in a neighborhood of $L$ for every 
$\th \in \Th$ (see \cite[Lemma 2]{IU1}). 
\begin{figure}[t] 
\begin{center}
\unitlength 0.1in
\begin{picture}( 20.6800, 21.4400)( 14.1600,-31.5400)
%
\special{pn 13}%
\special{ar 2450 2120 1034 1034  4.9299431 6.2831853}%
\special{ar 2450 2120 1034 1034  0.0000000 4.4841152}%
%
\special{pn 20}%
\special{pa 1834 2400}%
\special{pa 3104 2400}%
\special{fp}%
%
\special{pn 20}%
\special{pa 2644 1420}%
\special{pa 1934 2610}%
\special{fp}%
%
\special{pn 20}%
\special{pa 2254 1420}%
\special{pa 2964 2600}%
\special{fp}%
\put(22.9000,-11.8000){\makebox(0,0)[lb]{$\ol{\Sol}(\th)$}}%
\put(29.6000,-18.5000){\makebox(0,0)[lb]{$\Sol(\th)$}}%
\put(23.7000,-26.1000){\makebox(0,0)[lb]{$L_i$}}%
\put(20.3000,-21.0000){\makebox(0,0)[lb]{$L_j$}}%
\put(27.1000,-20.9000){\makebox(0,0)[lb]{$L_k$}}%
\put(30.2000,-25.8000){\makebox(0,0)[lb]{$p_j$}}%
\put(23.9000,-15.9000){\makebox(0,0)[lb]{$p_i$}}%
\put(17.9000,-25.7000){\makebox(0,0)[lb]{$p_k$}}%
%
\special{pn 20}%
\special{sh 1.000}%
\special{ar 2060 2400 40 40  0.0000000 6.2831853}%
%
\special{pn 20}%
\special{sh 1.000}%
\special{ar 2450 1740 40 40  0.0000000 6.2831853}%
%
\special{pn 20}%
\special{sh 1.000}%
\special{ar 2840 2400 40 40  0.0000000 6.2831853}%
\end{picture}%
\end{center}
\caption{Tritangent lines at infinity on $\ol{\Sol}(\th)$} 
\label{fig:cubic1} 
\end{figure}
The intersection point of $L_j$ and $L_k$ is denoted by $p_i$ 
for $\{i,j,k\}=\{1,2,3\}$. 
It is explicitly given by 
\[
p_1 = [0:1:0:0], \qquad p_2 = [0:0:1:0], \qquad 
p_3 = [0:0:0:1]. 
\]
\par
When the parameter $\th = \rh(\k)$ is non-generic, based on a 
method in \cite{BW} we construct an algebraic minimal 
resolution of the singular surface $\ol{\Sol}(\th)$ by 
considering the rational map: 
\begin{equation} \label{eqn:birational}
\tau : \P^2 \longrightarrow \P^3, \quad u=[u_1:u_2:u_3] 
\longmapsto [ \tau_0 (u) : \tau_1 (u) : \tau_2 (u) : 
\tau_3 (u) ],
\end{equation}
where the polynomials $\tau_0 (u)$, $\tau_1(u)$, $\tau_2(u)$, 
$\tau_3(u)$ are given by
\[
\left\{ 
\begin{array}{rcl}
\tau_0(u) &:=& -b_0^2 u_1 u_2 u_3, \\[2mm]
\tau_1(u) &:=& b_0^2 u_1 \{ b_0^2 u_1^2 + u_2^2 + u_3^2 + 
b_0^2 (b_1 b_2 + b_3 b_4) u_1 u_2 
+ b_0^2 (b_1 b_3 + b_2 b_4) u_1 u_3 \}, \\[2mm]
\tau_2(u) &:=& u_2 \{ b_0^4 u_1^2 + b_0^2 u_2^2 + u_3^2 
+ b_0^4 (b_1 b_2 + b_3 b_4) u_1 u_2 
+ b_0^2 (b_2 b_3 + b_1 b_4) u_2 u_3 \}, \\[2mm]
\tau_3(u) &:=& u_3\{ b_0^2 u_1^2 + b_0^2 u_2^2 + u_3^2 
+ b_0^2 (b_2 b_3 + b_1 b_4) u_2 u_3 
+ b_0^2 (b_1 b_3 + b_2 b_4) u_1 u_3 \}. \\[2mm]
\end{array}\right. 
\]
It is a birational map of $\P^2$ onto $\ol{\Sol}(\th)$ whose 
indeterminacy points are the six points
\begin{equation} \label{eqn:indeterminacy}
\begin{array}{rclrcl}
c_1 &:=& [0:- b_1 b_4:1], \qquad & c_4 &:=& [0:- b_2 b_3:1] \in l_1, \\[2mm]
c_2 &:=& [- b_1 b_3:0:1], \qquad & c_5 &:=& [- b_2 b_4:0:1] \in l_2, \\[2mm]
c_3 &:=& [- b_3 b_4:1:0], \qquad & c_6 &:=& [- b_1 b_2:1:0] \in l_3, 
\end{array}
\end{equation}
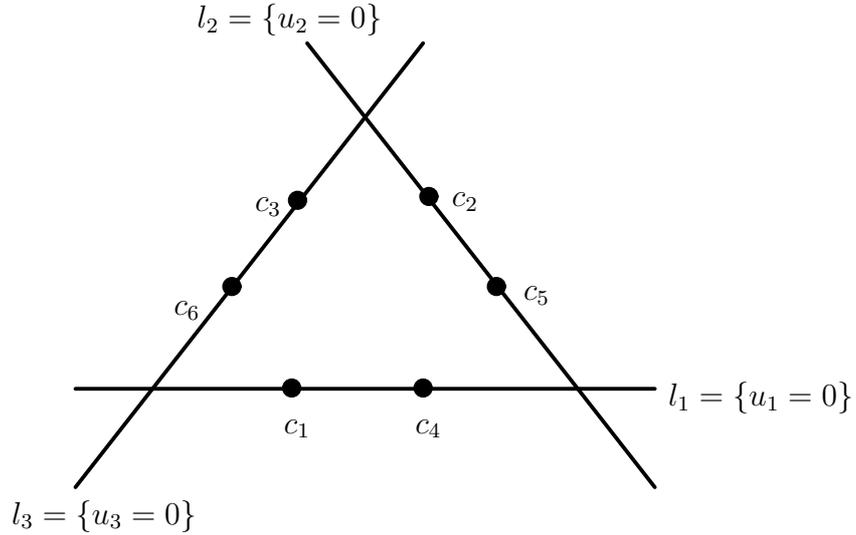
\begin{figure}[t] 
\begin{center}
\unitlength 0.1in
\begin{picture}( 34.0000, 25.8800)(  2.7000,-28.1000)
\put(16.8000,-24.9000){\makebox(0,0)[lb]{$c_1$}}%
\put(25.5000,-13.0600){\makebox(0,0)[lb]{$c_2$}}%
\put(15.3000,-13.3200){\makebox(0,0)[lb]{$c_3$}}%
\put(23.6000,-24.9000){\makebox(0,0)[lb]{$c_4$}}%
\put(29.2000,-17.9500){\makebox(0,0)[lb]{$c_5$}}%
\put(11.1000,-18.7300){\makebox(0,0)[lb]{$c_6$}}%
\put(36.7000,-23.6200){\makebox(0,0)[lb]{$l_1=\{ u_1 = 0 \}$}}%
\put(12.3000,-3.9200){\makebox(0,0)[lb]{$l_2=\{ u_2 = 0 \}$}}%
\put(2.7000,-29.8000){\makebox(0,0)[lb]{$l_3=\{ u_3 = 0 \}$}}%
%
\special{pn 20}%
\special{sh 1.000}%
\special{ar 2430 1230 40 40  0.0000000 6.2831853}%
%
\special{pn 20}%
\special{sh 1.000}%
\special{ar 2780 1700 40 40  0.0000000 6.2831853}%
%
\special{pn 20}%
\special{sh 1.000}%
\special{ar 1750 1250 40 40  0.0000000 6.2831853}%
%
\special{pn 20}%
\special{sh 1.000}%
\special{ar 1410 1700 40 40  0.0000000 6.2831853}%
%
\special{pn 20}%
\special{sh 1.000}%
\special{ar 1720 2230 40 40  0.0000000 6.2831853}%
%
\special{pn 20}%
\special{sh 1.000}%
\special{ar 2400 2230 40 40  0.0000000 6.2831853}%
%
\special{pn 20}%
\special{pa 600 2234}%
\special{pa 3600 2234}%
\special{fp}%
\special{pa 2400 430}%
\special{pa 600 2748}%
\special{fp}%
\special{pa 1800 430}%
\special{pa 3600 2748}%
\special{fp}%
\end{picture}%
\end{center}
\caption{Six indeterminacy points of $\tau$ in $\P^2$} 
\label{fig:indeterminacy}
\end{figure} 
where $l_i$ is the strict transform of the line $L_i$ under the 
map (\ref{eqn:birational}) and is given by 
\[
l_i = \{[u_1:u_2:u_3] \in \P^2 \, ; \, u_i=0 \} \qquad 
(i \in \{ 1,2,3 \}). 
\]
\par
Let $\rho : \widetilde{\Sol} (\th) \to \P^2$ be the blow-up 
of $\P^2$ at the six points $c_1,\dots,c_6$, and put 
\begin{equation} \label{eqn:des}
\pi := \rho \circ \tau : \widetilde{\Sol} (\th) \to \ol{\Sol}(\th).
\end{equation}
\begin{proposition} \label{prop:desing}
The birational morphism $(\ref{eqn:des})$ gives an algebraic 
minimal resolution of $\ol{\Sol}(\th)$. 
\end{proposition} 
{\it Proof}. 
Starting with $\Sol_0 := \P^2$ we inductively define 
$\rho_i : \Sol_i \to \Sol_{i-1}$ as the blow-up of 
$\Sol_{i-1}$ at the point $c_i$ for $i = 1,\dots,6$. 
The birational morphism $\rho : \widetilde{\Sol}(\th) \to 
\P^2$ then decomposes as 
\[ 
\rho = \rho_1 \circ \cdots \circ \rho_6. 
\]
Consider the linear system 
$\delta := \{ \, p_0 \tau_0 + p_1 \tau_1 + p_2 \tau_2 + 
p_3 \tau_3 \, : \, [p_0:p_1:p_2:p_3] \in \P^3 \, \} 
\subset |3H|$, where $H$ is a line in $\P^2$. 
Starting with $\delta_0 := \delta$ and $D_0 := 3H$ we 
inductively define 
$\delta_i := \rho_i^* \delta_{i-1} - E_i$ and 
$D_i := \rho_i^* D_{i-1} - E_i$ for $i = 1,\dots,6$, where 
$E_i := \rho_i^{*} (c_{i-1})$ is the exceptional curve of 
$\rho_i$ over the point $c_{i-1}$. 
In view of (\ref{eqn:birational}) the linear system 
$\rho_i^* \delta_{i-1} \subset|\rho_i^* D_{i-1}|$ admits 
$E_i$ as its fixed part and thus one has the inclusion 
$\delta_i \subset |D_i|$. 
Since the canonical divisor $K_{\Sol_0} = K_{\P^2}$ on 
$\P^2$ is linearly equivalent to the divisor 
$-3 H = -D_0$, we have 
\[ 
K_{\Sol_i} = \rho_i^* K_{\Sol_{i-1}} + E_i 
\sim -(\rho_i^* D_{i-1} - E_i) = - D_i \qquad 
(i = 1,\dots,6). 
\]
In particular the canonical divisor 
$K_{\widetilde{\Sol}(\th)} = K_{\Sol_6}$ on 
$\widetilde{\Sol}(\th)$ is linearly equivalent to $- D_{6}$. 
Therefore for any $(-1)$-curve $C$ on $\widetilde{\Sol}(\th)$ 
we have $C \cdot D_6 = -(C \cdot K_{\widetilde{\Sol}(\th)}) = 1$ 
from the adjunction formula. 
This means that the birational morphism (\ref{eqn:des}) 
never contracts any $(-1)$-curve $C$ to a point in $\ol{\Sol}(\th)$. 
Thus the proposition is established. \hfill $\Box$
\begin{remark} \label{rem:DO}
The surface $\ol{\Sol}(\th)$, if it is smooth, is known as a 
{\sl Del Pezzo surface} of degree $3$. 
On the other hand, if $\ol{\Sol}(\th)$ is singular, the 
desingularized surface $\wt{\Sol}(\th)$ is called a 
{\sl degenerate Del Pezzo} surface of degree $3$ in \cite{DO}. 
Any degenerate Del Pezzo surface of degree $3$ is obtained 
as a blow-up of $\P^2$ at six points, no two of which are 
infinitely near of order $1$ to the same point and no four 
of which are collinear. 
The image of a degenerate Del Pezzo surface under the map 
of its anti-canonical class is called an {\sl anti-canonical} 
Del Pezzo surface of degree $3$, an example of which is our 
cubic surface $\ol{\Sol}(\th)$. 
\end{remark}
\par
Since $\widetilde{\Sol}(\th)$ is a six-point blow-up of $\P^2$, 
its second cohomology group is expressed as 
\[
H^2(\widetilde{\Sol}(\th),\Z) = 
\Z E_0 \oplus \Z E_1 \oplus \Z E_2 \oplus \Z E_3 \oplus 
\Z E_4 \oplus \Z E_5 \oplus \Z E_6, 
\]
where $E_0$ is the class of the strict transform of a line in 
$\P^2$ and $E_i$ is the class of the exceptional curve over 
the point $c_i$ for $i = 1,2,3,4,5,6$. 
Their intersection relations are listed as 
\begin{equation} \label{eqn:in1}
(E_i, E_j) = \left\{ 
\begin{array}{ll}
{\mi}1 \quad & (i = j = 0), \\[1mm] 
-1 \quad & (i = j \neq 0), \\[1mm]
{\mi}0 \quad & (\text{otherwise}).
\end{array}\right. 
\end{equation}
If the line $L_i$ is identified with the cohomology class 
of its strict transform under $\pi$, we have 
\begin{equation} \label{eqn:infline}
L_1 = E_0-E_1-E_4, \quad 
L_2 = E_0-E_2-E_5, \quad
L_3 = E_0-E_3-E_6. 
\end{equation}
In view of formulas (\ref{eqn:in1}) and (\ref{eqn:infline}) 
there exists a direct sum decomposition 
\begin{equation} \label{eqn:decomp}
H^2(\widetilde{\Sol}(\th),\C) = V \oplus V^{\perp}, 
\end{equation} 
where $V$ is the subspace spanned by $L_1$, $L_2$, $L_3$ and 
$V^{\perp}$ is the orthogonal complement to $V$ with respect 
to the intersection form. 
It is easy to see that $V^{\perp}$ is spanned by the vectors 
\begin{equation} \label{eqn:bas}
L_4 := E_1-E_4, \quad
L_5 := E_1-E_5, \quad
L_6 := E_1-E_6, \quad
L_7 := 2 E_0-E_1-E_2-E_3-E_4-E_5-E_6.
\end{equation}
The decomposition (\ref{eqn:decomp}) will be important 
in the ergodic study of birational maps on 
$\widetilde{\Sol}(\th)$. 
\par
Now we investigate the structure of singularities on 
$\ol{\Sol}(\th)$ utilizing the rational map 
(\ref{eqn:birational}), especially observing the 
configuration of its indeterminacy points. 
For a case-by-case treatment we divide the parameter space 
$B$ into several pieces $\ol{B}(*)$ depending on abstract 
Dynkin types $* \in \mathcal{I}/S_4$, each of which is 
further decomposed into smaller pieces $\ol{B}_m(*)$ and 
then into even smaller pieces $B_m(*;\mbox{-})$. 
Here the definitions of $\ol{B}_m(*)$ and 
$\ol{B}_m(*;\mbox{-})$ are given in Table \ref{tab:type}. 
\begin{table}[p]
\begin{center} 
\begin{tabular}{|l||l|}
\hline 
\vspace{-4mm} & \\
parameter spaces $\ol{B}_m(*)$ & 
defining equations of $\ol{B}_m(*;\mbox{-})$
\\[1mm] \hline \hline 
\vspace{-4mm} & \\
$\ol{B}_1(D_4) = \ol{B}_1(D_4;\cdot)$ & 
$\ve_1 b_1 = \ve_2 b_2 = \ve_3 b_3 = \ve_4 b_4 = 1$ 
\\[1mm] \hline \hline 
\vspace{-4mm} & \\
$\ol{B}_1(A_1^{\oplus 4}) = \ol{B}_1(A_1^{\oplus 4};\cdot)$ & 
$\ve_1 b_1 = \ve_2 b_2 = \ve_3 b_3 = - \ve_4 b_4 \in \{1, \sqrt{-1}\}$ 
\\[1mm]\hline 
\vspace{-4mm} & \\
$\ol{B}_2(A_1^{\oplus 4}) = \ol{B}_2(A_1^{\oplus 4};\cdot)$ & 
$\ve_1 b_1 = \ve_2 b_2 = \ve_3 b_3 = \ve_4 b_4 = \sqrt{-1} $
\\[1mm] \hline \hline 
\vspace{-4mm} & \\
$\ol{B}_1(A_3) = \ds \bigcup_{1 \le i <j \le 4} \ol{B}_1(A_3;i,j)$ & 
$\ve_i b_i = \ve_j b_j = 1, \, \ve_k b_k = \ve_l b_l$
\\[5mm] \hline 
\vspace{-4mm} & \\
$\ol{B}_2(A_3) = \ds \bigcup_{1 \le i <j \le 4} \ol{B}_2(A_3;i,j)$ &
$ \ve_i b_i = \ve_j b_j = 1, \, \ve_k b_k = (\ve_l b_l)^{-1}$
\\[5mm] \hline \hline 
\vspace{-4mm} & \\
$\ol{B}_1(A_1^{\oplus 3}) =  \ds \bigcup_{1 \le i \le 4} 
\ol{B}_1(A_1^{\oplus 3};i)$ &
$\ve_j b_j = \ve_k b_k = \ve_l b_l = 1 $
\\[5mm] \hline 
\vspace{-4mm} & \\
$\ol{B}_2(A_1^{\oplus 3}) =  \ds \bigcup_{1 \le i \le 4} 
\ol{B}_2(A_1^{\oplus 3};i)$ &
$\ve_j b_j = \ve_k b_k = \ve_l b_l = (\ve_i b_i)^{-1}$ 
\\[5mm] \hline 
\vspace{-4mm} & \\
$\ol{B}_3(A_1^{\oplus 3}) = \ol{B}_3(A_1^{\oplus 3};\cdot)$ &
$\ve_1 b_1 = \ve_2 b_2 = \ve_3 b_3 = \ve_4 b_4$
\\[1mm] \hline 
\vspace{-4mm} & \\
$\ol{B}_4(A_1^{\oplus 3}) = \ds \bigcup_{1 \le i \le 3} 
\ol{B}_4(A_1^{\oplus 3};i)$ &
$ \ve_j b_j = \ve_k b_k = (\ve_i b_i)^{-1} = (\ve_4 b_4)^{-1}$
\\[5mm] \hline \hline 
\vspace{-4mm} & \\
$\ol{B}_1(A_2) = \ds \bigcup_{1 \le i \neq j \le 4} \ol{B}_1(A_2;i,j)$ &
$\ve_i b_i = 1, \, \ve_k b_k = b_l b_l$
\\[5mm] \hline 
\vspace{-4mm} & \\
$\ol{B}_2(A_2) = \ds \bigcup_{1 \le i \le 4} \ol{B}_2(A_2;i)$ &
$ \ve_i b_i = \ve_i b_j b_k b_l= 1$
\\[5mm] \hline \hline 
\vspace{-4mm} & \\
$\ol{B}_1(A_1^{\oplus 2}) = \ds \bigcup_{1 \le i < j \le 4} 
\ol{B}_1(A_1^{\oplus 2};i,j)$ &
$\ve_i b_i = \ve_j b_j = 1$ 
\\[5mm] \hline 
\vspace{-4mm} & \\
$\ol{B}_2(A_1^{\oplus 2}) = \ds \bigcup_{1 \le i < j \le 4} 
\ol{B}_2(A_1^{\oplus 2};i,j)$ &
$b_i b_j^{-1} =b_k b_l = \ve_i$ 
\\[5mm] \hline 
\vspace{-4mm} & \\
$\ol{B}_3(A_1^{\oplus 2}) = \ds \bigcup_{1 \le i \le 3} 
\ol{B}_3(A_1^{\oplus 2};i)$ &
$b_i b_4^{-1}= b_j b_k^{-1} = \ve_i$
\\[5mm] \hline 
\vspace{-4mm} & \\
$\ol{B}_4(A_1^{\oplus 2}) = \ds \bigcup_{1 \le i \le 3} 
\ol{B}_4(A_1^{\oplus 2};i)$ &
$b_i b_4= b_j b_k = \ve_i$
\\[5mm] \hline \hline 
\vspace{-4mm} & \\
$\ol{B}_1(A_1) = \ds \bigcup_{1 \le i \le 4} \ol{B}_1(A_1;i)$ &
$\ve_i b_i = 1$ 
\\[5mm] \hline 
\vspace{-4mm} & \\
$\ol{B}_2(A_1) = \ds \bigcup_{1 \le i \le 4} \ol{B}_2(A_1;i)$ &
$b_i = b_j b_k b_l$ 
\\[5mm] \hline 
\vspace{-4mm} & \\
$\ol{B}_3(A_1) = \ds \bigcup_{1 \le i \le 3} \ol{B}_3(A_1;i)$ &
$b_i b_4 = b_j b_k$
\\[5mm] \hline 
\vspace{-4mm} & \\
$\ol{B}_4(A_1) = \ol{B}_4(A_1;\cdot)$ &
$b_1 b_2 b_3 b_4 = 1 $
\\[1mm] \hline 
\end{tabular}
\end{center}
\caption{Parameter spaces $\ol{B}_m(*)$: $\ve_i \in\{\pm 1\}$ 
satisfy $\ve_1 \ve_2 \ve_3 \ve_4=1$.} 
\label{tab:type}
\end{table}
Put
\[
B(*) := \ol{B}(*) \setminus \bigcup_{* \to **} \ol{B}(**). 
\]
\begin{proposition} \label{prop:type}
Given any $b \in B$, put $\th = (\varphi \circ \alpha)(b)$. 
Then $\ol{\Sol}(\th)$ has simple singularities of abstract 
Dynkin type $* \in \mathcal{I}/S_4$ if and only if the parameter 
$b$ is an element of $B(*) := \cup_{m} B_m(*)$. 
\end{proposition}
{\it Proof}. 
A careful inspection of formulas (\ref{eqn:vD}) and 
(\ref{eqn:indeterminacy}) shows that the surface $\ol{\Sol}(\th)$ 
is singular if and only if the indeterminacy points 
$c_1,\dots,c_6 \in \P^2$ of the rational map (\ref{eqn:birational}) 
are {\sl not} in a general position, namely, if and only if 
at least one of the following conditions are satisfied. 
\begin{enumerate}
\item[(C1)] The six points lie on a (unique) conic; 
this condition is equivalent 
to $b_1 b_2 b_3 b_4 = 1$. 
\item[(C2)] Three of them, say, $c_i$, $c_j$ and $c_k$ are colinear; 
this condition is equivalent to 
\[
\left\{
\begin{array}{rcll}
b_4 - b_4^{-1} &=& 0 \quad & (\{i,j,k\}=\{1,2,3\}), \\[2mm]
b_1 b_2 b_3 b_4^{-1} &=& 1 \quad & (\{i,j,k\}=\{4,5,6\}), \\[2mm]
b_i - b_i^{-1} &=& 0 \quad & (i \in \{1,2,3\}, \, 
\{j,k\}=\{4,5,6\} \setminus \{i+3\}), \\[2mm]
b_i^{-1} b_j b_k b_4 &=& 1 \quad & (i \in \{4,5,6\}, \, 
\{j,k\}=\{1,2,3\} \setminus \{i-3\}). 
\end{array}
\right. 
\]
\item[(C3)] $c_i = c_{i+3}$ for some $i \in \{1,2,3\}$, 
that is, $b_i b_j^{-1} b_k^{-1} b_4 =1$ with 
$\{i,j,k\} = \{1,2,3\}$. 
\end{enumerate}
\par
The Dynkin graph of the singularities on $\ol{\Sol}(\th)$ 
appears as the dual graph of the $(-2)$-curves on 
the desingularized surface $\wt{\Sol}(\th)$. 
Here each $(-2)$-curve arises either as the strict transform 
of the conic in case (C1), or as the strict transform of 
the line through $c_i$, $c_j$, $c_k$ in case (C2), 
or as the exceptional curve over the 
``degenerate'' point $c_i$ in case (C3). 
Thus the Dynkin structure of the singularities on 
$\ol{\Sol}(\th)$ can be read off from the following data: 
\begin{itemize}
\item the six points themselves; 
\item the unique conic, if it exists, specified by condition (C1); 
\item all the lines specified by condition (C2);  
\item all the degenerate points specified by condition (C3). 
\end{itemize}
\par
A case-by-case check shows that all feasible data on 
varous parameter spaces $B_m(*)$ are depicted in 
Figures \ref{fig:A_14D_4}--\ref{fig:A_2}, where a degenerate 
indeterminacy point is marked by a white circle and a 
nondegenerate one is by a black-filled circle respectively. 
Each $(-2)$-curve arises either from the conic, if it exists, 
or from a line, or from a white circle in the figures. 
The $(-2)$-curve from a white circle intersects the one from 
a conic but none from a line. 
If two lines intersect in a black-filled circle, then the 
corresponding $(-2)$-curves are disjoint; otherwise 
they do meet in a single point. 
In this manner the data determines the dual graph of the 
$(-2)$-curves on $\wt{\Sol}(\th)$ and hence the Dynkin type 
of the singularities on $\ol{\Sol}(\th)$. \hfill $\Box$
\begin{figure}[!t] 
\begin{center}
\unitlength 0.1in
\begin{picture}( 60.8300, 14.1000)(  2.6000,-19.8000)
%
\special{pn 20}%
\special{pa 4424 1790}%
\special{pa 4440 1744}%
\special{pa 4456 1696}%
\special{pa 4470 1650}%
\special{pa 4486 1604}%
\special{pa 4502 1558}%
\special{pa 4518 1512}%
\special{pa 4534 1466}%
\special{pa 4550 1420}%
\special{pa 4566 1376}%
\special{pa 4582 1332}%
\special{pa 4598 1288}%
\special{pa 4614 1246}%
\special{pa 4630 1204}%
\special{pa 4644 1162}%
\special{pa 4660 1122}%
\special{pa 4676 1084}%
\special{pa 4692 1046}%
\special{pa 4708 1008}%
\special{pa 4724 972}%
\special{pa 4740 938}%
\special{pa 4756 904}%
\special{pa 4772 872}%
\special{pa 4788 842}%
\special{pa 4804 812}%
\special{pa 4820 784}%
\special{pa 4836 758}%
\special{pa 4852 734}%
\special{pa 4868 712}%
\special{pa 4884 692}%
\special{pa 4900 672}%
\special{pa 4916 656}%
\special{pa 4932 640}%
\special{pa 4948 628}%
\special{pa 4964 616}%
\special{pa 4980 608}%
\special{pa 4996 602}%
\special{pa 5012 598}%
\special{pa 5028 596}%
\special{pa 5044 598}%
\special{pa 5060 600}%
\special{pa 5076 606}%
\special{pa 5092 614}%
\special{pa 5108 626}%
\special{pa 5124 640}%
\special{pa 5140 656}%
\special{pa 5158 674}%
\special{pa 5174 696}%
\special{pa 5190 718}%
\special{pa 5206 744}%
\special{pa 5222 770}%
\special{pa 5238 798}%
\special{pa 5254 828}%
\special{pa 5272 860}%
\special{pa 5288 892}%
\special{pa 5304 926}%
\special{pa 5320 962}%
\special{pa 5336 996}%
\special{pa 5352 1034}%
\special{pa 5370 1070}%
\special{pa 5386 1108}%
\special{pa 5402 1146}%
\special{pa 5418 1184}%
\special{pa 5434 1220}%
\special{pa 5452 1258}%
\special{pa 5468 1296}%
\special{pa 5484 1334}%
\special{pa 5500 1370}%
\special{pa 5516 1406}%
\special{pa 5532 1442}%
\special{pa 5550 1476}%
\special{pa 5566 1510}%
\special{pa 5582 1542}%
\special{pa 5598 1572}%
\special{pa 5614 1602}%
\special{pa 5630 1630}%
\special{pa 5648 1656}%
\special{pa 5664 1680}%
\special{pa 5680 1702}%
\special{pa 5696 1722}%
\special{pa 5712 1740}%
\special{pa 5728 1756}%
\special{pa 5744 1768}%
\special{pa 5762 1778}%
\special{pa 5778 1786}%
\special{pa 5794 1790}%
\special{pa 5810 1792}%
\special{pa 5826 1790}%
\special{pa 5842 1786}%
\special{pa 5858 1780}%
\special{pa 5874 1770}%
\special{pa 5890 1758}%
\special{pa 5906 1744}%
\special{pa 5922 1726}%
\special{pa 5938 1708}%
\special{pa 5954 1686}%
\special{pa 5968 1664}%
\special{pa 5984 1638}%
\special{pa 6000 1612}%
\special{pa 6014 1584}%
\special{pa 6030 1554}%
\special{pa 6044 1522}%
\special{pa 6058 1490}%
\special{pa 6074 1456}%
\special{pa 6088 1422}%
\special{pa 6102 1386}%
\special{pa 6116 1348}%
\special{pa 6130 1312}%
\special{pa 6142 1274}%
\special{pa 6156 1234}%
\special{pa 6170 1196}%
\special{pa 6182 1156}%
\special{pa 6194 1116}%
\special{pa 6206 1076}%
\special{pa 6218 1036}%
\special{pa 6230 996}%
\special{pa 6242 956}%
\special{pa 6252 916}%
\special{pa 6264 878}%
\special{pa 6274 840}%
\special{pa 6284 802}%
\special{pa 6294 764}%
\special{pa 6304 728}%
\special{pa 6312 692}%
\special{pa 6322 658}%
\special{pa 6330 624}%
\special{pa 6338 594}%
\special{pa 6344 570}%
\special{sp}%
\put(53.6600,-10.2700){\makebox(0,0)[lb]{$c_2=c_5$}}%
\put(60.3400,-17.3300){\makebox(0,0)[lb]{$c_3=c_6$}}%
%
\special{pn 20}%
\special{pa 260 1196}%
\special{pa 1918 1196}%
\special{fp}%
\put(8.3100,-14.1500){\makebox(0,0)[lb]{$c_2=c_5$}}%
\put(3.2400,-11.2600){\makebox(0,0)[lb]{$c_1=c_4$}}%
\put(13.1200,-11.2600){\makebox(0,0)[lb]{$c_3=c_6$}}%
%
\special{pn 20}%
\special{sh 0}%
\special{ar 6044 1504 46 50  0.0000000 6.2831853}%
%
\special{pn 20}%
\special{sh 0}%
\special{ar 550 1196 46 50  0.0000000 6.2831853}%
%
\special{pn 20}%
\special{sh 0}%
\special{ar 1058 1196 46 50  0.0000000 6.2831853}%
%
\special{pn 20}%
\special{sh 0}%
\special{ar 1538 1196 46 50  0.0000000 6.2831853}%
\put(7.5900,-21.5000){\makebox(0,0)[lb]{$B(D_4)$}}%
\put(29.9000,-21.5000){\makebox(0,0)[lb]{$B_1(A_1^{\oplus 4})$}}%
\put(52.3900,-21.5000){\makebox(0,0)[lb]{$B_2(A_1^{\oplus 4})$}}%
%
\special{pn 20}%
\special{sh 0}%
\special{ar 4650 1146 46 50  0.0000000 6.2831853}%
\put(46.2200,-14.1500){\makebox(0,0)[lb]{$c_1=c_4$}}%
%
\special{pn 20}%
\special{sh 0}%
\special{ar 5394 1136 46 50  0.0000000 6.2831853}%
%
\special{pn 20}%
\special{pa 3072 600}%
\special{pa 3970 1792}%
\special{fp}%
\special{pa 3426 600}%
\special{pa 2520 1792}%
\special{fp}%
%
\special{pn 20}%
\special{pa 2528 998}%
\special{pa 3970 1594}%
\special{fp}%
\special{pa 3970 998}%
\special{pa 2520 1594}%
\special{fp}%
%
\special{pn 20}%
\special{sh 1.000}%
\special{ar 3246 818 46 50  0.0000000 6.2831853}%
%
\special{pn 20}%
\special{sh 1.000}%
\special{ar 2982 1176 46 50  0.0000000 6.2831853}%
%
\special{pn 20}%
\special{sh 1.000}%
\special{ar 2738 1504 46 50  0.0000000 6.2831853}%
%
\special{pn 20}%
\special{sh 1.000}%
\special{ar 3246 1286 46 50  0.0000000 6.2831853}%
%
\special{pn 20}%
\special{sh 1.000}%
\special{ar 3518 1176 44 50  0.0000000 6.2831853}%
%
\special{pn 20}%
\special{sh 1.000}%
\special{ar 3762 1504 44 50  0.0000000 6.2831853}%
\end{picture}%
\end{center}
\caption{On the strata of types $D_4$ and $A_1^{\oplus 4}$} 
\label{fig:A_14D_4} 
\vspace{2mm}
\begin{center}
\unitlength 0.1in
\begin{picture}( 61.6400, 14.4200)(  2.6000,-20.0000)
%
\special{pn 20}%
\special{pa 4516 1794}%
\special{pa 4532 1746}%
\special{pa 4548 1700}%
\special{pa 4564 1652}%
\special{pa 4580 1606}%
\special{pa 4596 1560}%
\special{pa 4612 1514}%
\special{pa 4628 1468}%
\special{pa 4642 1424}%
\special{pa 4658 1378}%
\special{pa 4674 1334}%
\special{pa 4690 1292}%
\special{pa 4706 1248}%
\special{pa 4722 1208}%
\special{pa 4738 1166}%
\special{pa 4754 1126}%
\special{pa 4770 1086}%
\special{pa 4786 1048}%
\special{pa 4802 1012}%
\special{pa 4818 976}%
\special{pa 4834 942}%
\special{pa 4850 908}%
\special{pa 4864 876}%
\special{pa 4880 846}%
\special{pa 4896 816}%
\special{pa 4912 790}%
\special{pa 4928 764}%
\special{pa 4944 740}%
\special{pa 4960 718}%
\special{pa 4976 696}%
\special{pa 4992 678}%
\special{pa 5008 662}%
\special{pa 5024 646}%
\special{pa 5040 634}%
\special{pa 5056 624}%
\special{pa 5072 616}%
\special{pa 5088 610}%
\special{pa 5104 606}%
\special{pa 5122 606}%
\special{pa 5138 606}%
\special{pa 5154 610}%
\special{pa 5170 616}%
\special{pa 5186 626}%
\special{pa 5202 638}%
\special{pa 5218 652}%
\special{pa 5234 670}%
\special{pa 5250 688}%
\special{pa 5266 710}%
\special{pa 5284 734}%
\special{pa 5300 758}%
\special{pa 5316 786}%
\special{pa 5332 814}%
\special{pa 5348 846}%
\special{pa 5364 876}%
\special{pa 5380 910}%
\special{pa 5398 944}%
\special{pa 5414 978}%
\special{pa 5430 1014}%
\special{pa 5446 1052}%
\special{pa 5462 1088}%
\special{pa 5480 1126}%
\special{pa 5496 1164}%
\special{pa 5512 1202}%
\special{pa 5528 1240}%
\special{pa 5544 1278}%
\special{pa 5560 1314}%
\special{pa 5578 1352}%
\special{pa 5594 1388}%
\special{pa 5610 1424}%
\special{pa 5626 1460}%
\special{pa 5642 1494}%
\special{pa 5658 1526}%
\special{pa 5676 1558}%
\special{pa 5692 1588}%
\special{pa 5708 1616}%
\special{pa 5724 1644}%
\special{pa 5740 1670}%
\special{pa 5756 1694}%
\special{pa 5772 1714}%
\special{pa 5790 1734}%
\special{pa 5806 1750}%
\special{pa 5822 1764}%
\special{pa 5838 1776}%
\special{pa 5854 1786}%
\special{pa 5870 1792}%
\special{pa 5886 1794}%
\special{pa 5902 1796}%
\special{pa 5918 1792}%
\special{pa 5934 1786}%
\special{pa 5950 1778}%
\special{pa 5966 1768}%
\special{pa 5982 1754}%
\special{pa 5998 1738}%
\special{pa 6012 1720}%
\special{pa 6028 1700}%
\special{pa 6044 1678}%
\special{pa 6060 1654}%
\special{pa 6074 1628}%
\special{pa 6090 1602}%
\special{pa 6104 1572}%
\special{pa 6120 1542}%
\special{pa 6134 1510}%
\special{pa 6148 1476}%
\special{pa 6162 1442}%
\special{pa 6176 1406}%
\special{pa 6190 1370}%
\special{pa 6204 1332}%
\special{pa 6218 1294}%
\special{pa 6230 1256}%
\special{pa 6244 1216}%
\special{pa 6256 1178}%
\special{pa 6270 1138}%
\special{pa 6282 1098}%
\special{pa 6294 1058}%
\special{pa 6306 1018}%
\special{pa 6318 978}%
\special{pa 6328 938}%
\special{pa 6340 900}%
\special{pa 6350 860}%
\special{pa 6360 822}%
\special{pa 6370 784}%
\special{pa 6380 748}%
\special{pa 6390 712}%
\special{pa 6400 678}%
\special{pa 6408 644}%
\special{pa 6416 612}%
\special{pa 6424 582}%
\special{pa 6424 580}%
\special{sp}%
\put(54.5400,-10.3400){\makebox(0,0)[lb]{$c_j=c_{j+3}$}}%
\put(61.1800,-17.3500){\makebox(0,0)[lb]{$c_{k+3}$}}%
\put(63.0700,-11.7200){\makebox(0,0)[lb]{$c_k$}}%
%
\special{pn 20}%
\special{pa 270 828}%
\special{pa 1918 828}%
\special{fp}%
\put(5.1200,-14.9800){\makebox(0,0)[lb]{$c_{i+3}$}}%
\put(9.9900,-13.0100){\makebox(0,0)[lb]{$c_{j+3}$}}%
\put(10.1700,-7.2800){\makebox(0,0)[lb]{$c_j$}}%
\put(4.7600,-7.2800){\makebox(0,0)[lb]{$c_i$}}%
\put(14.5900,-14.7900){\makebox(0,0)[lb]{$c_{k+3}$}}%
\put(14.8600,-7.3800){\makebox(0,0)[lb]{$c_k$}}%
%
\special{pn 20}%
\special{ar 2622 1242 46 50  0.0000000 6.2831853}%
\put(24.0600,-11.7200){\makebox(0,0)[lb]{$c_1=c_3$}}%
\put(29.1100,-14.8800){\makebox(0,0)[lb]{$c_2=c_5$}}%
\put(34.4300,-13.6000){\makebox(0,0)[lb]{$c_3=c_6$}}%
%
\special{pn 20}%
\special{sh 1.000}%
\special{ar 6334 936 46 50  0.0000000 6.2831853}%
%
\special{pn 20}%
\special{sh 1.000}%
\special{ar 6128 1508 46 50  0.0000000 6.2831853}%
%
\special{pn 20}%
\special{sh 1.000}%
\special{ar 558 828 44 50  0.0000000 6.2831853}%
%
\special{pn 20}%
\special{sh 1.000}%
\special{ar 1064 828 44 50  0.0000000 6.2831853}%
%
\special{pn 20}%
\special{sh 1.000}%
\special{ar 1540 828 46 50  0.0000000 6.2831853}%
%
\special{pn 20}%
\special{sh 1.000}%
\special{ar 558 1272 44 50  0.0000000 6.2831853}%
%
\special{pn 20}%
\special{sh 1.000}%
\special{ar 1064 1044 44 50  0.0000000 6.2831853}%
%
\special{pn 20}%
\special{sh 1.000}%
\special{ar 1540 1242 46 50  0.0000000 6.2831853}%
\put(4.0000,-21.5000){\makebox(0,0)[lb]{$B_1(A_1^{\oplus 3})$ or $B_2(A_1^{\oplus 3})$}}%
\put(30.2700,-21.6000){\makebox(0,0)[lb]{$B_3(A_1^{\oplus 3})$}}%
\put(52.2000,-21.7000){\makebox(0,0)[lb]{$B_4(A_1^{\oplus 3})$}}%
%
\special{pn 20}%
\special{pa 260 698}%
\special{pa 1920 1410}%
\special{fp}%
%
\special{pn 20}%
\special{sh 0}%
\special{ar 4742 1154 46 50  0.0000000 6.2831853}%
\put(47.1500,-14.1900){\makebox(0,0)[lb]{$c_i=c_{i+3}$}}%
%
\special{pn 20}%
\special{ar 3128 1302 46 50  0.0000000 6.2831853}%
%
\special{pn 20}%
\special{pa 1910 660}%
\special{pa 260 1410}%
\special{fp}%
%
\special{pn 20}%
\special{sh 0}%
\special{ar 5482 1144 46 50  0.0000000 6.2831853}%
%
\special{pn 20}%
\special{ar 3668 1440 46 50  0.0000000 6.2831853}%
\end{picture}%
\end{center}
\caption{On the stratum of type $A_1^{\oplus 3}$} 
\label{fig:A_1^3}
\vspace{2mm}
\begin{center}
\unitlength 0.1in
\begin{picture}( 62.0000, 14.3000)(  2.6000,-19.8000)
%
\special{pn 20}%
\special{pa 4550 1770}%
\special{pa 4566 1724}%
\special{pa 4582 1676}%
\special{pa 4598 1630}%
\special{pa 4614 1584}%
\special{pa 4628 1536}%
\special{pa 4644 1490}%
\special{pa 4660 1446}%
\special{pa 4676 1400}%
\special{pa 4692 1356}%
\special{pa 4708 1312}%
\special{pa 4724 1268}%
\special{pa 4738 1226}%
\special{pa 4754 1184}%
\special{pa 4770 1142}%
\special{pa 4786 1102}%
\special{pa 4802 1062}%
\special{pa 4818 1024}%
\special{pa 4834 988}%
\special{pa 4850 952}%
\special{pa 4864 916}%
\special{pa 4880 884}%
\special{pa 4896 852}%
\special{pa 4912 820}%
\special{pa 4928 792}%
\special{pa 4944 764}%
\special{pa 4960 738}%
\special{pa 4976 714}%
\special{pa 4992 692}%
\special{pa 5008 670}%
\special{pa 5024 652}%
\special{pa 5040 634}%
\special{pa 5056 620}%
\special{pa 5072 606}%
\special{pa 5088 596}%
\special{pa 5104 588}%
\special{pa 5120 582}%
\special{pa 5136 578}%
\special{pa 5152 576}%
\special{pa 5168 578}%
\special{pa 5184 580}%
\special{pa 5200 586}%
\special{pa 5216 596}%
\special{pa 5232 606}%
\special{pa 5248 620}%
\special{pa 5264 638}%
\special{pa 5282 656}%
\special{pa 5298 678}%
\special{pa 5314 700}%
\special{pa 5330 726}%
\special{pa 5346 752}%
\special{pa 5362 782}%
\special{pa 5380 812}%
\special{pa 5396 842}%
\special{pa 5412 876}%
\special{pa 5428 910}%
\special{pa 5444 944}%
\special{pa 5460 980}%
\special{pa 5478 1016}%
\special{pa 5494 1054}%
\special{pa 5510 1092}%
\special{pa 5526 1128}%
\special{pa 5542 1166}%
\special{pa 5560 1204}%
\special{pa 5576 1242}%
\special{pa 5592 1280}%
\special{pa 5608 1318}%
\special{pa 5624 1354}%
\special{pa 5640 1390}%
\special{pa 5658 1426}%
\special{pa 5674 1460}%
\special{pa 5690 1494}%
\special{pa 5706 1526}%
\special{pa 5722 1556}%
\special{pa 5738 1586}%
\special{pa 5754 1614}%
\special{pa 5770 1638}%
\special{pa 5786 1662}%
\special{pa 5802 1686}%
\special{pa 5818 1704}%
\special{pa 5834 1722}%
\special{pa 5850 1738}%
\special{pa 5866 1750}%
\special{pa 5882 1760}%
\special{pa 5898 1768}%
\special{pa 5914 1772}%
\special{pa 5930 1772}%
\special{pa 5946 1770}%
\special{pa 5962 1766}%
\special{pa 5978 1758}%
\special{pa 5994 1748}%
\special{pa 6008 1736}%
\special{pa 6024 1720}%
\special{pa 6040 1704}%
\special{pa 6054 1684}%
\special{pa 6070 1662}%
\special{pa 6084 1640}%
\special{pa 6100 1614}%
\special{pa 6114 1588}%
\special{pa 6130 1558}%
\special{pa 6144 1528}%
\special{pa 6158 1496}%
\special{pa 6172 1464}%
\special{pa 6186 1430}%
\special{pa 6200 1394}%
\special{pa 6214 1358}%
\special{pa 6228 1322}%
\special{pa 6242 1284}%
\special{pa 6254 1246}%
\special{pa 6268 1206}%
\special{pa 6280 1168}%
\special{pa 6294 1128}%
\special{pa 6306 1088}%
\special{pa 6318 1048}%
\special{pa 6330 1008}%
\special{pa 6342 968}%
\special{pa 6354 928}%
\special{pa 6364 888}%
\special{pa 6376 850}%
\special{pa 6386 812}%
\special{pa 6398 774}%
\special{pa 6408 738}%
\special{pa 6418 702}%
\special{pa 6428 666}%
\special{pa 6438 632}%
\special{pa 6446 600}%
\special{pa 6456 568}%
\special{pa 6460 550}%
\special{sp}%
\put(56.3200,-12.8500){\makebox(0,0)[lb]{$c_{j+3}$}}%
\put(54.2500,-8.7800){\makebox(0,0)[lb]{$c_j$}}%
\put(61.5200,-17.1300){\makebox(0,0)[lb]{$c_{k+3}$}}%
\put(63.4200,-11.4600){\makebox(0,0)[lb]{$c_k$}}%
%
\special{pn 20}%
\special{pa 270 828}%
\special{pa 1920 828}%
\special{fp}%
\put(5.1300,-16.2300){\makebox(0,0)[lb]{$c_{i+3}$}}%
\put(10.0100,-12.7500){\makebox(0,0)[lb]{$c_{j+3}$}}%
\put(10.1900,-7.3900){\makebox(0,0)[lb]{$c_j$}}%
\put(5.0300,-7.3900){\makebox(0,0)[lb]{$c_i$}}%
\put(14.6100,-14.8400){\makebox(0,0)[lb]{$c_{k+3}$}}%
\put(14.9700,-7.3900){\makebox(0,0)[lb]{$c_k$}}%
%
\special{pn 20}%
\special{ar 2652 838 46 50  0.0000000 6.2831853}%
\put(24.3600,-7.4900){\makebox(0,0)[lb]{$c_i=c_{i+3}$}}%
\put(29.5100,-11.2600){\makebox(0,0)[lb]{$c_j=c_{j+3}$}}%
\put(36.8200,-8.4800){\makebox(0,0)[lb]{$c_k$}}%
\put(36.4600,-16.3300){\makebox(0,0)[lb]{$c_{k+3}$}}%
%
\special{pn 20}%
\special{sh 1.000}%
\special{ar 6370 908 46 50  0.0000000 6.2831853}%
%
\special{pn 20}%
\special{sh 1.000}%
\special{ar 6162 1484 46 50  0.0000000 6.2831853}%
%
\special{pn 20}%
\special{sh 1.000}%
\special{ar 5592 1286 46 50  0.0000000 6.2831853}%
%
\special{pn 20}%
\special{sh 1.000}%
\special{ar 5430 918 46 50  0.0000000 6.2831853}%
%
\special{pn 20}%
\special{sh 1.000}%
\special{ar 558 828 46 50  0.0000000 6.2831853}%
%
\special{pn 20}%
\special{sh 1.000}%
\special{ar 1064 828 46 50  0.0000000 6.2831853}%
%
\special{pn 20}%
\special{sh 1.000}%
\special{ar 1542 828 46 50  0.0000000 6.2831853}%
%
\special{pn 20}%
\special{sh 1.000}%
\special{ar 558 1386 46 50  0.0000000 6.2831853}%
%
\special{pn 20}%
\special{sh 1.000}%
\special{ar 1064 1048 46 50  0.0000000 6.2831853}%
%
\special{pn 20}%
\special{sh 1.000}%
\special{ar 1542 1246 46 50  0.0000000 6.2831853}%
%
\special{pn 20}%
\special{sh 1.000}%
\special{ar 3728 958 46 50  0.0000000 6.2831853}%
%
\special{pn 20}%
\special{sh 1.000}%
\special{ar 3728 1386 46 50  0.0000000 6.2831853}%
\put(4.0000,-21.5000){\makebox(0,0)[lb]{$B_1(A_1^{\oplus 2})$ or $B_2(A_1^{\oplus 2})$}}%
\put(30.7700,-21.5000){\makebox(0,0)[lb]{$B_3(A_1^{\oplus 2})$}}%
\put(51.8100,-21.5000){\makebox(0,0)[lb]{$B_4(A_1^{\oplus 2})$}}%
%
\special{pn 20}%
\special{pa 260 700}%
\special{pa 1922 1416}%
\special{fp}%
%
\special{pn 20}%
\special{sh 0}%
\special{ar 4776 1126 46 50  0.0000000 6.2831853}%
\put(47.4800,-13.9500){\makebox(0,0)[lb]{$c_i=c_{i+3}$}}%
%
\special{pn 20}%
\special{ar 3158 898 46 50  0.0000000 6.2831853}%
\end{picture}%
\end{center}
\caption{On the stratum of type $A_1^{\oplus 2}$} 
\label{fig:A_1^2}
\vspace{2mm}
\begin{center}
\unitlength 0.1in
\begin{picture}( 61.3800, 14.3000)(  2.7000,-19.8000)
%
\special{pn 20}%
\special{pa 4504 1770}%
\special{pa 4520 1724}%
\special{pa 4534 1676}%
\special{pa 4550 1630}%
\special{pa 4566 1584}%
\special{pa 4582 1536}%
\special{pa 4598 1490}%
\special{pa 4614 1444}%
\special{pa 4628 1400}%
\special{pa 4644 1356}%
\special{pa 4660 1310}%
\special{pa 4676 1268}%
\special{pa 4692 1224}%
\special{pa 4708 1182}%
\special{pa 4724 1142}%
\special{pa 4738 1102}%
\special{pa 4754 1062}%
\special{pa 4770 1024}%
\special{pa 4786 986}%
\special{pa 4802 950}%
\special{pa 4818 916}%
\special{pa 4834 882}%
\special{pa 4850 850}%
\special{pa 4866 820}%
\special{pa 4882 790}%
\special{pa 4896 764}%
\special{pa 4912 738}%
\special{pa 4928 714}%
\special{pa 4944 690}%
\special{pa 4960 670}%
\special{pa 4976 652}%
\special{pa 4992 634}%
\special{pa 5008 620}%
\special{pa 5024 606}%
\special{pa 5040 596}%
\special{pa 5056 588}%
\special{pa 5072 582}%
\special{pa 5088 578}%
\special{pa 5104 576}%
\special{pa 5120 578}%
\special{pa 5136 580}%
\special{pa 5152 586}%
\special{pa 5168 596}%
\special{pa 5184 608}%
\special{pa 5200 622}%
\special{pa 5218 638}%
\special{pa 5234 656}%
\special{pa 5250 678}%
\special{pa 5266 702}%
\special{pa 5282 726}%
\special{pa 5298 754}%
\special{pa 5314 782}%
\special{pa 5330 812}%
\special{pa 5348 844}%
\special{pa 5364 876}%
\special{pa 5380 910}%
\special{pa 5396 946}%
\special{pa 5412 982}%
\special{pa 5428 1018}%
\special{pa 5446 1056}%
\special{pa 5462 1092}%
\special{pa 5478 1130}%
\special{pa 5494 1168}%
\special{pa 5510 1206}%
\special{pa 5526 1244}%
\special{pa 5542 1282}%
\special{pa 5560 1320}%
\special{pa 5576 1356}%
\special{pa 5592 1392}%
\special{pa 5608 1428}%
\special{pa 5624 1462}%
\special{pa 5640 1494}%
\special{pa 5656 1526}%
\special{pa 5672 1558}%
\special{pa 5688 1586}%
\special{pa 5704 1614}%
\special{pa 5722 1640}%
\special{pa 5738 1664}%
\special{pa 5754 1686}%
\special{pa 5770 1706}%
\special{pa 5786 1724}%
\special{pa 5802 1738}%
\special{pa 5818 1750}%
\special{pa 5832 1760}%
\special{pa 5848 1768}%
\special{pa 5864 1772}%
\special{pa 5880 1772}%
\special{pa 5896 1770}%
\special{pa 5912 1766}%
\special{pa 5928 1758}%
\special{pa 5942 1748}%
\special{pa 5958 1734}%
\special{pa 5974 1720}%
\special{pa 5988 1702}%
\special{pa 6004 1682}%
\special{pa 6020 1662}%
\special{pa 6034 1638}%
\special{pa 6050 1612}%
\special{pa 6064 1586}%
\special{pa 6078 1556}%
\special{pa 6094 1526}%
\special{pa 6108 1494}%
\special{pa 6122 1462}%
\special{pa 6136 1428}%
\special{pa 6150 1392}%
\special{pa 6164 1356}%
\special{pa 6178 1320}%
\special{pa 6190 1282}%
\special{pa 6204 1242}%
\special{pa 6216 1204}%
\special{pa 6230 1164}%
\special{pa 6242 1124}%
\special{pa 6254 1084}%
\special{pa 6268 1044}%
\special{pa 6280 1004}%
\special{pa 6290 964}%
\special{pa 6302 926}%
\special{pa 6314 886}%
\special{pa 6324 848}%
\special{pa 6336 808}%
\special{pa 6346 772}%
\special{pa 6356 734}%
\special{pa 6366 698}%
\special{pa 6376 664}%
\special{pa 6386 630}%
\special{pa 6396 598}%
\special{pa 6404 566}%
\special{pa 6408 550}%
\special{sp}%
\put(48.3200,-11.2600){\makebox(0,0)[lb]{$c_1$}}%
\put(46.5200,-15.9300){\makebox(0,0)[lb]{$c_4$}}%
\put(54.3600,-15.4400){\makebox(0,0)[lb]{$c_5$}}%
\put(52.8200,-11.7600){\makebox(0,0)[lb]{$c_2$}}%
\put(61.0200,-17.1300){\makebox(0,0)[lb]{$c_6$}}%
\put(62.9000,-11.4600){\makebox(0,0)[lb]{$c_3$}}%
%
\special{pn 20}%
\special{pa 270 828}%
\special{pa 1916 828}%
\special{fp}%
\put(5.0400,-17.9200){\makebox(0,0)[lb]{$c_{i+3}$}}%
\put(10.0900,-16.0300){\makebox(0,0)[lb]{$c_{j+3}$}}%
\put(10.1800,-7.3900){\makebox(0,0)[lb]{$c_j$}}%
\put(4.9600,-7.3900){\makebox(0,0)[lb]{$c_i$}}%
\put(14.7600,-17.1300){\makebox(0,0)[lb]{$c_{k+3}$}}%
\put(14.8600,-7.3900){\makebox(0,0)[lb]{$c_k$}}%
%
\special{pn 20}%
\special{ar 2666 858 46 50  0.0000000 6.2831853}%
\put(24.4900,-7.6900){\makebox(0,0)[lb]{$c_i=c_{i+3}$}}%
\put(31.0700,-7.9800){\makebox(0,0)[lb]{$c_j$}}%
\put(36.9100,-8.6800){\makebox(0,0)[lb]{$c_k$}}%
\put(30.9800,-17.3300){\makebox(0,0)[lb]{$c_{j+3}$}}%
\put(36.4700,-16.3300){\makebox(0,0)[lb]{$c_{k+3}$}}%
%
\special{pn 20}%
\special{sh 1.000}%
\special{ar 6318 908 44 50  0.0000000 6.2831853}%
%
\special{pn 20}%
\special{sh 1.000}%
\special{ar 6110 1484 46 50  0.0000000 6.2831853}%
%
\special{pn 20}%
\special{sh 1.000}%
\special{ar 5544 1286 46 50  0.0000000 6.2831853}%
%
\special{pn 20}%
\special{sh 1.000}%
\special{ar 5382 918 44 50  0.0000000 6.2831853}%
%
\special{pn 20}%
\special{sh 1.000}%
\special{ar 4842 858 46 50  0.0000000 6.2831853}%
%
\special{pn 20}%
\special{sh 1.000}%
\special{ar 4644 1336 46 50  0.0000000 6.2831853}%
%
\special{pn 20}%
\special{sh 1.000}%
\special{ar 558 828 46 50  0.0000000 6.2831853}%
%
\special{pn 20}%
\special{sh 1.000}%
\special{ar 1064 828 46 50  0.0000000 6.2831853}%
%
\special{pn 20}%
\special{sh 1.000}%
\special{ar 1540 828 44 50  0.0000000 6.2831853}%
%
\special{pn 20}%
\special{sh 1.000}%
\special{ar 558 1534 46 50  0.0000000 6.2831853}%
%
\special{pn 20}%
\special{sh 1.000}%
\special{ar 1064 1356 46 50  0.0000000 6.2831853}%
%
\special{pn 20}%
\special{sh 1.000}%
\special{ar 1530 1494 46 50  0.0000000 6.2831853}%
%
\special{pn 20}%
\special{sh 1.000}%
\special{ar 3170 878 46 50  0.0000000 6.2831853}%
%
\special{pn 20}%
\special{sh 1.000}%
\special{ar 3738 978 44 50  0.0000000 6.2831853}%
%
\special{pn 20}%
\special{sh 1.000}%
\special{ar 3180 1474 46 50  0.0000000 6.2831853}%
%
\special{pn 20}%
\special{sh 1.000}%
\special{ar 3738 1406 44 50  0.0000000 6.2831853}%
\put(4.0000,-21.5000){\makebox(0,0)[lb]{$B_1(A_1)$ or $B_2(A_1)$}}%
\put(30.8900,-21.5000){\makebox(0,0)[lb]{$B_3(A_1)$}}%
\put(52.0600,-21.5000){\makebox(0,0)[lb]{$B_4(A_1)$}}%
\end{picture}%
\end{center}
\caption{On the stratum of type $A_1$} 
\label{fig:A_1}
\vspace{2mm}
\end{figure}
\begin{figure}[t]
\begin{center}
\unitlength 0.1in
\begin{picture}( 47.8000, 12.0000)(  2.5000,-13.9000)
%
\special{pn 20}%
\special{pa 3120 440}%
\special{pa 5030 440}%
\special{fp}%
\put(40.1000,-9.5000){\makebox(0,0)[lb]{$c_{j+3}$}}%
\put(40.1000,-3.6000){\makebox(0,0)[lb]{$c_j$}}%
\put(31.4000,-6.7000){\makebox(0,0)[lb]{$c_i=c_{i+3}$}}%
\put(45.3000,-11.9000){\makebox(0,0)[lb]{$c_{k+3}$}}%
\put(45.3000,-3.6000){\makebox(0,0)[lb]{$c_k$}}%
%
\special{pn 20}%
\special{sh 1.000}%
\special{ar 4060 450 50 50  0.0000000 6.2831853}%
%
\special{pn 20}%
\special{sh 1.000}%
\special{ar 4590 450 50 50  0.0000000 6.2831853}%
%
\special{pn 20}%
\special{sh 1.000}%
\special{ar 4060 710 50 50  0.0000000 6.2831853}%
%
\special{pn 20}%
\special{sh 1.000}%
\special{ar 4590 970 50 50  0.0000000 6.2831853}%
\put(37.2000,-15.5000){\makebox(0,0)[lb]{$B_2(A_3)$}}%
%
\special{pn 20}%
\special{pa 250 440}%
\special{pa 2160 440}%
\special{fp}%
\put(4.1000,-3.6000){\makebox(0,0)[lb]{$c_i=c_{i+3}$}}%
\put(16.9000,-10.3000){\makebox(0,0)[lb]{$c_{k+3}$}}%
\put(16.7000,-3.6000){\makebox(0,0)[lb]{$c_k$}}%
%
\special{pn 20}%
\special{sh 0}%
\special{ar 640 450 50 50  0.0000000 6.2831853}%
%
\special{pn 20}%
\special{sh 0}%
\special{ar 1200 450 50 50  0.0000000 6.2831853}%
%
\special{pn 20}%
\special{sh 1.000}%
\special{ar 1730 450 50 50  0.0000000 6.2831853}%
%
\special{pn 20}%
\special{sh 1.000}%
\special{ar 1730 810 50 50  0.0000000 6.2831853}%
\put(8.6000,-15.6000){\makebox(0,0)[lb]{$B_1(A_3)$}}%
%
\special{pn 20}%
\special{pa 3110 260}%
\special{pa 5030 1180}%
\special{fp}%
\put(9.7000,-6.5000){\makebox(0,0)[lb]{$c_j=c_{j+3}$}}%
%
\special{pn 20}%
\special{sh 0}%
\special{ar 3500 440 50 50  0.0000000 6.2831853}%
\end{picture}%
\end{center}
\caption{On the stratum of type $A_3$} 
\label{fig:A_3}
\begin{center}
\unitlength 0.1in
\begin{picture}( 49.0000, 13.5000)(  3.6000,-14.6000)
%
\special{pn 20}%
\special{pa 3010 360}%
\special{pa 5260 360}%
\special{fp}%
\put(36.7000,-8.2000){\makebox(0,0)[lb]{$c_{i+3}$}}%
\put(42.5000,-10.0000){\makebox(0,0)[lb]{$c_{j+3}$}}%
\put(42.4000,-2.8000){\makebox(0,0)[lb]{$c_j$}}%
\put(36.6000,-2.8000){\makebox(0,0)[lb]{$c_i$}}%
\put(47.6000,-12.1000){\makebox(0,0)[lb]{$c_{k+3}$}}%
\put(47.6000,-2.8000){\makebox(0,0)[lb]{$c_k$}}%
%
\special{pn 20}%
\special{sh 1.000}%
\special{ar 3730 370 50 50  0.0000000 6.2831853}%
%
\special{pn 20}%
\special{sh 1.000}%
\special{ar 4290 370 50 50  0.0000000 6.2831853}%
%
\special{pn 20}%
\special{sh 1.000}%
\special{ar 4820 370 50 50  0.0000000 6.2831853}%
%
\special{pn 20}%
\special{sh 1.000}%
\special{ar 3730 540 50 50  0.0000000 6.2831853}%
%
\special{pn 20}%
\special{sh 1.000}%
\special{ar 4290 760 50 50  0.0000000 6.2831853}%
%
\special{pn 20}%
\special{sh 1.000}%
\special{ar 4820 990 50 50  0.0000000 6.2831853}%
\put(39.5000,-16.2000){\makebox(0,0)[lb]{$B_2(A_2)$}}%
%
\special{pn 20}%
\special{pa 3000 240}%
\special{pa 5260 1170}%
\special{fp}%
%
\special{pn 20}%
\special{pa 360 360}%
\special{pa 2270 360}%
\special{fp}%
\put(12.9000,-11.1000){\makebox(0,0)[lb]{$c_{j+3}$}}%
\put(12.6000,-2.8000){\makebox(0,0)[lb]{$c_j$}}%
\put(5.2000,-2.8000){\makebox(0,0)[lb]{$c_i=c_{i+3}$}}%
\put(18.0000,-9.5000){\makebox(0,0)[lb]{$c_{k+3}$}}%
\put(17.8000,-2.8000){\makebox(0,0)[lb]{$c_k$}}%
%
\special{pn 20}%
\special{sh 0}%
\special{ar 750 370 50 50  0.0000000 6.2831853}%
%
\special{pn 20}%
\special{sh 1.000}%
\special{ar 1310 370 50 50  0.0000000 6.2831853}%
%
\special{pn 20}%
\special{sh 1.000}%
\special{ar 1840 370 50 50  0.0000000 6.2831853}%
%
\special{pn 20}%
\special{sh 1.000}%
\special{ar 1320 870 50 50  0.0000000 6.2831853}%
%
\special{pn 20}%
\special{sh 1.000}%
\special{ar 1840 730 50 50  0.0000000 6.2831853}%
\put(9.7000,-16.3000){\makebox(0,0)[lb]{$B_1(A_2)$}}%
\end{picture}%
\end{center}
\caption{On the stratum of type $A_2$} 
\label{fig:A_2}
\end{figure} 
\begin{theorem} \label{thm:param}
Given a parameter $\k \in \K$, put $\th=\mathrm{rh}(\k)$. 
If $\k \in \K(*)$ with $* \in \mathcal{I}/S_4$, then the 
affine cubic surface $\Sol(\th)$ has simple 
singularities of abstract Dynkin type $*$. 
\end{theorem} 
{\it Proof}. 
First, notice that $\ol{\Sol}(\th)$ has all its singularities 
within its affine part $\Sol(\th)$, since $\ol{\Sol}(\th)$ is 
smooth around the tritangent lines at infinity $L$. 
Thanks to Proposition \ref{prop:type}, in order to establish 
the theorem, it suffices to prove $\K(*)=\beta^{-1}(B(*))$ 
for every $* \in \mathcal{I}/S_4$. 
To this end we use the $W(D_4^{(1)})$-action on the parameter 
space $B$, where the action $w_i : b \mapsto b'$ is given by 
\[
b'_0=
\left\{ 
\begin{array}{cl}
{\-}b_0^{-1} \quad & (i=0), \\[1mm] 
b_0 b_i \quad & (i=1,2,3), \\[1mm]
-b_0 b_4 \quad & (i=4),
\end{array}\right. 
\qquad b_j'= 
\left\{ 
\begin{array}{cl}
b_0 b_j \quad & (i=0, \, j \in \{1,2,3,4\}), \\[1mm] 
b_i^{-1} \quad & (i, j \in \{1,2,3,4\}, \, i = j), \\[1mm]
b_j \quad & (i, j \in \{1,2,3,4\}, \, i \neq j). 
\end{array}\right. 
\]
\par
Let $\k \in \K(*)$. 
From the definition of $\K(*)$ there exist 
$w \in W(D_4^{(1)})$ and $\k' \in H(*)$ such that 
$\k = w(\k')$, where $H(*)$ is given by (\ref{eqn:F4H}). 
An inspection of Table \ref{tab:type} shows that 
$\beta(\k') \in B(*)$; actually $B(*)$ has been defined 
so that this is the case. 
It implies that 
\[
\beta(\k) = \beta(w(\k')) = w(\beta(\k')) 
\in w(B(*)) = B(*), 
\]
where the last equality follows from the 
$W(D_4^{(1)})$-invariance of the set $B(*)$. 
Thus we have $\k \in \beta^{-1}(B(*))$ and hence 
the inclusion $\K(*) \subset \beta^{-1}(B(*))$. 
The proof of the reverse inclusion $\beta^{-1}(B(*)) \subset 
\K(*)$ relies on two claims, which are presented in the 
next two paragraphs. 
\par
The first claim asserts that for any $b \in B(*)$ there exist 
$b' \in \beta(H(*))$ and $w \in W(D_4^{(1)})$ such that 
$b=w(b')$. 
We see this for the case $* = A_1$ dividing it into several 
subcases. 
In the subcase $b \in B_4(A_1)$ the claim is true because any 
$b \in B_4(A_1)$ is either of the forms $(\pm 1,b_1,b_2,b_3,b_4)$, 
where the negative pattern $(-1,b_1,b_2,b_3,b_4)$ is recast 
to a positive pattern $(1,b_1^{-1},b_2^{-1},b_3^{-1},b_4^{-1})$ 
by the action $w_1 w_2 w_3 w_4$, and any positive pattern 
certainly belongs to $\beta(H(*))$. 
The remaining subcases $b \in B_1(A_1)$, $B_2(A_1)$, $B_3(A_1)$ 
can be treated along the same line with the help of actions 
\[
\begin{array}{l}
(\pm 1, b_i, b_j, b_k, b_l) \overset{w_i}{\longmapsto} 
(\pm b_i, b_i^{-1}, b_j, b_k, b_l) \overset{w_j}{\longmapsto} 
(\pm b_i b_j, b_i^{-1}, b_j^{-1}, b_k, b_l), \\[2mm]
(\pm b_i,b_i^{-1},b_j,b_k,b_l) \overset{w_0}{\longmapsto} 
(\pm b_i^{-1}, \pm 1, \pm b_i b_j, \pm b_i b_k, \pm b_i b_l). 
\end{array}
\]
In a similar manner the first claim is valid for every 
abstract Dynkin type $* \in \mathcal{I}/S_4$. 
\par
The second claim is that $\beta^{-1}(\beta(\K(*))) = \K(*)$. 
Indeed any element of $\beta^{-1}(\beta(\K (*)))$ is expressed 
as $(\k_0+n_0,\k_1+ 2 n_1,\k_2+2 n_2,\k_3+2 n_3,\k_4+ 2n_4)$ 
for some $(\k_0,\k_1,\k_2,\k_3,\k_4) \in \K(*)$ and some 
$n_i \in \Z$ such that $n_0+n_1+n_2+n_3+n_4=0$. 
Since there exists an action 
\[
w_i w_0 (w_j w_k w_l w_0)^2 : (\k_0,\k_i,\k_j,\k_k,\k_l) \mapsto 
(\k_0+1,\k_i-2,\k_j,\k_k,\k_l) 
\quad (\{i,j,k,l\}=\{1,2,3,4\}), 
\]
and $\K(*)$ is $W(D_4^{(1)})$-invariant, one has 
$\beta^{-1}(\beta(\K(*))) \subset \K(*)$ and so 
$\beta^{-1}(\beta(\K(*))) = \K(*)$. 
\par
Now let $\k \in \beta^{-1}(B(*))$, that is, $\beta(\k) \in B(*)$. 
From the first claim there exist $b' \in \beta(H(*))$ and 
$w \in W(D_4^{(1)})$ such that 
$\beta(\k)=w(b') \in w(\beta(H(*))) = \beta(w(H(*))) \subset 
\beta(\K(*))$. 
Therefore the second claim implies $\k \in \K(*)$ and thus 
$\beta^{-1}(B(*)) \subset \K(*)$. 
The proof is complete. \hfill $\Box$ 
\section{Dynamics on Cubic Surfaces} \label{sec:dynamics} 
In this section we discuss the polynomial automorphisms on the 
cubic surface mentioned in the Introduction 
and in particular investigate their dynamical properties. 
They were intently studied in \cite{Cantat,CL,Iwasaki1,IU1} 
and the expositions below are largely based on \cite{IU1}. 
\par
Since $\Sol(\th)$ has the structure of a $(2,2,2)$-surface, namely, 
the defining function $f(x,\th)$ of $\Sol(\th)$ is quadratic in each 
variable $x_i$, the line through $x \in \Sol(\th)$ parallel to the 
$x_i$-axis passes through a unique second point $x' \in \Sol(\th)$. 
This defines an involutive automorphism 
\[
\sigma_i : \Sol(\th) \to \Sol(\th), \quad x \mapsto x'. 
\]
The surface $\Sol(\th)$ admits a natural complex area-form called 
the {\sl Poincar\'e residue}: 
\[
\omega(\th) := \dfrac{dx_1 \wedge dx_2 \wedge dx_3}{d_x f(x,\th)} 
\qquad \mbox{restricted to} \quad \Sol(\th). 
\]
It pulls back to the natural $2$-form $\omega_z(\k)$ on $\M_z(\k)$ 
via the Riemann-Hilbert correspondence (\ref{eqn:RH2}) 
(see \cite{Iwasaki0}). 
Moreover it is almost preserved by the map $\sigma_i$, namely, 
it is sent to its negative 
\begin{equation} \label{eqn:pre}
\sigma_i^* \omega(\th) = - \omega(\th) \qquad (i = 1,2,3). 
\end{equation}
\par
Let $G$ be the group generated by three involutions 
$\sigma_1$, $\sigma_2$, $\sigma_3$, and $G(2)$ its index-two 
subgroup generated by three elements $\si_1 \si_2$, $\si_2 \si_3$, 
$\si_3 \si_1$. 
It is known that $G$ is of finite index in the group of all 
polynomial automorphisms of $\Sol(\th)$ (see \cite{CL,El-Huti}). 
Although our main interest is in an element of $G(2)$, we spend 
a short while working with a general element of $G$. 
Each element $\si \in G$ extends to a birational map on 
$\ol{\Sol}(\th)$ and it in turn lifts to a one on 
$\wt{\Sol}(\th)$. 
For the biregular map $\si : \Sol(\th) \carl$ the induced 
birational maps on $\ol{\Sol}(\th)$ and $\wt{\Sol}(\th)$ are 
repersented by the same symbol $\si$. 
Note that the birational map $\si : \wt{\Sol}(\th) \carl$ restricts 
to an automorphism of $\wt{\Sol}(\th) \setminus L$, still denoted by 
$\si$. 
The area-form $\omega(\th)$ induces a meromorphic 
$2$-from $\wt{\omega}(\theta)$ on $\wt{\Sol}(\th)$, 
whose pole divisor is the sum $L_1+L_2+L_3$ of the 
three lines at infinity. 
\par
Recall that the concept of a non-elementary loop in 
$\pi_1(Z,z)$ was defined in Definition \ref{def:elementary}. 
Its counterpart in the group $G$ is defined in the following 
manner, whose relation with the original concept will be 
discussed in Section \ref{sec:MT}. 
\begin{definition} \label{def:elementary2}
An AS element $\si \in G$ is said to be {\sl elementary} if 
$\si = (\si_i \si_j)^n$ for some $\{i,j,k\} = \{1,2,3\}$ and 
$n \in \Z$; otherwise, $\si$ is said to be {\sl non-elementary}. 
\end{definition}
We describe how an element $\si \in G$ acts on the subspace 
$V \subset H^2(\widetilde{\Sol}(\th),\C)$ spanned by 
$L_1$, $L_2$, $L_3$. 
This was done in \cite{IU1} when $\Sol(\th)$ is smooth and 
it carries over when $\Sol(\th)$ is singular. 
\begin{enumerate}
\item For each $i \in \{1, 2, 3\}$, $\si_i$ blows down the 
line $L_i$ to the point $p_i$, which is the unique indeterminacy 
point of $\si_i$. 
Moreover $\si_i$ restricts to an automorphism of $L_j$ that 
exchanges $p_i$ and $p_k$, where $\{i,j,k\}=\{1, 2, 3\}$, 
as in Figure \ref{fig:cubic2} (see \cite[Lemma 3]{IU1}). 
Thus the endomorphisms $\si_1^*$, $\si_2^*$, 
$\si_3^* :H^2(\widetilde{\Sol}(\th), \Z) \carl$ map the 
subspace $V$ into itself and their restrictions to $V$ are 
represented by the matrices 
\begin{equation} \label{eqn:matrices}
s_1 = \begin{pmatrix} 0 & 1 & 1 \\ 0 & 1 & 0 \\ 0 & 0 & 1 
\end{pmatrix}, \qquad 
s_2 = \begin{pmatrix} 1 & 0 & 0 \\ 1 & 0 & 1 \\ 0 & 0 & 1 
\end{pmatrix}, \qquad 
s_3 = \begin{pmatrix} 1 & 0 & 0 \\ 0 & 1 & 0 \\ 1 & 1 & 0 
\end{pmatrix}, 
\end{equation}
respectively, relative to the basis $L_1$, $L_2$, $L_3$ 
(see \cite[Lemma 10]{IU1}). 
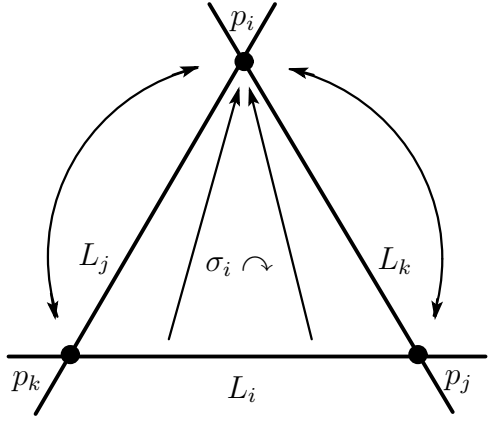
\begin{figure}[t] 
\begin{center}
\unitlength 0.1in
\begin{picture}( 24.7000, 21.7000)(  6.0000,-28.7000)
\put(17.3000,-28.1000){\makebox(0,0)[lb]{$L_i$}}%
\put(9.6000,-21.4000){\makebox(0,0)[lb]{$L_j$}}%
\put(25.1000,-21.4000){\makebox(0,0)[lb]{$L_k$}}%
\put(28.6000,-27.6000){\makebox(0,0)[lb]{$p_j$}}%
\put(17.5000,-8.7000){\makebox(0,0)[lb]{$p_i$}}%
\put(6.2000,-27.5000){\makebox(0,0)[lb]{$p_k$}}%
\put(16.2000,-21.5000){\makebox(0,0)[lb]{$\si_i \car$}}%
%
\special{pn 20}%
\special{sh 1.000}%
\special{ar 1820 1030 40 40  0.0000000 6.2831853}%
%
\special{pn 20}%
\special{sh 1.000}%
\special{ar 2720 2560 40 40  0.0000000 6.2831853}%
%
\special{pn 20}%
\special{sh 1.000}%
\special{ar 920 2560 40 40  0.0000000 6.2831853}%
%
\special{pn 20}%
\special{pa 1640 730}%
\special{pa 2900 2860}%
\special{fp}%
%
\special{pn 20}%
\special{pa 1980 730}%
\special{pa 740 2870}%
\special{fp}%
%
\special{pn 20}%
\special{pa 600 2570}%
\special{pa 3070 2570}%
\special{fp}%
%
\special{pn 13}%
\special{ar 1830 2060 1016 1016  4.9362935 6.2831853}%
\special{ar 1830 2060 1016 1016  0.0000000 0.2970642}%
%
\special{pn 13}%
\special{ar 1820 2050 1016 1016  2.8289604 4.4722180}%
%
\special{pn 13}%
\special{pa 2820 2270}%
\special{pa 2810 2360}%
\special{fp}%
\special{sh 1}%
\special{pa 2810 2360}%
\special{pa 2838 2296}%
\special{pa 2816 2308}%
\special{pa 2798 2292}%
\special{pa 2810 2360}%
\special{fp}%
%
\special{pn 13}%
\special{pa 820 2240}%
\special{pa 840 2360}%
\special{fp}%
\special{sh 1}%
\special{pa 840 2360}%
\special{pa 850 2292}%
\special{pa 832 2308}%
\special{pa 810 2298}%
\special{pa 840 2360}%
\special{fp}%
%
\special{pn 13}%
\special{pa 2180 1100}%
\special{pa 2070 1070}%
\special{fp}%
\special{sh 1}%
\special{pa 2070 1070}%
\special{pa 2130 1108}%
\special{pa 2122 1084}%
\special{pa 2140 1068}%
\special{pa 2070 1070}%
\special{fp}%
%
\special{pn 13}%
\special{pa 1450 1100}%
\special{pa 1570 1060}%
\special{fp}%
\special{sh 1}%
\special{pa 1570 1060}%
\special{pa 1500 1062}%
\special{pa 1520 1078}%
\special{pa 1514 1100}%
\special{pa 1570 1060}%
\special{fp}%
%
\special{pn 13}%
\special{pa 1430 2480}%
\special{pa 1790 1180}%
\special{fp}%
\special{sh 1}%
\special{pa 1790 1180}%
\special{pa 1754 1240}%
\special{pa 1776 1232}%
\special{pa 1792 1250}%
\special{pa 1790 1180}%
\special{fp}%
%
\special{pn 13}%
\special{pa 2170 2480}%
\special{pa 1850 1180}%
\special{fp}%
\special{sh 1}%
\special{pa 1850 1180}%
\special{pa 1848 1250}%
\special{pa 1864 1232}%
\special{pa 1886 1240}%
\special{pa 1850 1180}%
\special{fp}%
\end{picture}%
\end{center}
\caption{The birational map $\si_i$ restricted to $L$}
\label{fig:cubic2} 
\end{figure}
\item Given any element $\si \in G$ other than the unit element, 
we can write 
\begin{equation} \label{eqn:reduced}
\si = \si_{i_1} \si_{i_2} \cdots \si_{i_m},  
\end{equation} 
for some $m \in \N$ and some $m$-tuple of indices 
$(i_1,\dots,i_m) \in \{1,2,3\}^m$. 
Here we may assume that every neighboring indices 
$i_{\nu}$ and $i_{\nu+1}$ are distinct, 
because $\si_i$ is an involution. 
The expression (\ref{eqn:reduced}) with this condition is 
unique; it is called the reduced expression of $\si$. 
Thus $G$ is isomorphic to the universal Coxeter group of 
rank $3$ (see \cite[Theorem 4]{IU1}). 
From what we mentioned in item (1), the following hold for 
the expression (\ref{eqn:reduced}): 
\begin{equation} \label{eqn:EI}
\si^n I(\si^{-1}) = \{p_{i_1}\} \quad (n \ge 0), \qquad 
\si^{-n}I(\si^{-1}) = \bigcup_{\nu=1}^m L_{i_{\nu}} 
\quad (n \ge 1), 
\end{equation}
where $I(\si)$ stands for the indeterminacy set of the birational 
map $\si : \wt{\Sol}(\th) \circlearrowleft$. 
\item Let $f$ and $g$ be bimeromorphic maps on a compact 
K\"{a}hler surface $X$. 
For the induced actions $f^*$ and $g^*$ on $H^{1,1}(X)$ 
the composition rule $(f \circ g)^* = g^* \circ f^*$ is not 
always true. 
This is true if and only if $g$ blows down no curve into 
a point of $I(f)$. 
It follows from item (1) that every neighboring pair 
$\si_{i_{\nu}}$ and $\si_{i_{\nu+1}}$ satisfies this 
condition so that 
\begin{equation} \label{eqn:comp}
\si^{*} = \si_{i_m}^{*} \cdots \si_{i_2}^{*} 
\si_{i_1}^{*} : H^{1,1}(\wt{\Sol}(\th))=H^2(\wt{\Sol}(\th),\C) 
\circlearrowleft, 
\end{equation}
provided that (\ref{eqn:reduced}) is a reduced expression 
(see \cite[Lemma 8]{IU1}). 
\item It is easily seen from formula (\ref{eqn:EI}) that an 
element $\si \in G$ is AS if and only if the initial index 
$i_1$ and the terminal index $i_m$ are distinct in 
expression $(\ref{eqn:reduced})$. 
Moreover any element is conjugate to some AS element 
in $G$ (see \cite[Lemma 12 and page 324]{IU1}). 
In what follows we may and shall assume that $\si$ is AS, 
namely, that $i_1 \neq i_m$. 
\item By formula (\ref{eqn:comp}) and the matrix 
representations (\ref{eqn:matrices}), the eigenvalues of 
$\si^*|_{V}$ are $0$ and the two roots of the quadratic 
equation 
\begin{equation} \label{eqn:quadratic}
\l^2 - \a(\si) \l + (-1)^m = 0, 
\end{equation}
where $\a(\si)$ is the trace of the matrix 
$s := s_{i_m} \cdots s_{i_2} s_{i_1}$, 
which takes an even positive integer. 
Moreover $\a(\si) > 2$ if and only if 
$\si$ is non-elementary (see \cite[Lemma 13]{IU1}). 
\item Assume that $\si$ is non-elementary. 
Then for any $n \in \N$ the $n$-th iterate $\si^n$ has 
exactly two fixed points $p_{i_1}$ and $p_{i_m}$ on $L$ 
in the sense of Definition \ref{def:fixed}, where 
$p_{i_1} \in X_0^{\circ}(\si^n)$ and 
$p_{i_m} \in X_0^{\circ}(\si^{-n})$ are superattracting 
fixed points of $\si^n$ and $\si^{-n}$ respectively in 
the usual sense. 
So their indices are $\nu_{p_{i_1}}(\si^n) = 
\nu_{p_{i_m}}(\si^n) = 1$ (see \cite[Lemma 15]{IU1}). 
Moreover one has $\{i_1, i_2, \dots, i_m\}=\{1,2,3\}$ 
and thus $\si^{-n} I(\si^{-1}) = L$ from formula 
(\ref{eqn:EI}). 
In particular $\si : \wt{\Sol}(\th) \carl$ contracts 
$L$ to the point $p_{i_1}$. 
\end{enumerate}
\par
From now on we focus our attention on an even element 
$\si \in G(2)$. 
It follows from (\ref{eqn:pre}) that the birational map 
$\si : \wt{\Sol}(\th) \carl$ preserves the meromorphic 
$2$-form $\wt{\omega}(\th)$, that is, 
$\si^* \wt{\omega}(\th)=\wt{\omega}(\th)$. 
\begin{theorem} \label{thm:FDD} 
Assume that $\si \in G(2)$ is an AS element. 
Then the first dynamical degree $\lambda(\si)$ is a quadratic 
unit that appears as the largest root of the quadratic 
equation $(\ref{eqn:quadratic})$ with $m$ even. 
Moreover $\lambda(\si)$ is strictly greater than $1$ 
if and only if $\si$ is non-elementary. 
\end{theorem} 
{\it Proof}. 
Consider the action $\si^* : H^2(\wt{\Sol}(\th), \C) \carl$. 
We know that $\si^*$ preserves the subspace $V$. 
It also preserves its orthogonal complement $V^{\perp}$. 
Indeed, for any $v \in V^{\perp}$ and $v' \in V$ 
one has $(\si^* v, v')=(v, (\si^{-1})^* v')=0$, since 
$(\si^{-1})^*$ is the adjoint of $\si^*$ relative to 
the intersection form and preserves $V$. 
This shows that $\si^*$ preserves $V^{\perp}$. 
Now we claim that the operator $\si^*|_{V^{\perp}}$ is unitary. 
Indeed, a corollary to the push-pull formula 
(see \cite[Corollary 3.4]{DF}) yields 
\[
(\si^*v_1, \si^*v_2) = (v_1,v_2) + Q(v_1,v_2) \qquad 
(v_1, v_2 \in H^2(\widetilde{\Sol}(\th), \C)), 
\]
where $Q(v_1,v_2)$ is a nonnegative Hermitian form 
that can be expressed as 
\[
Q(v_1,v_2) = \sum_{i=1}^{3} k_i 
\cdot (v_1,L_i) \cdot (v_2,L_i),
\]
with some positive integers $k_1$, $k_2$, $k_3 \in \N$. 
Thus if $v_1$ and $v_2$ are in $V^{\perp}$ then 
$Q(v_1,v_2)$ vanishes and so $\si^*|_{V^{\perp}}$ 
preserves the intersection form on $V^{\perp}$. 
Recall that the vectors $L_4$, $L_5$, $L_6$, $L_7$ in 
(\ref{eqn:bas}) form a basis of $V^{\perp}$, whose 
intersection relations are known to be $(L_i, L_j) = 
-2 \delta_{ij}$ from formula (\ref{eqn:in1}), 
where $\delta_{ij}$ is Kronecker's delta.  
Thus the intersection form on $V^{\perp}$ is negative 
definite. 
Since $\si^*|_{V^{\perp}}$ preserves a negative 
definite Hermitian form, it must be unitary. 
In particular all of its eigenvalues are of modulus 
$1$. 
\par
On the other hand, since $\si$ is assumed to be AS, the 
eigenvalues of $\si|_{V}$ consist of $0$ and the two 
roots of quadratic equation (\ref{eqn:quadratic}) 
with $m$ even. 
These three numbers and the four numbers of modulus $1$ 
in the last paragraph constitute all the seven eigenvalues 
of $\si^* : H^2(\wt{\Sol}(\th), \C) \carl$. 
Note that equation (\ref{eqn:quadratic}) has a real root 
$\ge 1$ because $\a(\si) \ge 2$. 
Thus the first dynamical degree $\lambda(\si)$, which is 
the spectral radius of $\si^* : H^{1,1}(\wt{\Sol}(\th)) 
= H^2(\wt{\Sol}(\th), \C) \carl$, is given by the 
largest root of equation (\ref{eqn:quadratic}). 
Moreover $\lambda(\si) > 1$ if and only if $\a(\si) > 2$, 
which is the case precisely when $\si$ is non-elementary. 
The proof is complete. \hfill $\Box$ \par\medskip 
We now apply the construction in Remark \ref{rem:nonwander} 
to $X = \wt{\Sol}(\th)$ and $f = \si$. 
The birational map $\si : \wt{\Sol}(\th) \carl$ restricts 
to an automorphism $\si : \wt{\Sol}(\th)_{\si} \carl$, where 
$\wt{\Sol}(\th)_{\si}$ designates the space $X_f$ of 
definition (\ref{eqn:invset}) adapted in the present setting. 
This space can be identified in the following. 
\begin{lemma} \label{lem:SS}
For any non-elementary AS element $\si \in G(2)$, 
we have $\wt{\Sol}(\th)_{\si} = \wt{\Sol}(\th) \setminus L$. 
\end{lemma}
{\it Proof}. This readily follows from what we have mentioned 
in item (6). \hfill $\Box$ 
\begin{theorem} \label{thm:nonwander} 
If $\si \in G(2)$ is a non-elementary AS element, then the 
nonwandering set $\Omega(\si)$ of the birational map 
$\si : \wt{\Sol}(\th) \carl$ is compact in 
$\wt{\Sol}(\th) \setminus L$ and the trajectory of 
each point $x \in \wt{\Sol}(\th) \setminus \Omega(\si)$ 
tends to infinity $L$ under the iterations of $\si$. 
\end{theorem}
{\it Proof}. 
Put $\si=\si_{i_1} \si_{i_2} \cdots \si_{i_{m}}$ as in 
(\ref{eqn:reduced}). 
Since $\si : \wt{\Sol}(\th) \carl$ contracts $L$ into the 
superattracting fixed point $p_{i_1} \in L$, there exists a 
neighborhood $U$ of $L$ in $\wt{\Sol}(\th)$ such that each 
point of $U$ is attracted to $p_{i_1}$ under the iterations 
of $\si$. 
Hence the nonwandering set $\Omega(\si) \subset 
\wt{\Sol}(\th)_{\si} = \wt{\Sol}(\th) \setminus L$ of $
\si : \wt{\Sol}(\th) \carl$ is contained in $\wt{\Sol}(\th) 
\setminus U$ and thus compact in $\wt{\Sol}(\th) \setminus L$ 
because $\Omega(\si)$ is closed. \hfill $\Box$
\begin{theorem} \label{thm:chaos} 
Assume that $\si \in G(2)$ is a non-elementary AS element. 
Then the birational map $\si : \wt{\Sol}(\th) \carl$ admits 
a $\si$-invariant Borel probability measure $\nu_{\si}$ with 
support in $\Omega(\ga)$ that satisfies the conditions 
in Definition $\ref{def:chaos}$. 
Moreover, 
\begin{enumerate}
\item the measure-theoretic entropy $h_{\nu_{\si}}(\si)$ 
and the topological entropy $h_{\mathrm{top}}(\si)$ are 
expressed as
\begin{equation} \label{eqn:top}
h_{\nu_{\si}}(\si) = h_{\mathrm{top}}(\si) = \log \lambda(\si), 
\end{equation}
\item the measure $\nu_{\si}$ puts no mass on any algebraic 
curve on $\wt{\Sol}(\th)$, 
\item there exists a set $\mathcal{P}_n(\si) 
\subset \mathrm{supp}\, \nu_{\si}$ of saddle periodic points 
of period $n$ such that 
\[
\# \mathcal{P}_n(\si) \sim \lambda(\si)^n, \qquad 
\dfrac{1}{\lambda(\si)^n} \sum_{p \in \mathcal{P}_n(\si)} 
\delta_p \to \nu_{\si}, \qquad \mbox{as} \quad n \to \infty. 
\]
\end{enumerate}
\end{theorem}
{\it Proof}. 
Put $\si=\si_{i_1} \si_{i_2} \cdots \si_{i_{m}}$ as 
in (\ref{eqn:reduced}). 
From Theorem \ref{thm:FDD} the first dynamical degree 
$\lambda(\si)$ is strictly greater than $1$. 
Since $I(\si)=\{p_{i_m}\}$ and 
$\si^N(I(\si^{-1})) =\{p_{i_1}\}$ for any $N \ge 0$,  
\[
\sum_{N = 0}^{\infty} \l(\si)^{-N} 
\log \, \mathrm{dist}(\si^N I(\si^{-1}), I(\si))
= \log \, \mathrm{dist}(p_{i_1}, p_{i_n}) 
\sum_{N = 0}^{\infty} \l(\si)^{-N}
> - \infty, 
\]
and hence condition (\ref{eqn:AS3}) is satisfied. 
Therefore Theorem \ref{thm:BD} implies that there exists a 
$\si$-invariant Borel probability measure $\nu_{\si}$ 
that satisfies all conditions of the theorem. 
\hfill $\Box$ \par\medskip 
We turn our attention to the second topic, that is, estimating 
the number of isolated periodic points of the birational map 
$\si : \wt{\Sol}(\th) \carl$ for a given element $\si \in G(2)$. 
Let $\mathrm{Per}_n^i(\si \setminus L)$ denote the set of 
all isolated periodic points of period $n$ that 
lie in $\wt{\Sol}(\th) \setminus L$. 
\begin{theorem} \label{thm:per} 
Let $\si \in G(2)$ be a non-elementary AS element. 
Then any irreducible periodic curve of the map 
$\si : \wt{\Sol}(\th) \carl$ must lie in the exceptional 
set $\E(\th)$ of the minimal resolution  $\pi : \wt{\Sol}(\th) 
\to \ol{\Sol}(\th)$ in $(\ref{eqn:des})$. 
Moreoveor for every $n \in \N$ the set 
$\mathrm{Per}_n^i(\si \setminus L)$ is finite and its 
cardinality counted with multiplicity is estimated as 
\begin{equation} \label{eqn:esti}
|\# \mathrm{Per}_n^i(\si \setminus L) - \lambda(\si)^n| 
\le O(1) \qquad \mbox{as} \quad n \to \infty. 
\end{equation}
\end{theorem} 
{\it Proof}. 
We begin with the assertion about periodic curves. 
Let $\si = \si_{i_1} \si_{i_2} \cdots \si_{i_m}$ be the 
reduced expression as in (\ref{eqn:reduced}). 
In view of item (6), since $\si \in G(2)$ is non-elementary, 
for any $n \in \N$ the $n$-th iterate $\si^n$ contracts $L$ 
to the superattracting fixed point $p_{i_1} \in L$ of $\si$. 
We claim that any periodic point of the map 
$\si : \ol{\Sol}(\th) \carl$ must be isolated. 
Indeed, assume the contrary that $\si$ admits a periodic 
curve $C \subset \ol{\Sol}(\th)$. 
Since $\Sol(\th)$ is affine, the compact curve $C$ 
must intersect the lines at infinity $L$. 
If $n \in \N$ is the primitive period of $C$, then the 
fixed curve $C$ of $\si^n$ must meet $L$ in $p_{i_1}$, 
because $L$ is contracted into the superattracting 
fixed point $p_{i_1}$ by $\si^{n}$. 
This contradicts the fact that $p_{i_1}$ is an isolated 
fixed point of $\si^{n}$. 
Therefore any irreducible periodic curve of 
$\si : \wt{\Sol}(\th) \carl$ must be contracted 
to a point by the map (\ref{eqn:des}), that is, 
it must lie in the exceptional set $\E(\th)$ 
of the resolution (\ref{eqn:des}). 
\par
The above argument shows in particular that no 
irreducible component of the pole divisor 
$(\wt{\omega}(\theta))_{\infty} = L_1 + L_2 + L_3$ 
is a periodic curve of $\si$. 
Thus Theorem \ref{thm:esti} gives an estimate
\[
|\# \mathrm{Per}_n^{i}(\si) - \lambda(\si)^n| \le 
O(1), 
\]
since $\wt{\Sol}(\th)$ is not birationally equivalent 
to any Abelian surface. 
On the other hand the number 
$\# \mathrm{Per}_n^i(\si \setminus L)$ is obtained from 
$\# \mathrm{Per}_n^{i}(\si)$ by subtracting the sum of 
local indices at the isolated fixed points on $L$ 
for the map $\si^n$. 
In view of item (6) those fixed points are just $p_{i_1}$ 
and $p_{i_m}$, each having local index $1$. 
Therefore the number $\# \mathrm{Per}_n^i(\si \setminus L)$ 
is given by 
\begin{equation} \label{eqn:RP}
\# \mathrm{Per}_n^i(\si \setminus L) = 
\# \mathrm{Per}_n^{i}(\si)-2, 
\end{equation}
so that the estimate (\ref{eqn:esti}) is established by 
combining all these observations. 
\hfill $\Box$ \par\medskip 
By virtue of Corollary \ref{cor:esti}, Theorems \ref{thm:per} 
has the following. 
\begin{corollary} \label{cor:asymp} 
Assume that $\si \in G(2)$ is a non-elementary AS element. 
Then we have 
\[
\# \mathrm{Per}_n^i(\si \setminus L) \sim 
\# \mathrm{HPer}_n(\si) 
\sim \lambda(\si)^n \qquad \mbox{as} \quad n \to \infty, 
\]
so that asymptotically almost all points in 
$\mathrm{Per}_n^i(\si \setminus L)$ belong to 
$\mathrm{HPer}_n(\si)$. 
\end{corollary}
\begin{remark} \label{rem:dynkin}
From Theorem \ref{thm:per} any periodic curve of a 
non-elementary map $\si : \wt{\Sol}(\th) \carl$ must 
be an irreducible component of $\E(\th)$. 
On the other hand $\E(\th)$ is the exceptional set 
of a minimal resolution of simple singularities on 
$\ol{\Sol}(\th)$. 
Thus no two of the irreducible components of $\E(\th)$ 
are tangent and no three of them meet in a single point. 
Therefore every iterate of the map $\si$ satisfies 
Assumption \ref{ass:fix} and it is why this assumption 
was made in Section \ref{sec:Pre}. 
\end{remark}
\section{Proofs of the Main Theorems} \label{sec:MT} 
In this section we prove our main theorems 
combining all the previous discussions. 
As is mentioned in the Introduction, the Riemann-Hilbert 
correspondence (\ref{eqn:RH2}) recasts the monodromy map on 
the moduli space $\M_z(\k)$ to a biregular map on the cubic 
surface $\Sol(\th)$, where the latter map was studied in 
Section \ref{sec:dynamics}. 
After a brief review of how these two maps are related, we 
establish our main theorems by translating the results on 
$\Sol(\th)$ back to $\M_z(\k)$. 
\par
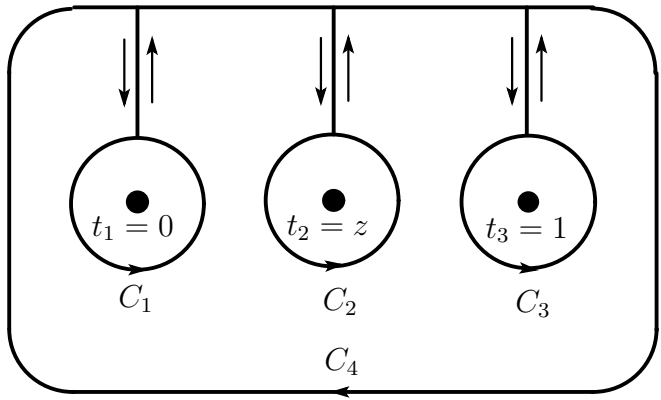
\begin{figure}[t]
\begin{center}
\unitlength 0.1in
\begin{picture}( 33.5200, 20.1300)(  2.1000,-26.1300)
%
\special{pn 20}%
\special{ar 538 936 328 328  3.1177876 4.6856097}%
%
\special{pn 20}%
\special{pa 210 936}%
\special{pa 210 2280}%
\special{fp}%
%
\special{pn 20}%
\special{ar 538 2280 328 328  1.5469913 3.1415927}%
%
\special{pn 20}%
\special{ar 1882 1608 344 344  0.0000000 6.2831853}%
%
\special{pn 20}%
\special{ar 2898 1616 346 346  0.0000000 6.2831853}%
%
\special{pn 20}%
\special{pa 546 600}%
\special{pa 3234 600}%
\special{fp}%
%
\special{pn 20}%
\special{pa 538 2608}%
\special{pa 3234 2608}%
\special{fp}%
%
\special{pn 20}%
\special{ar 3234 2280 328 328  6.2831853 6.2831853}%
\special{ar 3234 2280 328 328  0.0000000 1.5951817}%
%
\special{pn 20}%
\special{pa 874 600}%
\special{pa 874 1272}%
\special{fp}%
%
\special{pn 20}%
\special{pa 1890 600}%
\special{pa 1890 1264}%
\special{fp}%
%
\special{pn 20}%
\special{pa 2890 600}%
\special{pa 2890 1264}%
\special{fp}%
%
\special{pn 20}%
\special{pa 3562 944}%
\special{pa 3562 2288}%
\special{fp}%
%
\special{pn 20}%
\special{ar 3226 936 330 330  4.7123890 6.2831853}%
\special{ar 3226 936 330 330  0.0000000 0.0244600}%
%
\special{pn 20}%
\special{ar 874 1626 344 344  0.0000000 6.2831853}%
%
\special{pn 13}%
\special{sh 1.000}%
\special{ar 2898 1616 50 50  0.0000000 6.2831853}%
%
\special{pn 20}%
\special{sh 1.000}%
\special{ar 1890 1608 50 50  0.0000000 6.2831853}%
%
\special{pn 20}%
\special{sh 1.000}%
\special{ar 874 1616 50 50  0.0000000 6.2831853}%
%
\special{pn 13}%
\special{pa 806 768}%
\special{pa 806 1104}%
\special{fp}%
\special{sh 1}%
\special{pa 806 1104}%
\special{pa 826 1038}%
\special{pa 806 1052}%
\special{pa 786 1038}%
\special{pa 806 1104}%
\special{fp}%
%
\special{pn 13}%
\special{pa 1824 760}%
\special{pa 1824 1096}%
\special{fp}%
\special{sh 1}%
\special{pa 1824 1096}%
\special{pa 1844 1030}%
\special{pa 1824 1044}%
\special{pa 1804 1030}%
\special{pa 1824 1096}%
\special{fp}%
%
\special{pn 13}%
\special{pa 2814 760}%
\special{pa 2814 1096}%
\special{fp}%
\special{sh 1}%
\special{pa 2814 1096}%
\special{pa 2834 1030}%
\special{pa 2814 1044}%
\special{pa 2794 1030}%
\special{pa 2814 1096}%
\special{fp}%
%
\special{pn 13}%
\special{pa 950 1096}%
\special{pa 950 768}%
\special{fp}%
\special{sh 1}%
\special{pa 950 768}%
\special{pa 930 836}%
\special{pa 950 822}%
\special{pa 970 836}%
\special{pa 950 768}%
\special{fp}%
%
\special{pn 13}%
\special{pa 1966 1088}%
\special{pa 1966 760}%
\special{fp}%
\special{sh 1}%
\special{pa 1966 760}%
\special{pa 1946 828}%
\special{pa 1966 814}%
\special{pa 1986 828}%
\special{pa 1966 760}%
\special{fp}%
%
\special{pn 13}%
\special{pa 2966 1088}%
\special{pa 2966 760}%
\special{fp}%
\special{sh 1}%
\special{pa 2966 760}%
\special{pa 2946 828}%
\special{pa 2966 814}%
\special{pa 2986 828}%
\special{pa 2966 760}%
\special{fp}%
%
\special{pn 13}%
\special{pa 832 1970}%
\special{pa 908 1970}%
\special{fp}%
\special{sh 1}%
\special{pa 908 1970}%
\special{pa 840 1950}%
\special{pa 854 1970}%
\special{pa 840 1990}%
\special{pa 908 1970}%
\special{fp}%
%
\special{pn 13}%
\special{pa 1840 1952}%
\special{pa 1924 1944}%
\special{fp}%
\special{sh 1}%
\special{pa 1924 1944}%
\special{pa 1856 1930}%
\special{pa 1872 1950}%
\special{pa 1860 1970}%
\special{pa 1924 1944}%
\special{fp}%
%
\special{pn 13}%
\special{pa 2864 1962}%
\special{pa 2932 1970}%
\special{fp}%
\special{sh 1}%
\special{pa 2932 1970}%
\special{pa 2868 1942}%
\special{pa 2880 1964}%
\special{pa 2864 1982}%
\special{pa 2932 1970}%
\special{fp}%
%
\special{pn 13}%
\special{pa 1966 2608}%
\special{pa 1890 2608}%
\special{fp}%
\special{sh 1}%
\special{pa 1890 2608}%
\special{pa 1958 2628}%
\special{pa 1944 2608}%
\special{pa 1958 2588}%
\special{pa 1890 2608}%
\special{fp}%
\put(8.6500,-21.2000){\makebox(0,0){$C_1$}}%
\put(19.2400,-21.3700){\makebox(0,0){$C_2$}}%
\put(29.2300,-21.4600){\makebox(0,0){$C_3$}}%
\put(19.3200,-24.4800){\makebox(0,0){$C_4$}}%
\put(8.5000,-17.5000){\makebox(0,0){$t_1=0$}}%
\put(18.6000,-17.5000){\makebox(0,0){$t_2=z$}}%
\put(28.9000,-17.6000){\makebox(0,0){$t_3=1$}}%
\end{picture}%
\end{center}
\caption{Four loops in $\P^1 \setminus \{0,z,1,\infty\}$; 
the fourth point $t_4=\infty$ is outside $C_4$, invisible.} 
\label{fig:loop2}
\end{figure}

 
Given a pair $(z,a) \in Z \times A$, let $\mathcal{R}_z(a)$ 
be the moduli space of Jordan equivalence classes of 
representations $\rho : \pi_1(\P^1 \setminus \{t_1,t_2,t_3,t_4\},*) 
\to SL_2(\C)$ such that $\Tr\, \rho(C_i) = a_i$ for 
$i \in \{1,2,3,4\}$, where $t_1 = 0$, $t_2 = z$, $t_3 = 1$, 
$t_4 = \infty$ and $C_i$ is a loop surrounding $t_i$ once 
anti-clockwise as in Figure \ref{fig:loop2}. 
The space $\mathcal{R}_z(a)$ is called a relative 
$SL_2(\C)$-character variety of the quadruply punctured 
sphere $\P^1 \setminus \{t_1,t_2,t_3,t_4\}$. 
Then the Riemann-Hilbert correspondence is the map 
\begin{equation} \label{eqn:RH}
\RH_{z,\k} : \M_z(\k) \to \mathcal{R}_z(a), \quad 
Q \mapsto \rho 
\end{equation} 
sending each stable parabolic connnection $Q$ to the Jordan 
equivalence class $\rho$ of its monodromy representation, 
where $\k \mapsto a$ is the composition of the maps $\a$ and 
$\b$ in (\ref{eqn:KBATh}). 
On the other hand, there exists an isomorphism of affine 
algebraic surfaces 
\[
\mathcal{R}_z(a) \to \Sol(\th), \quad \rho \mapsto x 
= (x_1,x_2,x_3), 
\]
where $x_i = \Tr\,\rho(C_jC_k)$ for $\{i,j,k\} = \{1,2,3\}$ 
and $a \mapsto \th$ is the map $\varphi$ in (\ref{eqn:KBATh}). 
It enables us to identify the character variety 
$\mathcal{R}_z(a)$ with the affine cubic surface $\Sol(\th)$, 
so that the Riemann-Hilbert correspondence (\ref{eqn:RH}) 
can be reformulated as the map in (\ref{eqn:RH2}) with 
$\th = \rh(\k)$ in (\ref{eqn:rhp}). 
\begin{theorem}[\cite{IIS1,IIS2,IIS3}]  \label{thm:SolRHP} 
Given any $\k \in \K_I$, put $\th = \rh(\k) \in \Th$. 
If $I = \emptyset$ then the surface $\Sol(\th)$ is smooth 
and the Riemann-Hilbert correspondence $(\ref{eqn:RH2})$ is 
a biholomorphism. 
Otherwise, $\Sol(\th)$ has simple singularities of Dynkin 
type $D_I$ and the Riemann-Hilbert correspondence 
$(\ref{eqn:RH2})$ is a proper surjective map that is an 
analytic minimal resolution of the singular surface 
$\Sol(\th)$. 
\end{theorem}
\par
On the other hand, we have an algebraic minimal resolution 
$\pi : \wt{\Sol}(\th) \setminus L \to \Sol(\th)$ as the 
restriction of (\ref{eqn:des}) to the affine surface 
$\Sol(\th)$. 
Since the minimal resolution is unique up to isomorphisms, 
the Riemann-Hilbert correspondence (\ref{eqn:RH2}) lifts 
to a biholomorphism 
\begin{equation} \label{eqn:LRH}
\wt{\RH}_{z,\k} : \M_z(\k) \to \wt{\Sol}(\th) \setminus L 
\end{equation}
sending the exceptional set $\E_z(\k)$ of $\M_z(\k)$ to the 
exceptional set $\E(\th)$ of $\wt{\Sol}(\th) \setminus L$. 
It is known that for each $i \in \Z/3\Z$ the monodromy map 
$\ga_{i*} : \M_z(\k) \carl$ along the basic loop 
$\ga_i \in \pi_1(Z,z)$ is strictly conjugate to the automorphism 
$\si_i \si_{i+1} : \wt{\Sol}(\th) \setminus L \carl$ via the 
lifted Riemann-Hilbert correspondence (\ref{eqn:LRH}) 
(see \cite{IIS1}). 
Moreover there exists an isomorphism of groups 
\begin{equation} \label{eqn:isom} 
\Phi : \pi_1(Z,z) \to G(2) 
\end{equation} 
sending $\ga_i$ to $\si_i \si_{i+1}$ for each $i \in \Z/3\Z$ 
(see \cite{IU2}). 
Thus the monodromy map $\ga_* : \M_z(\k) \circlearrowleft$ 
along a general loop $\ga \in \pi_1(Z,z)$ is strictly conjugate 
to the automorphism 
$\si := \Phi(\ga) : \wt{\Sol}(\th) \setminus L \carl$. 
\par
Any loop $\ga \in \pi_1(Z,z)$ can be written as a word in the 
alphabet $\ga_1^{\pm1}$, $\ga_2^{\pm1}$, $\ga_3^{\pm1}$. 
Such a word of minimal length is called a reduced expression 
of $\ga$ and this minimal length is by definition 
the {\sl length} of $\ga$. 
A loop $\ga$ is said to be {\sl minimal} if it is of minimal 
length in the conjugacy class of $\ga$. 
It suffices to consider minimal loops only, because 
conjugate loops induce conjugate monodromy maps 
which are dynamically the same. 
Any minimal loop $\ga \in \pi_1(Z,z)$ is sent to an AS 
element $\si \in G(2)$ by the isomorphism (\ref{eqn:isom}) 
and vice versa. 
Moreover $\ga$ is non-elementary in the sense of Definition 
\ref{def:elementary} if and only if the AS element $\si$ 
is non-elementary in the sense of Definition \ref{def:elementary2}. 
The first dynamical degree of $\ga$ is that of the 
birational map $\si : \wt{\Sol}(\th) \carl$, namely, 
$\l(\ga) := \l(\si)$. 
Thus Theorem \ref{thm:FDD} implies that $\l(\ga) >1$ 
if and only if $\ga$ is non-elementary. 
\par
We are now in a position to establish our main theorems. 
\par\medskip\noindent
{\it Proofs of Theorems $\ref{thm:main1}$ and $\ref{thm:main2}$}. 
We may assume that $\ga \in \pi_1(Z,z)$ is a non-elementary 
minimal loop. 
Then $\si := \Phi(\ga) \in G(2)$ is a non-elementary AS element. 
First we prove Theorem \ref{thm:main1}. 
From Theorem \ref{thm:nonwander} the nonwandering set 
$\Omega(\si)$ of $\si : \wt{\Sol}(\th) \carl$ are compact 
in $\wt{\Sol}(\th) \setminus L$ and the trajectory of 
any point $x \in \wt{\Sol}(\th) \setminus \Omega(\si)$ 
tends to infinity $L$ under the iterations of $\si$. 
Since the lifted Riemann-Hilbert correspondence (\ref{eqn:LRH}) is 
proper, the nonwandering set $\Omega(\ga)$ of $\ga_* : \M_z(\k) \carl$ 
is also compact in $\M_z(\k)$ as the inverse image 
of $\Omega(\si)$ by the map (\ref{eqn:LRH}). 
For the same reason the trajectory of any point $Q \in \M_z(\k) 
\setminus \Omega(\ga)$ tends to infinity $\mathcal{Y}_{z}(\k)$ 
under the iterations of $\ga_*$. 
Let $\nu_{\si}$ be the $\si$-invariant Borel probability measure 
on $\wt{\Sol}(\th)$ mentioned in Theorem \ref{thm:chaos}, 
which has support in $\Omega(\si) \subset \wt{\Sol}(\th) \setminus L$. 
It pulls back to a $\ga_*$-invariant probability measure 
$\mu_{\ga} := \wt{\RH}_{z,\k}^*(\nu_{\si})$ on $\M_z(\k)$ 
through the lifted Riemann-Hilbert correspondence (\ref{eqn:LRH}). 
The measure $\mu_{\ga}$ has support in $\Omega(\ga)$ 
and all the assertions in Theorem \ref{thm:main1} follow from 
those in Theorem \ref{thm:chaos}. 
In particular formula (\ref{eqn:ent}) comes from 
formula (\ref{eqn:top}). 
We proceed to the proof of Theorem \ref{thm:main2}. 
Since the lifted Riemann-Hilbert correspondence (\ref{eqn:LRH}) 
is biholomorphic, one has 
\begin{equation} \label{eqn:RHPer}
\# \mathrm{Per}_n^i(\ga;\kappa)=\# \mathrm{Per}_n^i(\si \setminus L), 
\end{equation}
and all the assertions in Theorem \ref{thm:main2} follow from 
those in Theorem \ref{thm:per} and Corollary \ref{cor:asymp}. 
\hfill$\Box$
\begin{remark} \label{rem:elem} 
Let $\ga \in \pi_1(Z,z)$ be an elementary loop conjugate to the loop 
$\ga_i^m$ for some index $i \in \{1,2,3\}$ and some integer $m \in \Z$. 
Then the monodromy map $\ga_* : \M_z(\k) \carl$ is semi-conjugate to 
the map $(\si_i \si_{i+1})^m : \Sol(\th) \carl$ via the Riemann-Hilbert 
correspondence (\ref{eqn:RH2}). 
Notice that the latter map preserves the projection 
$\Sol(\th) \to \C$, $x =(x_1,x_2,x_3) \mapsto x_k$, where $(i,j,k)$ 
is a cyclic permutation of $(1, 2, 3)$. 
Through the map (\ref{eqn:RH2}) this projection is pulled back to an 
analytic fibration $\M_z(\k) \to \C$ which is preserved by the 
monodromy map $\ga_* : \M_z(\k) \carl$. 
\end{remark}
\section{Periodic Solutions along a Pochhammer Loop} \label{sec:Ex} 
We illustrate the power of Theorem \ref{thm:formula} by 
calculating for various $\k \in \K$ the explicit values of 
the number $\# \mathrm{Per}_n^i(\wp;\kappa)$ of isolated 
periodic solutions to $\PVI(\k)$ along a Pochhammer loop 
$\wp$  (cf. Example \ref{ex:card}). 
It is interesting that the result strongly depends on 
the value of $\k \in\K$. 
\par 
Some subspaces of the parameter spaces $B$ and $\K$ are 
introduced to facilitate an efficient 
case-by-case discussion. 
For each $n \in \N$ let $B^{(n)}$ be the subspace of 
$B$ defined by 
\[
\begin{array}{rcll}
B^{(1)} &:=& \bigl\{\, b \in B \,:\, b_0 = 1, \, 
b_1^2 = b_2^2 = b_3^2 = \pm\sqrt{-1} \, \bigr\},  & \\[3mm]
B^{(n)} &:=& \bigl\{\, b \in B \,:\, b_0 = 1, \, 
R_n(b_1,b_2,b_3) = 0 \, \bigr\} \qquad\qquad & (n \ge 2), 
\end{array}
\]
where
\[
\begin{array}{rcl}
R_n(b_1,b_2,b_3) &:=& 
\ds \prod_{\substack{1 \le m <n \\ (m,n)=1}} 
\Bigl( r(b_1,b_2,b_3) + 2 b_1^2 b_2^2 b_3^2 \, 
\cos \frac{\pi m}{n} \Bigr), \\[9mm] 
r(b_1,b_2,b_3) &:=& 1 - 3 b_1^2 b_2^2 b_3^2 + 
\bigl(b_1^2 b_2^2 b_3^2-1 \bigr) \, 
\Bigl\{ b_1^2 b_2^2 b_3^2 + \ds \sum_{i\in \Z/ 3 \Z}
 \bigl(b_i^2 - b_{i+1}^2 b_{i+2}^2 \bigr) \Bigr\}. 
\end{array}
\]
Moreover let $\K^{(n)}$ be the $W(F_4^{(1)})$-translates 
of the set $\beta^{-1}(B^{(n)})$, where $\beta : \K \to B$ 
is the map defined by (\ref{eqn:b}). 
For each abstract Dynkin type 
$* \in \mathcal{I}/S_4$, we define 
\[
\K^{(n)}(*) := \K(*) \cap \K^{(n)}, \qquad 
\K^{(\star)}(*) := \K(*) \setminus 
\bigcup_{n \ge 2} \K^{(n)}(*). 
\]
\par
Recall from formula (\ref{eqn:DDEP}) that the first 
dynamical degree of $\wp$ is given by 
$\l(\wp) = 9 + 4 \sqrt{5}$ and from Theorem \ref{thm:main2} 
that any periodic curve along $\wp$ must be a Riccati curve. 
For any $\k \in \K(*)$ the Riccati curves on $\M_z(\k)$ 
has the dual graph of abstract Dynkin type $* \in 
\mathcal{I}/S_4$. 
\begin{theorem} \label{thm:Poch}
Along a Pochhammer loop $\wp \in \pi_1(Z,z)$ the following hold: 
\begin{enumerate}
\item If $\k \in \K(A_1)$, then the unique Riccati curve on 
$\M_z(\k)$ is a periodic curve of primitive period $n \ge 1$ 
precisely when $\k \in \K^{(n)}(A_1)$. 
If moreover $\k \in \K^{(\star)}(A_1)$, then we have 
\[
\# \mathrm{Per}_n^i(\wp;\k) =
\left\{ 
\begin{array}{ll}
\lambda(\wp)^n + \lambda(\wp)^{-n} -10 \qquad & 
(\k \in \K^{(1)}(A_1)), \\[2mm] 
\lambda(\wp)^n + \lambda(\wp)^{-n} +4 \qquad & 
\bigl(\k \in \K^{(\star)}(A_1) \setminus \K^{(1)}(A_1)
\bigr). 
\end{array}\right. 
\]
\item If $\k \in \K(A_2)$, then neither of the two Riccati 
curves on $\M_z(\k)$ is a periodic curve of any period, and 
\[
\# \mathrm{Per}_n^i(\wp;\k) =
\lambda(\wp)^n + \lambda(\wp)^{-n} + 4. 
\]
\item If $\k \in \K(A_1^{\oplus 2})$, then neither of the 
two Riccati curves on $\M_z(\k)$ is a fixed curve. 
If moreover $\k \in \K^{(\star)}(A_1^{\oplus 2})$, then 
neither of them is a periodic curve of any period, and 
\[
\# \mathrm{Per}_n^i(\wp;\k) = 
\lambda(\wp)^n + \lambda(\wp)^{-n} + 4. 
\]
If $\k \in \K^{(n)}(A_1^{\oplus 2})$ with $n \ge 2$, then 
both of them are periodic curves of primitive period $n$. 
\item If $\k \in \K(A_3)$, then the Riccati curve corresponding 
to the central node of the Dynkin diagram of type $A_3$ is a fixed 
curve, but neither of the other two Riccati curves on 
$\M_z(\k)$ is a periodic curve of any pariod, and 
\[
\# \mathrm{Per}_n^i(\wp;\k) =
\lambda(\wp)^n + \lambda(\wp)^{-n} -2. 
\]
\item If $\k \in \K(A_1^{\oplus 3})$, then none of the three 
Riccati curves on $\M_z(\k)$ is a fixed curve. 
If moreover $\k \in \K^{(\star)}(A_1^{\oplus 3})$, then none 
of them is a periodic curve of any period, and 
\[
\# \mathrm{Per}_n^i(\wp;\k) =
\lambda(\wp)^n + \lambda(\wp)^{-n} + 4. 
\]
If $\k \in \K^{(n)}(A_1^{\oplus 3})$ with $n \ge 2$, then 
all of them are periodic curves of primitive period $n$. 
\item If $\k \in \K(D_4)$, then all of the four Riccati 
curves on $\M_z(\k)$ are fixed curves, and 
\[
\# \mathrm{Per}_n^i(\wp;\k) =
\lambda(\wp)^n + \lambda(\wp)^{-n} - 4. 
\]
\item If $\k \in \K(A_1^{\oplus 4})$, then none of the 
four Riccati curves on $\M_z(\k)$ is a periodic curve of any 
period, and
\[
\# \mathrm{Per}_n^i(\wp;\k) =
\lambda(\wp)^n + \lambda(\wp)^{-n} + 4. 
\]
\end{enumerate}
\end{theorem}
\begin{remark} \label{rem:Poch}
The number $\# \mathrm{Per}_n^i(\wp;\k)$ is yet to be determined 
for $\k \in \K^{(n)}(*)$ with $* = A_1$, $A_1^{\oplus 2}$, 
$A_1^{\oplus 3}$ and $n \ge 2$, in which cases periodic 
curves of higher periods occur and things are much subtler. 
In this paper we content ourselves with the cases allowing at 
most fixed curves. 
\end{remark}
\par
The first step toward the proof of Theorem \ref{thm:Poch} is to 
calculate the actions of the basic elements $\si_i$ on the cohomology 
group $H^2(\wt{\Sol}(\th),\C) = V \oplus V^{\perp}$, where the subspaces 
$V$ and $V^{\perp}$ are spanned by the vectors $L_1$, $L_2$, $L_3$ 
in (\ref{eqn:infline}) and by the vectors $L_4$, $L_5$, $L_6$, $L_7$ 
in (\ref{eqn:bas}) respectively. 
It is convenient to introduce another basis of $V^{\perp}$ defined by 
\[
\begin{array}{rclrcl}
L_4' &:=& E_0-E_1-E_5-E_6, \qquad & L_5' &:=& E_0-E_2-E_4-E_6, \\[2mm]
L_6' &:=& E_0-E_3-E_4-E_5, \qquad & L_7' &:=& E_0-E_1-E_2-E_3. 
\end{array}
\]
\begin{lemma} \label{lem:action}
For each $i = 1,2,3$ the action $\si_i^* : H^2(\wt{\Sol}(\th),\C) \carl$ 
preserves the subspaces $V$ and $V^{\perp}$. 
Its restriction to $V$ is represented by the matrix $s_i$ in 
$(\ref{eqn:matrices})$ relative to the basis $L_1$, $L_2$, $L_3$, 
while its restriction to ${V^{\perp}}$ has eigenvalues $\pm 1$ whose 
eigenspaces $V_{\pm 1}$ have bases in 
Table $\ref{tab:basis}$. 
\end{lemma}
\begin{table}[t]
\begin{center} 
\begin{tabular}{|c||c|c|}
\hline 
\vspace{-0.42cm} ~ & ~ & ~ \\
strata & basis of $V_{1}$ & basis of $V_{-1}$
\\[1mm] \hline \hline 
\vspace{-0.42cm} ~ & ~ & ~ \\
$A_1$, $A_2$ & $L_7$ 
& $L_4$, $L_5$, $L_6$
\\[1mm] \hline 
\vspace{-0.42cm} ~ & ~ & ~ \\
$A_1^{\oplus 2}$, $A_3$ & $L_{i+3}$, $L_7$ & 
$L_{j+3}$, $L_{k+3}$ 
\\[1mm] \hline 
\vspace{-0.42cm} ~ & ~ & ~ \\
$A_1^{\oplus 3}$ & $L_4'$, $L_5'$, 
$L_6'$ & $L_7'$
\\[1mm] \hline 
\vspace{-0.42cm} ~ & ~ & ~ \\
$D_4$, $A_1^{\oplus 4}$ & $L_4$, $L_5$, $L_6$, $L_7$ & none 
\\[1mm] \hline 
\end{tabular}
\end{center}
\caption{Bases of the eigenspaces $V_{\pm 1}$, where 
$\{i,j,k\}=\{1,2,3\}$.} 
\label{tab:basis}
\end{table}
{\it Proof.}
We only deal with the stratum $B(A_1)$ as the other strata can 
be treated in similar manners. 
Moreover it suffices to consider the case $b \in B_4(A_1)$, namely, 
the case where the six points $c_1, \dots, c_6$ in 
(\ref{eqn:indeterminacy}) lie on a conic $C$, because the 
entire $B(A_1)$ is covered by the $W(D_4^{(1)})$-translates of 
$B_4(A_1)$. 
Since $b_1b_2b_3b_4 = 1$ on $B_4(A_1)$, 
formula (\ref{eqn:indeterminacy}) reads: 
\[
\begin{array}{rclrcl}
c_1 &=& [0:1:-b_2b_3], \qquad & c_4 &=& [0:-b_2b_3:1], \\[2mm]
c_2 &=& [-b_3b_1:0:1], \qquad & c_5 &=& [1:0:-b_3b_1], \\[2mm]
c_3 &=& [1:-b_1b_2:0], \qquad & c_6 &=& [-b_1b_2:1:0], 
\end{array}
\]
and the conic passing through these points is given by 
\begin{equation} \label{eqn:conic}
C:= \Bigl\{\, [u_1:u_2:u_3] \, : \, \sum_{j \in \Z/3\Z} \{ u_j^2
+(b_j^{-1}b_{j+1}^{-1}+b_jb_{j+1})u_ju_{j+1} \}=0 \, \Bigr\}.
\end{equation}
\par
For each $i = 1,2,3$, the birational map $\phi_i := \tau 
\circ \si_i \circ \tau^{-1} : \P^2 \to \P^2$ is expressed as 
\[
\begin{array}{rcl}
\phi_1[u_1:u_2:u_3] &=& [(b_2u_2 + b_3^{-1}u_3)
(b_2^{-1}u_2 + b_3u_3) : u_1u_2 : u_1u_3], \\[2mm]
\phi_2[u_1:u_2:u_3] &=& [u_1u_2 : 
(b_3u_3 + b_1^{-1}u_1)(b_3^{-1}u_3 + b_1u_1) : u_2u_3], \\[2mm]
\phi_3[u_1:u_2:u_3] &=& [u_1u_3 : 
u_2u_3 : (b_1u_1 + b_2^{-1}u_2)(b_1^{-1}u_1 + b_2u_2)]. 
\end{array}
\]
where $\tau : \P^2 \to \P^3$ is defined in 
(\ref{eqn:birational}). 
The birational map $\phi_i$ has the indeterminacy set 
\[
I(\phi_i) = \{c_i, c_{i+3}, e_i \} \qquad \mbox{with} 
\quad 
e_1:=[1:0:0], \quad e_2:=[0:1:0], \quad e_3:=[0:0:1], 
\]
and sends the six points as 
\[
\phi_i : \left\{
\begin{array}{rclrcl}
c_i &\longleftrightarrow& C_i, \quad & 
c_{i+3} &\longleftrightarrow& C_{i+3}, \\[2mm]
c_j &\longleftrightarrow& c_{j+3}, \quad &
c_k &\longleftrightarrow& c_{k+3}, 
\end{array}
\right.
\]
with $\{i,j,k\} = \{1,2,3\}$, where $C_i$ and $C_{i+3}$ are 
lines defined by 
\[
C_i := \bigl\{ b_{i+1}u_{i+1}+b_{i+2}^{-1}u_{i+2}=0 \bigr\}, 
\quad 
C_{i+3} := \bigl\{ b_{i+1}^{-1}u_{i+1}+b_{i+2}u_{i+2}=0 \bigr\}
\qquad (i \in \Z/3\Z).
\]
The lines $C_i$ and $C_{i+3}$ pass through the points 
$c_i$ and $c_{i+3}$ respectively. 
Moreover $\phi_i$ maps a generic line in $\P^2$ to a conic 
passing through $c_i$ and $c_{i+3}$. 
Thus $\si_i^* : H^2(\wt{\Sol}(\th),\Z) \carl$ is given by 
\[
\si_i^* : \left\{
\begin{array}{rclrcl}
L &\longmapsto& 2L - E_i - E_{i+3}, \quad &  \\[2mm]
E_i &\longmapsto& L- E_i, \quad & E_{i+3} &\longmapsto& L- E_{i+3}, \\[2mm]
E_j &\longmapsto& E_{j+3}, & E_{j+3} &\longmapsto& E_j, \\[2mm]
E_k &\longmapsto& E_{k+3}, & E_{k+3} &\longmapsto& E_k. \\[2mm]
\end{array}
\right.
\]
This formula readily leads to the statement of the lemma for 
the stratum $B(A_1)$. \hfill $\Box$ \par\medskip 
The nonlinear monodromy map $\wp_* : \M_z(\k) \carl$ is strictly 
conjugate to the automorphism 
\[
\si := (\si_1 \circ \si_2 \circ \si_3)^2 : 
\wt{\Sol}(\th) \setminus L \carl 
\]
through the lifted Riemann-Hilbert correspondence 
(\ref{eqn:LRH}) (cf. \cite[Section 8]{IU2}). 
It follows from Lemma \ref{lem:action} and some calculations 
that $\si^* : V \carl$ has three simple eigenvalues $0$, 
$\lambda(\wp)$ and $\lambda(\wp)^{-1}$, while 
$\si^* : V^{\perp} \carl$ has only a quadruplicate eigenvalue $1$. 
\begin{corollary} \label{cor:Poch}
The $n$-th iterate of the birational map 
$\si : \wt{\Sol}(\th) \carl$ has Lefschetz number
\begin{equation} \label{eqn:Lef}
L(\si^n)= 6 + \lambda(\wp)^n+\lambda(\wp)^{-n} \qquad (n \in \N). 
\end{equation}
\end{corollary}
{\it Proof.} 
Since $\wt{\Sol}(\th)$ is a smooth rational surface, one has 
\[
H^q(\wt{\Sol}(\th), \C) \cong \left\{ 
\begin{array}{cl}
\C \qquad & (q = 0,4), \\[1mm] 
0  \qquad & (q = 1,3). 
\end{array}
\right. 
\]
Trivially $(\si^n)^*$ is identity on $H^0(\wt{\Sol}(\th), \C)$. 
It is also identity on $H^4(\wt{\Sol}(\th), \C)$, since $\si$ 
and so $\si^n$ are birational. 
Moreover it has three simple eigenvalues $0$, $\l(\wp)^n$, 
$\l(\wp)^{-n}$ and a quadruplicate eigenvalue $1$ on 
$H^2(\wt{\Sol}(\th),\C)$. 
So the Lefschetz number of $\si^n$ is given as in (\ref{eqn:Lef}). 
\hfill $\Box$ \par\medskip\noindent 
{\it Proof of Theorem $\ref{thm:Poch}$.} 
As in the proof of Lemma \ref{lem:action}, we consider the 
$A_1$-stratum only, since the remaining strata can be treated in 
similar manners. 
Again we may assume that $b \in B_4(A_1)$. 
Then the unique $(-2)$-curve $E$ on $\wt{\Sol}(\th)$ is 
the strict transform of the conic $C \subset \P^2$ in (\ref{eqn:conic}). 
This conic has a parametrization $g : \P^1 \ni z \mapsto 
[g_1(z):g_2(z):g_3(z)] \in C$ with 
\[
\left\{
\begin{array}{rcl}
g_1(z) &:=& b_1 b_2 (1 - b_2^2 b_3^2 ) z + 
(1-b_1^2) (1- b_2^2) + b_1^2 (1 -  b_2^2 b_3^2 ), \\[2mm]
g_2(z) &:=& -b_2 (b_2 z +b_1) (b_1 b_2 b_3^2 z+1),   \\[2mm]
g_3(z) &:=& b_2 b_3 (z +b_1 b_2) (b_1 b_2 z+1), 
\end{array}
\right.
\]
in terms of which $\si = (\phi_1 \circ \phi_2 \circ \phi_3)^2$ 
acts on $C$ and so on $E$ as the M\"obius transformation
\[
z \mapsto 
-\frac{(r(b_1,b_2,b_3)-b_1^2 b_2^4 b_3^4) z + 
b_1 b_2 b_3^2 (b_1^2 b_2^2 b_3^2 - b_1^2 b_2^2 +b_2^2 -1)}{ 
b_1 b_2^3 b_3^2 (b_1^2 b_2^2 b_3^2 - b_1^2 b_3^2 +b_3^2 -1) z + 
b_1^2 b_2^4 b_3^4 }. 
\]
Therefore $E$ is a fixed curve of $\si$ if and only if 
\[
\left\{
\begin{array}{rcl}
b_1 b_2 b_3^2 (b_1^2 b_2^2 b_3^2 - b_1^2 b_2^2 +b_2^2 -1) &=& 0, 
\\[2mm]
b_1 b_2^3 b_3^2 (b_1^2 b_2^2 b_3^2 - b_1^2 b_3^2 +b_3^2 -1) &=& 0, 
\\[2mm]
(r(b_1,b_2,b_3)-b_1^2 b_2^4 b_3^4) &=& -b_1^2 b_2^4 b_3^4. 
\end{array}
\right.
\]
A little calculation shows that the above equations hold 
if and only if $b \in B^{(1)}$. 
On the other hand $E$ is a periodic curve of primitive period 
$n \ge 2$ for $\si$ if and only if the matrix 
\[
Q =
\begin{pmatrix} 
r(b_1,b_2,b_3)-b_1^2 b_2^4 b_3^4 & 
b_1 b_2 b_3^2 (b_1^2 b_2^2 b_3^2 - b_1^2 b_2^2 +b_2^2 -1) \\[2mm] 
b_1 b_2^3 b_3^2 (b_1^2 b_2^2 b_3^2 - b_1^2 b_3^2 +b_3^2 -1) & 
b_1^2 b_2^4 b_3^4  & 
\end{pmatrix} 
\]
has eigenvalues $c \cdot \exp\bigl(\pm \frac{m}{n} \sqrt{-1} \pi \bigr)$ 
for some $c \in \C^{\times}$ and some $1 \le m < n$ such that $(m,n)=1$. 
It is easy to see that this is the case precisely when 
$b \in B^{(n)} \cap B(A_1)$, since the eigenvalues of $Q$ are the 
roots of quadratic equation 
$\l^2- r(b_1,b_2,b_3) \, \l + b_1^4 b_2^4 b_3^4=0$. 
Under the biholomorphism (\ref{eqn:LRH}) the Riccati curve on 
$\M_z(\k)$ is sent to the $(-2)$-curve $E$ on 
$\wt{\Sol}(\th) \setminus L$, so that the Riccati curve is a 
periodic curve of primitive period $n$ along $\wp$ 
if and only if $\k \in \K^{(n)}(A_1)$. 
\par
Finally we apply Theorem \ref{thm:formula} to calculate 
the exact value of $\# \mathrm{Per}_n^{i}(\wp;\kappa)$ for 
$\k \in \K^{(\star)}(A_1)$. 
First, if $\k \in \K^{(\star)}(A_1) \setminus \K^{(1)}(A_1)$, 
then the map $\si : \wt{\Sol}(\th) \carl $ has no periodic curves 
and thus $P_n(\si) = \emptyset$ for any $n \ge 1$. 
Theorem \ref{thm:formula} and Corollary \ref{cor:Poch} imply 
$\# \mathrm{Per}_n^{i}(\si) = 
L(\si^n)= 6 + \lambda(\wp)^n+\lambda(\wp)^{-n}$, 
which together with formulas (\ref{eqn:RP}) and 
(\ref{eqn:RHPer}) yields 
\[
\# \mathrm{Per}_n^{i}(\wp;\kappa) = 
\# \mathrm{Per}_n^{i}(\si \setminus L) = 
\# \mathrm{Per}_n^{i}(\si) -2 = 
\lambda(\wp)^n+\lambda(\wp)^{-n}+4. 
\]
\par
Secondly, if $\k \in \K^{(1)}(A_1)$ then the $(-2)$-curve 
$E$ is the unique fixed curve of $\si$ and $P_n(\si)=\{1\}$. 
It turns out that there are exactly six points on $E$ at which 
the local index is positive. 
Denote them by $y_1,\dots,y_6$. 
In terms of suitable local coordinates $(x_1,x_2)$ around 
$y_i$ such that $E = \{ x_1=0 \}$, the local endomorphism 
$\si_{y_i}^* : A_{y_i} \cong \C[\![ x_1,x_2]\!] \carl$ can be 
expressed as 
\[
\left\{
\begin{array}{rclcl}
\si_{y_i}^*(x_1) &=& x_1+x_1^2 \, h_1(x_1,x_2) &=& 
x_1+x_1^3 \, \wt{h}_1(x_1,x_2), \\[2mm]
\si_{y_i}^*(x_2) &=& x_2+x_1^2 \, h_2(x_1,x_2), & & 
\end{array}
\right.
\]
where $\wt{h}_1(0,0) \neq 0$ and 
$h_2(0,x_2) = 3 x_2 (1+x_2^2)^6$. 
From formulas (\ref{eqn:nup}) and (\ref{eqn:nuA}), 
\[
\nu_{E}(\si) = \nu_{(x_1)}(\si_{y_i}^*) = 2, \qquad 
\nu_{y_i}(\si)=\nu_{A_{y_i}}(\si_{y_i}^*) =  1 + 2 \cdot 1=3. 
\]
Noticing $C_1(\si)=\{y_1, \dots , y_6 \}$, 
$\mathrm{PC}_1(\si)=\{ E \}$ and $\tau_{E}=-2$, we have 
\[
\xi_1(\si) = \sum_{x \in C_1(\si)} \nu_x(\si) + 
\sum_{C \in \mathrm{PC}_1(\si)} \tau_{C} \cdot \nu_{C}(\si) = 
6 \cdot 3 + (-2) \cdot 2= 14. 
\]
Since $P_n(\si)=\{1\}$ for every $n \in \N$, 
Theorem \ref{thm:formula} and Corollary \ref{cor:Poch} imply that 
\[
\# \mathrm{Per}_n^{i}(\si) = L(\si^n) - \xi_1(\si)
= \lambda(\wp)^n+\lambda(\wp)^{-n} - 8, 
\]
which together with formulas (\ref{eqn:RP}) and (\ref{eqn:RHPer}) 
yields
\[
\# \mathrm{Per}_n^{i}(\wp;\kappa) = 
\# \mathrm{Per}_n^{i}(\si \setminus L) = 
\# \mathrm{Per}_n^{i}(\si) -2 = \l(\wp)^n + \l(\wp)^{-n} - 10. 
\]
Therefore the theorem is established for the $A_1$-stratum. 
The remaining strata can be treated in similar manners 
(we refer to \cite[Section 8]{IU2} for the $D_4$-stratum). 
\hfill $\Box$ 

\end{document}